\documentclass{amsart}
\usepackage{graphicx}
\usepackage{latexsym}
\usepackage{epsf}
\usepackage{epsfig}
\usepackage{color}
\usepackage{amsmath} 
\usepackage{amssymb}
\usepackage{amsxtra}
\usepackage{fancyhdr}
\usepackage{booktabs}
\usepackage{enumitem}
\usepackage{here}
\usepackage{amsmath,amssymb}
\usepackage{graphics,epsfig,times,epstopdf,color}
\usepackage{ae,aecompl}
\usepackage{a4wide}
\usepackage{url}
\usepackage{float}
\usepackage{caption} 
\usepackage{subcaption}
\usepackage{tikz}
\usepackage{tikz-3dplot}
\usepgflibrary{arrows}
\usepackage[utf8]{inputenc}
\usepackage[english]{babel}
\usepackage{lineno,hyperref}
\modulolinenumbers[5]

\newtheorem{remark}{Remark}

\newcommand{\OO}{\mathcal{O}}

\newcommand{\dsp}{\displaystyle}
\newcommand{\RR}{\mathbb{R}}

\newcommand{\zetabar}{\underline{\zeta}}
\newcommand{\hbarr}{\underline{h}}
\newcommand{\vbar}{\underline{v}}
\newcommand{\eps}{\varepsilon}
\newcommand{\nn}{\nonumber}
\newcommand{\red}{\textcolor{black}}
\renewcommand{\t}{\tilde}
\makeindex
\begin{document}
		\title[{\tiny Higher-ordered/extended Boussinesq system for efficient numerical simulations}]{A new class of higher-ordered/extended Boussinesq system
			for efficient numerical simulations by splitting operators}
		\date{}
		\author{Ralph Lteif}
		\address{Lebanese American University (LAU), Graduate Studies and Research (GSR) office, School of Arts and Sciences, Computer Science and Mathematics Department, Beirut, Lebanon}
		\email[Corresponding author]{ralph.lteif@lau.edu.lb}
		\author{St\'ephane Gerbi}
		\address{Laboratoire de Math\'ematiques UMR 5127 CNRS \& Universit\'e de Savoie Mont Blanc, Campus scientifique, 73376 Le Bourget du Lac Cedex, France}
       \email{stephane.gerbi@univ-smb.fr}

\begin{abstract}In this work, we numerically study the higher-ordered/extended Boussinesq system describing the propagation of water-waves over flat topography. A reformulation of the same order of precision that avoids the calculation of high order derivatives on the surface deformation is proposed. We show that this formulation enjoys an extended range of applicability while remaining stable. Moreover, a significant improvement in terms of linear dispersive properties in high frequency regime is made due to the suitable adjustment of a dispersion correction parameter. We develop a second order splitting scheme where the hyperbolic part of the system is treated with a high-order finite volume scheme and the dispersive part is treated with a finite difference approach. Numerical simulations are then performed under two main goals: validating the model and the numerical methods and assessing the potential need of such higher-order model. \red{The applicability of the proposed model and numerical method in practical problems is illustrated by a comparison with experimental data.}
\end{abstract}

\maketitle
\tableofcontents
	\section{Introduction}
	\subsection{Motivation}
	In this paper, we numerically study higher-ordered/extended Boussinesq type models. These equations describe the one-dimensional flow of the free surface of a homogeneous, immiscible fluid moving above a flat topography. They are obtained from the free surface water-waves (Euler) equations~\cite{Lannes} (see Eq. (1.3) therein) for an irrotational and incompressible fluid. 
	
	The derivation classically relies on a re-scaling of the system in order to reveal small dimensionless parameters which allow to perform asymptotic expansions of non-local operators (Dirichlet-Neumann), thus ignoring the terms whose influence is minimal. We start by introducing respectively the commonly known nonlinear and shallowness parameters:
	\begin{equation*}
		0\leq \varepsilon=\frac{a}{h_0}=\frac{\text{wave amplitude}}{\text{reference depth}}\leq 1 \; ,\end{equation*}
		
		\begin{equation*} 0\leq \sqrt{\mu}=\frac{h_0}{\lambda}=\frac{\text{reference depth}}{\text{wave-length of the wave}}<1 \; .
	\end{equation*}
The order of magnitude of these parameters makes it possible to identify the considered physical regime.

In a specific long wave regime, $\eps$ is considered of the same order as $\mu$ ($\eps \sim \mu$). In this regime, Boussinesq derived in~\cite{Boussi1871,Boussi1872} a weakly nonlinear model bearing his name. In what follows we refer to it as the ``original" or ``standard" Boussinesq system. Using the  horizontal depth-mean velocity $v$ and the free surface parametrization $\zeta$,
	the standard Boussinesq (sB) equations reads:
	\begin{equation}\label{standard-bouss}
		\left\{
		\begin{array}{lcl}
			\displaystyle\partial_t\zeta+\partial_x\big( (1+\eps\zeta )v \big)=0\vspace{1mm}\; ,\\
			\displaystyle \big(  1 - \eps \dfrac{1}{3} \partial_x^2 \big)  \partial_t v + \partial_x\zeta + \varepsilon v\partial_x v  =\mathcal{O} (\eps^2) \; .
		\end{array}
		\right.
	\end{equation}
	This model can be derived from the \emph{Green-Naghdi} (GN) equations (see~\cite{GreenNaghdi76}) by neglecting all terms of order $\mathcal{O}(\eps^2,\mu\eps,\mu^2)$. Equivalent Boussinesq systems enjoying a better mathematical structure or physical properties have been studied and derived extensively in the literature, see for instance~\cite{LPS12,MSZ12,SWX17,SX12}. 
The sB equations are restricted by containing only weak dispersion and non-linearity (only $\OO(\mu,\eps)$ terms are retained). This normally limits precise applications to a small zone moderately exterior to the surf zone. 

Significant improvement have been made in recent years to expand the application range and cover the range fully from deep water into the surf zone. Madsen et al.~\cite{MMS91,Madsen92} reached this goal by rearranging the dispersive terms (or $\OO(\mu)$ terms) in order to improve linear dispersion properties. On the other hand, Nwogu~\cite{Nwogu93} achieved the same result by redefining the dependent velocity variable. 
These models have been extensively examined for their utility in the prediction of near-shore problems (wave breaking, run-up, wave-induced circulation) as detailed in~\cite{MSS97a,MSS97b,Sorensen98}. However, these Boussinesq-type models which assume the velocity profiles to be second order polynomials in the vertical coordinate induce inaccuracies near wave breaking~\cite{Gobbi00}. 
Several mechanisms exists in order to handle wave breaking that occurs as waves approach the shore. For instance, a hybrid method consisting in suppressing the dispersive terms in breaking regions was initially suggested by Tonelli and Petti~\cite{Tonelli2009}. Another strategy consists of on an eddy viscosity approach based on the solution of a turbulent kinetic energy following early work by Nwogu~\cite{Nwogu93}. The interested reader is referred to~\cite{KAZOLEA2018} for a comparison between the hybrid and eddy viscosity strategies. Other efficient mechanisms that allow solving the wave breaking problem are worth mentioning, see for instance~\cite{Kirby2016,ESCALANTE2019} and references therein.  %

Many attempts have been made to extend Boussinesq-type models in order to offer better dispersive properties. To this end, in order to incorporate high-order dispersive and nonlinear effects in Boussinesq-type models, one should include some high-order terms in the asymptotic expansion of the velocity potential.
The first attempt to derive higher-order Boussinesq-like equations (retaining $\OO(\mu^2)$ terms) was performed by Dingemans in~\cite{Dingemans73}. Two versions of equations were given, one based on the depth-averaged velocity
and one based on the velocity at the still water level. Dingemans did not provide analyses or computations based on these equations. One can see also the review papers written by Kirby~\cite{Kirby96} and Madsen \& Schäffer~\cite{Madsen98,Madsen99} where Boussinesq-type equations of higher order in dispersion as well as in non-linearity are derived, intensified and analyzed with emphasis on linear dispersion, shoaling and nonlinear properties for large wave numbers. A fully non-linear Boussinesq-type model retaining $\OO(\mu^2)$ terms was derived for an horizontal bottom in~\cite{Gobbi00} and examined for its
ability to represent weakly nonlinear wave evolution in intermediate depth and its numerical properties of solitary wave solutions in shallow water. Algorithms for the numerical solution of the latter model are described in~\cite{Gobbi99}, where the model is applied to the study of wave shoaling
and harmonic generation in the problem of waves propagating over an isolated step. 
Including higher order terms in the model is not the only method to improve the dispersion relation. In fact, there are other strategies to improve the disperive properties, for instance, improving the degree of freedom of the velocity unknown~\cite{Castro2014,DucheneIsrawiTalhouk16,kazakova2019}. The interested reader is also referred to~\cite{FernndezNieto2018} where a hierarchy of new models is derived with a layer-wise approach incorporating non-hydrostatic effects to approximate the Euler equations. The
linear dispersion relation of these models is analyzed therein and proved to converge to the dispersion relation for the Euler equations when the number of layers goes to infinity. One can see also~\cite{Escalante2019_2} where a two-layer non-hydrostatic model with improved
dispersive properties is derived and~\cite{Escalante2021} where a numerical scheme is designed for theses models. %
 Other higher-order asymptotic shallow-water models were derived in the literature. For instance, the \emph{extended Green-Naghdi} (eGN) equations (accurate up to the order $\OO(\mu^3)$ while the full non-linearity is preserved) were firstly derived in their Hamiltonian formulation by Matsuno in~\cite{Matsuno15,Matsuno16}. More recently, Khorbatly et al. derived the eGN equations in~\cite{KZI18} by performing an asymptotic analysis of the Dirichlet–Neumann operator that originates from the formulation of the water wave problem~\cite{Lannes}. In the aforementioned papers of Khorbatly et al., the mathematical analysis is addressed which is mainly devoted to the well-posedness of the equations.
 
Note, however, that less determined efforts were made to numerically study higher-ordered asymptotic equations for the water-waves problem. Due to their extensive length, the high-order equations incorporating very high order derivatives (see for example fifth-order derivatives in~\cite{Gobbi00}) may not seem viable as a basis for a numerical model. 
 At this point, it is worth mentioning several recent advances presenting novel hyperbolic reformulation of Serre-Green-Naghdi~\cite{Bassi2020} and Boussinesq-type~\cite{Escalante2020} models, see also~\cite{Favrie2017} for the first derivation of hyperbolic reformulation of a dispersive system from variational principles in the 
flat bottom case. Those new first-order reformulations are based on a relaxed augmented system in which the divergence constraints of the velocity flow variables are coupled with the other conservation laws via an evolution equation for the depth-averaged non-hydrostatic pressures. They avoid the use of high order derivatives which are not easy to treat numerically due to the large stencil usually needed. Moreover allow to overcome the numerical difficulties and the severe time step restrictions arising from higher order terms~\cite{Busto2021}.%

Bearing these facts in mind, we develop in this paper a numerical model solving a class of higher-ordered/extended Boussinesq type models where we discuss their usefulness for practical applications. We propose a new reformulation of the model with improved linear dispersive properties and an extended range of applicability. 
 This new reformulation is suitable to the implementation of a hybrid scheme splitting the hyperbolic and dispersive parts of the equations. This strategy has been initially introduced for Boussinesq-like and Green-Naghdi equations in order to handle correctly wave breaking, see~\cite{Tonelli2009,BCLMT,LannesMarche14}. The hyperbolic part of the system is treated with a high-order finite volume scheme whereas the dispersive part is treated with a finite difference method at the same order. This splitting strategy allowed us to overcome the severe time step restriction induced due to the presence of high order derivatives by calculating the time step in the first finite-volume sub-step. The numerical investigations show that the time step restriction from the CFL condition (according to which the time step must be chosen proportional to the mesh spacing) of the finite-volume step is enough to ensure stability for the whole numerical method. %
  As a preliminary step, we treat in this paper the case of flat bottom topography, investigating the variable topography case can be based on the results of this paper and will be faced in a forthcoming work. We would like to emphasize that this paper provides sufficient grounds for treating the more complex variable topography case. 
The goal of this paper is to give a numerical assessment of the potential need of such a new formulation of higher-order model \red{that avoids the calculation of higher order derivatives
   and evaluate the effect of adding factorized high order terms to standard models.}
\subsection{Higher-ordered/extended Boussinesq equations}
Let us start first by introducing the higher-ordered \emph{extended Boussinesq} (eB) equations that naturally show up through the asymptotic approximation of the Dirichlet–Neumann operator. These equations are a straightforward extension of the sB equations~\eqref{standard-bouss}. Neglecting the terms of order $\OO(\eps^3)$ while keeping the $\OO(\eps^2)$ terms in the equations one gets the eB equations. Alternatively, one can easily recover the eB equations from the eGN equations derived in~\cite{KZI18}, by considering weak non-linearity ($\eps \sim \mu$) and dropping all terms of order $\mathcal{O}(\mu\eps^2,\mu^2\eps,\eps^3)$ therein. Thus, one can write the weakly nonlinear Boussinesq system including higher order dispersive effects as follows:
	\begin{equation}\label{ex-boussinesq}
		\left\{
		\begin{array}{lcl}
			\displaystyle\partial_t\zeta+\partial_x(hv)=0\vspace{1mm}\; ,\\
			\displaystyle (  1+\eps\mathcal{T}[\eps\zeta]+\eps^2\mathfrak{T} )\partial_t v + \partial_x\zeta+\eps  v \partial_x v   +\eps^2 \mathcal{Q}v =\mathcal{O} (\eps^3) \; ,
		\end{array}
		\right.
	\end{equation}
	where $h=1+\eps\zeta$ is the non-dimensionalised height of the fluid and denote by
	\begin{equation*}
		\mathcal{T}[\eps\zeta]w =-\frac{1}{3h}\partial_x\big(h^3\partial_xw\big), \quad\mathfrak{T} w = -\frac{1}{45}\partial_x^4w ,  \quad\mathcal{Q}v = -\frac{1}{3}\partial_x\big(vv_{xx}-v_x^2\big) \; .
	\end{equation*}
	Actually, neglecting terms of order $\OO(\eps^2)$ in~\eqref{ex-boussinesq},  the second order differential operator becomes $\mathcal{T}[\eps\zeta]w =-\frac{1}{3}\partial^2_xw$ and one can easily recover~\eqref{standard-bouss}. The eB model~\eqref{ex-boussinesq} can be found in~\cite{Madsen98} (see Eq. (3.11) therein after neglecting all terms of order $\OO(\eps\mu^4, \eps^2\mu^2)$). 
 Unfortunately, the eB model~\eqref{ex-boussinesq} seems to suffer from instabilities that turn out to be fatal for any practical use.
This is due to the positive sign in front of the elliptic fourth-order linear operator $\mathfrak{T}$ which also prevent the invertibility of the factorized operator, see~\cite[Section 3.1]{KLIG}. In fact, linearizing the eB model~\eqref{ex-boussinesq} around some rest state solution, one gets the following unstable dispersion relation:
\begin{equation}\label{rdeB}
	w^2=\dfrac{k^2}{\Big(1+\dfrac{1}{3}\eps k^2-\dfrac{1}{45}\eps^2 k^4\Big)},
\end{equation}
where $k$ is the spatial wave number and $w$ represents the time frequency. Note that one can recover the dispersion relation associated to~\eqref{standard-bouss} by neglecting $\eps^2$ terms in~\eqref{rdeB}. As expected, Madsen \& Schäffer~\cite{Madsen98} noticed an improved accuracy for small wave number values $k$ when comparing with the lower-order equations (see Figure 1 therein). However, a quick functional study shows that the denominator of the right hand side of~\eqref{rdeB} becomes negative whenever $k \gg 1$, preventing any hope concerning the well-posedness of the initial value problem. In fact, a singularity occurs in~\eqref{rdeB} for $k^2=\frac{1}{2}(15+9\sqrt{5})$, i.e. $k \approx 4.2$. Despite being a large wave number value, this singularity prove\red{s} to be inoperable for any reasonable application of the eB equations~\eqref{ex-boussinesq}, see~\cite{Madsen98}. 

To overcome this problem and in order to gain some confidence into the model validity, an enhanced set of equations of same order of precision, without instabilities, has been derived in~\cite{KLIG} by replacing the left most term of the second equation of~\eqref{ex-boussinesq}, $(  1+\eps\mathcal{T}[\eps\zeta]+\eps^2\mathfrak{T} )\partial_t v$, by $(  1+\eps\mathcal{T}[\eps\zeta]-\eps^2\mathfrak{T} )(\partial_t v) +2 \eps^2\mathfrak{T} (\partial_t v)$ and making use of a $
BBM$ \emph{trick} 
 (Benjamin-Bona-Mahony)~\cite{BBM1972} %
   represented in the following approximate equation $\partial_tv= -
	\partial_x\zeta+O(\eps)$, see~\cite{KZI18} for more details. In light of these remarks and after setting $\pm\eps^2\mathcal{T}[\eps\zeta](vv_x)$ in the second equation of~\eqref{ex-boussinesq}, one obtains the following model:
	\begin{equation}\label{eq:eGNLW1}\left\{ \begin{array}{l}
			\partial_{ t}\zeta +\partial_x\left(hv\right)\ =\  0,\\ \\
			\mathfrak{J}\left( \partial_{ t} v + \eps{v } \partial_x {v} \right) +\partial_x \zeta +\dfrac{2}{45} \eps^2 \partial_x^5\zeta+  \dfrac{2}{3} \eps^2 \partial_x((\partial_x v)^2)=  \OO(\eps^3).
		\end{array} \right. \end{equation}
	where $h=1+\eps\zeta$ and
	\begin{equation*}  \mathfrak{J}=1+\eps\mathcal{T}[\eps\zeta] -\eps^2 \mathfrak T.\end{equation*}
	The benefit of the new formulation~\eqref{eq:eGNLW1} is in replacing  $(1 +\eps\mathcal{T}[\eps\zeta]+\eps^2\mathfrak{T}) $ by a new operator $\mathfrak{J}$. This replacement induce\red{s} a fifth order derivative term on $\zeta$, namely $\dfrac{2}{45} \eps^2 \partial_x^5 \zeta$, but the invertibility of the operator $\mathfrak{J}$ is now earned. In fact, this technique modifies the dispersion relation~\eqref{rdeB} into:
	\begin{equation}\label{rdeB1}
	w^2=\dfrac{k^2 \Big(1+\dfrac{2}{45} \eps^2 k^4\Big)}{\Big(1+\dfrac{1}{3}\eps k^2+\dfrac{1}{45}\eps^2 k^4\Big)}.
\end{equation}
This has an important side effect of removing 
every 
singularity found in~\eqref{rdeB}, thus making a set of useless high-order equations applicable. An equivalent formulation of model~\eqref{eq:eGNLW1}\footnote{The equivalent formulation is obtained by multiplying both sides of the second equation of system~\eqref{eq:eGNLW1} by the water
height function, $h$.} was fully justified recently in~\cite{KLIG}. In fact, a unique solution of the model~\eqref{eq:eGNLW1} exist over the time scale of order $\frac{1}{\sqrt{\eps}}$ and stay close  the solution of the \emph{full Euler} system. Although every singularity found in~\eqref{rdeB} is now removed in~\eqref{rdeB1}, the phase velocity associated to~\eqref{rdeB1} has the same classical velocity limit of long waves as $k \rightarrow 0$, but a finite limit equal to $\sqrt{2}$ as $k \rightarrow \infty$ instead of expected zero limit.

In order to improve the model linear dispersive properties in intermediate regime of wave numbers and extend the range of applicability without changing the basic limit properties of the dispersion relation, we derive in this paper a new reformulation~\eqref{eq:eGNLW3f5} of the same order of precision of~\eqref{eq:eGNLW1} that allows the adjustment of a dispersion correction parameter $\alpha$ (see Section~\ref{IeGNeq}) and prevents the calculation of high order derivatives on $\zeta$ (see Section~\ref{reformsec}).
The dispersion relation obtained around some rest state solution associated to the new formulation~\eqref{eq:eGNLW3f5} is the following (see detailed calculation in Appendix~\ref{appendix}):
	\begin{equation}\label{rdeGN3int}
		w^2=\dfrac{k^2\Big(1+\dfrac{\eps(\alpha-1)k^2}{3} +\dfrac{\eps^2(\alpha-1)k^4}{45} +\dfrac{(7-5\alpha)\eps^2k^4}{45(1+\frac{\eps \alpha}{3}k^2)} \Big)}{\Big(1+\dfrac{\eps \alpha}{3}k^2+\dfrac{\eps^2 \alpha}{45}k^4\Big)}.
	\end{equation}
The corresponding phase velocity has in the limit when $k \rightarrow 0$ a classical velocity of long waves, and in the limit of short waves when $k \rightarrow \infty$ the phase velocity has a finite limit equal to $\dfrac{\alpha-1}{\alpha}$. Consequently, setting $\alpha=1$, the phase velocity vanishes in the limit of short waves when $k \rightarrow \infty$ as in the case of the exact linear dispersion relation~\eqref{rdFE} for the full Euler equation. This paper is devoted to the numerical study of the eB system~\eqref{eq:eGNLW3f5}. We will show that thanks to the proper choice of the dispersion correction parameter $\alpha$ (see Section~\ref{Secalphachoice}) the eB model with factorized high order derivatives~\eqref{eq:eGNLW3f5} have better dispersive properties in intermediate regime of wave numbers than the eB model without factorization~\eqref{eq:eGNLW3}. Moreover, we will show that this newly derived formulation is stable with respect to high frequency perturbations (see Section~\ref{Stabilitysec}).

The paper is organized as follows, we firstly propose a reformulation of the same order of precision of the extended Boussinesq model \eqref{eq:eGNLW1} up to the third order. This reformulation makes the model more appropriate for the numerical implementation and significantly improved in terms of linear dispersive properties due to the suitable adjustment of a dispersion correction parameter. The reformulation is performed then via the factorization of high order derivatives on the surface deformation $\zeta$. We will show that the improvement is significant in the dispersive properties of the model with factorization of high order derivatives on the surface deformation together with an appropriate choice of an optimal value of the dispersion correction parameter $\alpha$. 
	We then study the stability of two models, with and without factorization, and we will show that factorizing only the fifth order derivative presented in the second model equation 
	induces 
	a destabilizing effect :
    we need to factor 
    every 
     high order derivative on $\zeta$.
    
    Secondly, we propose a suitable Strang splitting of operators to solve the improved model : a hyperbolic part representing the Nonlinear Shallow Water system and a dispersive part representing the high order derivatives. The hyperbolic part of the system is treated with a high-order finite volume scheme whereas the dispersive part is treated with a finite difference method at the same order. This splitting strategy allowed us to calculate the time step in the first finite-volume sub-step allowing to overcome the time step restrictions induced by high order derivatives existing in the dispersive part of the model and the numerical investigations show that the time step restriction from the CFL condition (according to which the time step must be chosen proportional to the mesh spacing) of the finite-volume step is enough to ensure stability for the whole numerical method.  Moreover, a reconstruction  of nodal unknowns and 
    centered 
     unknowns is presented. 
    
    Finally, numerical validations are presented under two main goals: showing the interest of the proposed formulation of the extended Boussinesq model as well as the good behavior of the numerical scheme and assessing the potential need of such higher-order models.

	\section{Reformulation of the extended Boussinesq system}\label{REFsec}
	The system~\eqref{eq:eGNLW1} is much easier to solve numerically than the standard formulation~\eqref{ex-boussinesq}. In fact, the operator $\mathfrak{J}$ has an appropriate structure allowing its inversion. 
	 Using straightforward asymptotic expansions, the left-most term of the second equation of~\eqref{eq:eGNLW1} can be written under the form:
	\begin{equation*}\mathfrak{J}(\partial_t v + \eps v \partial_x v)= (1+\eps \mathcal{T}[0] -\eps^2\mathfrak{T})(\partial_t v+ \eps v\partial_x v) - \dfrac{2}{3} \eps^2 \zeta \partial_x^2 (\partial_t v) -\eps^2 \partial_x \zeta \partial_x (\partial_t v) + \mathcal{O}(\eps^3).\end{equation*}
Now using the fact that $\partial_{ t}  v = -\partial_x \zeta+\OO(\eps)$, one can deduce that the above equation can be recast under the following form:
	$$\mathfrak{J}(\partial_t v + \eps v \partial_x v)= (1+\eps \mathcal{T}[0] -\eps^2\mathfrak{T})(\partial_t v+ \eps v\partial_x v) + \dfrac{2}{3} \eps^2 \zeta \partial_x^3 \zeta +\eps^2 \partial_x \zeta \partial_x^2 \zeta + \mathcal{O}(\eps^3).$$
	where $\mathcal{T}[0]w=  -\dfrac{1}{3}\partial_x^2 w$ and $\mathfrak{T} w=-\dfrac{1}{45} \partial_x^4 w$.
	 Hence, system~\eqref{eq:eGNLW1} becomes:
	\begin{equation}\label{eq:eGNLW2}
		\left\{ \begin{array}{l}
			\dsp \partial_{ t}\zeta +\partial_x\big(h v\big)\ =\  0,\\ \\
			\dsp \Big(1+\eps \mathcal{T}[0] -\eps^2\mathfrak{T}\Big) \Big( \partial_{ t}  v + \eps v \partial_x v  \Big) + \partial_x \zeta+\dfrac{2}{45} \eps^2 \partial_x^5\zeta+ \dfrac{2}{3} \eps^2 \partial_x((\partial_x v)^2)\\+\dfrac{2}{3} \eps^2 \zeta \partial_x^3\zeta +\eps^2\partial_x\zeta \partial_x^2 \zeta=  \OO(\eps^3),
		\end{array} \right. \end{equation}
	where $h=1+\eps\zeta$.
	The left-most factorized operator of the second equation of~\eqref{eq:eGNLW2} also enjoys a structure allowing its inversion.  For the proof of the invertibility, one has to apply a Lax-Milgram theorem where the coercivity condition of the bilinear form is satisfied (see~\cite[Lemma 1]{KZI18}). Moreover, this operator can be inverted once for all numerical time steps because of its time-independent form. The simple one-dimensional structure of the model with a time-independent operator reduce slightly the computational time. In fact, the strategy of removing time-dependency from the left-most factorized operator was originally initiated for numerical simulations of the fully nonlinear and weakly dispersive GN models in the two-dimensional case~\cite{LannesMarche14}, in order to reduce significantly the computational time.

	\subsection{A one-parameter family of extended Boussinesq equations}\label{IeGNeq}
	The eB equations are significantly improved in terms of linear dispersive properties due to the higher-order terms existing in these equations, see~\cite{MBS03}. Additional improvement providing a finer characterization in high frequency regimes can be brought by adjusting a dispersion correction parameter $\alpha$. Following the lines in~\cite{CBB07,MMS91} and without affecting the accuracy of the model, we improve the frequency dispersion of problem~\eqref{eq:eGNLW2}. This is possible, if one adds to the second equation of~\eqref{eq:eGNLW2} some terms of the same order as the equation precision and adjusts the parameter $\alpha$ in an appropriate way. See section~\ref{Secalphachoice} for the discussion on the choice of the parameter $\alpha$.
	From the second equation of~\eqref{eq:eGNLW2}, one deduce the following approximation:
	\begin{equation}\label{eq1}\partial_t v +\eps v \partial_x v+\partial_x \zeta+\eps\mathcal{T}[0](\partial_t v) =\OO(\eps^2),\end{equation}
	where $\mathcal{T}[0]w=  -\dfrac{1}{3}\partial_x^2 w$.
	Using again the fact that $\partial_{ t}  v = -\partial_x \zeta+\OO(\eps)$, thus approximation~\eqref{eq1} can be written as:
		\begin{equation*}\partial_t v +\eps v \partial_x v+\partial_x \zeta+\dfrac{\eps}{3}\partial_x^3 \zeta =\OO(\eps^2),\end{equation*}
	and hence, for any $\alpha \in \RR^*_+$:
	\begin{equation}\label{eq2}\partial_{t}  v=\alpha\partial_{ t}  v+(\alpha-1)[\partial_x \zeta +\eps v\partial_x v+\dfrac{\eps}{3}\partial_x^3 \zeta] +(1-\alpha)\OO(\eps^2).\end{equation}
	The second equation of~\eqref{eq:eGNLW2} can be recast after substituting $\partial_t v$ by its approximation given in~\eqref{eq2}:
		\begin{multline*}
		\Big(1+\eps \mathcal{T}[0] -\eps^2\mathfrak{T}\Big) \Big(\alpha\partial_{ t}  v+\alpha \eps v \partial_x v+(\alpha-1)\partial_x \zeta +\dfrac{\eps}{3}(\alpha-1)\partial_x^3 \zeta +(1-\alpha)\OO(\eps^2)\Big)\\ +  \big(\dfrac{\alpha-1}{\alpha}+\dfrac{1}{\alpha} \big)\partial_x \zeta+\dfrac{2}{45} \eps^2 \partial_x^5\zeta +  \dfrac{2}{3}\eps^2 \partial_x((\partial_x v)^2)+\dfrac{2}{3} \eps^2 \zeta \partial_x^3 \zeta +\eps^2\partial_x\zeta \partial_x^2 \zeta=  \OO(\eps^3).
	\end{multline*}
After neglecting the terms $(\eps \mathcal{T}[0]-\eps^2\mathfrak{T})(1-\alpha)\OO(\epsilon^2)$ and $-\eps^2\mathfrak{T} \Big (\dfrac{\eps}{3}(\alpha-1) \partial_x^3\zeta\Big)$ of order $\OO(\eps^3)$, one has:		
				\begin{multline*}
\alpha\partial_{ t}  v+\alpha\eps v \partial_x v+	(\alpha-1)(\partial_x \zeta +\dfrac{\eps}{3}\partial_x^3 \zeta)+(1-\alpha)\OO(\eps^2)\\+	\Big(\eps \mathcal{T}[0] -\eps^2\mathfrak{T}\Big) \Big(\alpha\partial_{ t}  v+\alpha \eps v \partial_x v+(\alpha-1)\partial_x \zeta\Big) +\dfrac{\eps^2}{3}(\alpha-1)\mathcal{T}[0]\partial_x^3 \zeta \\+  \big(\dfrac{\alpha-1}{\alpha}+\dfrac{1}{\alpha} \big)\partial_x \zeta+\dfrac{2}{45} \eps^2 \partial_x^5\zeta +  \dfrac{2}{3}\eps^2 \partial_x((\partial_x v)^2)+\dfrac{2}{3} \eps^2 \zeta \partial_x^3 \zeta +\eps^2\partial_x\zeta \partial_x^2 \zeta=  \OO(\eps^3).
	\end{multline*}
Using the fact that $\alpha\partial_{ t}  v+\alpha\eps v \partial_x v= \partial_{ t}  v+\eps v \partial_x v+(\alpha-1)(\partial_{ t}  v+\eps v \partial_x v)$ one has:
					\begin{multline*}
\partial_{ t}  v+\eps v \partial_x v+	(\alpha-1)(\partial_{ t}  v+\eps v \partial_x v +\partial_x \zeta +\dfrac{\eps}{3}\partial_x^3 \zeta)+(1-\alpha) \OO(\eps^2)\\+	\Big(\eps\alpha \mathcal{T}[0] -\eps^2\alpha\mathfrak{T}\Big) \Big(\partial_{ t}  v+ \eps v \partial_x v +\dfrac{(\alpha-1)}{\alpha}\partial_x \zeta\Big)  +\dfrac{\eps^2}{3}(\alpha-1)\mathcal{T}[0]\partial_x^3 \zeta\\+  \big(\dfrac{\alpha-1}{\alpha}+\dfrac{1}{\alpha} \big)\partial_x \zeta+\dfrac{2}{45} \eps^2 \partial_x^5\zeta + \dfrac{2}{3} \eps^2 \partial_x((\partial_x v)^2)+\dfrac{2}{3} \eps^2 \zeta \partial_x^3 \zeta +\eps^2\partial_x\zeta \partial_x^2 \zeta=  \OO(\eps^3).
	\end{multline*}
Following straightforward computations and using the fact that $\dfrac{\eps^2}{3}(\alpha-1)\mathcal{T}[0]\partial_x^3 \zeta=-\dfrac{\eps^2}{9}(
\alpha-1)\partial_x^5 \zeta $, one has:
	\begin{multline*}
		\Big(1+\eps \alpha \mathcal{T}[0] -\eps^2\alpha\mathfrak{T}\Big) \Big(\partial_{ t}  v + \eps v \partial_x v +\dfrac{\alpha-1}{\alpha} \partial_x \zeta \Big)+(\alpha-1)(\partial_t v + \epsilon v\partial_x v +\partial_x \zeta +\dfrac{\eps}{3} \partial_x^3 \zeta)+(1-\alpha)\OO(\eps^2)\\+\dfrac{1}{\alpha} \partial_x \zeta+\dfrac{(7-5\alpha)}{45} \eps^2 \partial_x^5\zeta+ \dfrac{2}{3}  \eps^2\partial_x((\partial_x v)^2)+\dfrac{2}{3} \eps^2 \zeta \partial_x^3 \zeta +\eps^2\partial_x\zeta \partial_x^2 \zeta=  \OO(\eps^3).
	\end{multline*}
	Finally using~\eqref{eq2}, one has $(\alpha-1)(\partial_t v + \epsilon v\partial_x v +\partial_x \zeta +\dfrac{\eps}{3} \partial_x^3 \zeta)+(1-\alpha)\OO(\eps^2)=0$ and thus system~\eqref{eq:eGNLW2} with improved frequency dispersion can be written as:
	\begin{equation}\label{eq:eGNLW3}
		\left\{ \begin{array}{l}
			\dsp \partial_{ t}\zeta +\partial_x\big(h v\big)\ =\  0,\\ \\
			\dsp \Big(1+\eps \alpha \mathcal{T}[0] -\eps^2\alpha\mathfrak{T}\Big) \Big(\partial_{ t}  v + \eps v \partial_x v +\dfrac{\alpha-1}{\alpha} \partial_x \zeta \Big) +\dfrac{1}{\alpha} \partial_x \zeta+\dfrac{(7-5\alpha)}{45} \eps^2 \partial_x^5\zeta \\+ \dfrac{2}{3}\eps^2  \partial_x((\partial_x v)^2)+\dfrac{2}{3}\eps^2 \zeta \partial_x^3 \zeta +\eps^2 \partial_x\zeta \partial_x^2 \zeta=  \OO(\eps^3).
		\end{array} \right. \end{equation}
Similarly, a significant improvement of the dispersive properties has been attained in the derivation of a
three-parameter family of GN equations, see~\cite{CLM}. In here, we will limit ourselves to the one-parameter family of eB equations~\eqref{eq:eGNLW3} for the sake of simplicity. 

	\subsection{Reformulation of the extended Boussinesq equations~\eqref{eq:eGNLW3}}\label{reformsec}
	In what follows, we derive an equivalent model to~\eqref{eq:eGNLW3} (in the sense of precision), \textit{i.e.} $\mathcal{O}(\eps^3)$, that prevents the calculation of high order derivatives on $\zeta$. 
	To this effect, 
	we call such a model \emph{eB with factorized high order derivatives}. Certainly, the model enclose high order derivatives on $\zeta$, but we make it possible not to compute them by factoring them out by $(1+\eps \alpha \mathcal{T}[0])$. The price to pay is an increase in computational cost, since one needs to solve an extra linear system but the gain is significant in extending the range of applicability. 

Using the fact that $(1+\eps \alpha \mathcal{T}[0]) (\partial_x\zeta) = \partial_x\zeta + \OO(\eps)$, one has $\partial_x \zeta =(1+\eps \alpha \mathcal{T}[0] )^{-1} (\partial_x \zeta) + \OO(\eps)$, and thus the terms $\partial_x^2 \zeta$, $\partial_x^3 \zeta$ and $\partial_x^5 \zeta$ become respectively:
	\begin{equation}\label{dx2zeta}
		\partial_{x}^2  \zeta =\partial_x\Big((1+\eps \alpha \mathcal{T}[0] )^{-1} (\partial_x \zeta)\Big)+\OO(\eps),
	\end{equation} \begin{equation}\label{dx3zeta}
		\partial_{x}^3  \zeta =\partial_x^2\Big((1+\eps \alpha \mathcal{T}[0] )^{-1} (\partial_x \zeta)\Big)+\OO(\eps),
	\end{equation}
	\begin{equation}\label{dx5zeta}
		\partial_{x}^5  \zeta =\partial_x^4\Big((1+\eps \alpha \mathcal{T}[0] )^{-1} (\partial_x \zeta)\Big)+\OO(\eps).
	\end{equation}
	Replacing $ \partial_x^2  \zeta $, $ \partial_x^3  \zeta $ and $ \partial_x^5  \zeta $ by their expression obtained in~\eqref{dx2zeta},~\eqref{dx3zeta} and~\eqref{dx5zeta} respectively  in the second equation of~\eqref{eq:eGNLW3}, one can write the eB equations with improved dispersion and factorized high order derivatives on $\zeta$ as:
	\begin{equation}\label{eq:eGNLW3f5}
		\left\{ \begin{array}{l}
			\dsp \partial_{ t}\zeta +\partial_x\big(h v\big)\ =\  0,\\ \\
			\dsp \Big(1+\eps \alpha \mathcal{T}[0] -\eps^2\alpha\mathfrak{T}\Big) \Big(\partial_{ t}  v + \eps v \partial_x v +\dfrac{\alpha-1}{\alpha} \partial_x \zeta \Big) +\dfrac{1}{\alpha} \partial_x \zeta\\+\dfrac{(7-5\alpha)}{45} \eps^2 \partial_x^4 \Big((1+\eps \alpha \mathcal{T}[0] )^{-1} (\partial_x \zeta) \Big) +  \dfrac{2}{3} \eps^2 \partial_x((\partial_x v)^2)
			\\+\dfrac{2}{3}\eps^2 \zeta \partial_x^2 \Big((1+\eps \alpha \mathcal{T}[0] )^{-1} (\partial_x \zeta) \Big) +\eps^2 \partial_x\zeta\partial_x \Big((1+\eps \alpha \mathcal{T}[0] )^{-1} (\partial_x \zeta) \Big)=  \OO(\eps^3),
		\end{array} \right. \end{equation}
	where $h=1+\eps\zeta$.
	Note, that this formulation avoids the calculation of high order derivatives, in particular on the surface deformation $\zeta$, but there are still fourth and second order derivatives (this time not on $\zeta$) in the second equation of~\eqref{eq:eGNLW3f5} and a large stencil is still needed though. Significant interest is behind the derivation of the eB formulation~\eqref{eq:eGNLW3f5}. In fact, we believe that factorizing high order derivatives, namely on the surface deformation $\zeta$, will extend the range of applicability to high frequency regimes	(see discussion in section~\ref{Secalphachoice}), while remaining stable (see section~\ref{stabsec}). In the following section, we will highlight the advantages of working with the factorized eB model~\eqref{eq:eGNLW3f5} rather than~\eqref{eq:eGNLW3}.
	\subsection{Choice of the parameter $\alpha$}\label{Secalphachoice}
	The main comparison between any asymptotic model and the \emph{full Euler} equations is performed at the stage of linear periodic plane wave solutions.
	At this point, a part of the model's properties~\cite{Stoker57} are summed up in the dispersion relation, relating the spatial wave number $k$ and the time frequency $w$.
	It comes from the earlier linearisation of the system around some rest state. Improving the dispersive characteristics of our model require a suitable choice of the parameter $\alpha$
	so that the dispersion characteristics of the \emph {full Euler system} corresponds with those of the improved eB systems at the dispersion relation level.
Following~\cite{BCLMT}, we adjust this parameter so that both phase and group velocities are minimized over a range of values of $k \in [0,K]$. This can be done by minimizing a weighted averaged
	error (see for instance~\cite{CBB07}) introduced for this reason.
	
In what follows, we will show that the eB model with factorized high order derivatives have better dispersive properties \red{when compared to models including the eB model without factorization and other lower-order models.}
	The dispersion relation corresponding to~\eqref{eq:eGNLW3} can be derived by investigating the linear behavior of small perturbation to a constant state solution $(\zetabar,\vbar)$ and then looking for the corresponding plane wave solutions of the form $(\zeta^0,v^0)e^{i(kx-wt)}$:
	\begin{equation}\label{rdeGN}
		\dfrac{(w-\eps k\vbar)^2}{\hbarr k^2}=\dfrac{\Big(1+\dfrac{\eps(\alpha-1)k^2}{3}  +\dfrac{ (6-4\alpha)\eps^2k^4}{45}-\dfrac{ 2\eps^2k^2\zetabar}{3} \Big)}{\Big(1+\dfrac{\eps \alpha}{3}k^2+\dfrac{\eps^2 \alpha}{45}k^4\Big)}.
	\end{equation}
	The choice of $\alpha$ is classically made to obtain a good matching with the dispersion relation of the full Euler equations around the rest state $(\zetabar,\vbar)=(0,0)$.
	The exact dispersion relation for the \emph{full Euler} system is recalled below:
	\begin{equation}\label{rdFE}
		w^2_{S}=\dfrac{|k|}{\sqrt{\eps}}\tanh(\sqrt{\eps}|k|).
	\end{equation}
	For small wave numbers, the Taylor expansions of~\eqref{rdeGN} (with $(\zetabar,\vbar)=(0,0)$) and~\eqref{rdFE} are equivalent and the choice of $\alpha$ does not play any role in the leading terms (see numerical test~\ref{BRHWswn}). Indeed, one has:
	\begin{align*}w^2_{\alpha,eB}=w^2_{S} & \approx  \Big( k^2-\dfrac{\eps k^4}{3} + \dfrac{2}{15}\eps^2k^6+\OO(\eps^3 k^8)\Big).\end{align*}
	%
	
%

Classically, finding an optimal value of $\alpha$ for a range of values of $k$ requires the minimization of the squared relative weighted error defined below:
	\begin{equation}\label{wae}\text{Err}= \sqrt{\int_0^{K} \dfrac{1}{k} \Big(\dfrac{ C^p_{eB}- C^p_{S}}{ C^p_{S}}
		+ \dfrac{ C^g_{eB}- C^g_{S}}{ C^g_{S}}\Big)^2 dk},    
	\end{equation} 
	over some range $k \in [0,K]$, where $C^p_S(k)$ and $C^g_S(k)$ are respectively the reference phase and group velocities associated with the Stokes linear theory. The division by $k$ will emphasize the importance of keeping errors to a
minimum in shallow water.  The linear phase and group velocities associated to~\eqref{rdeGN} are defined as:
$$C^p_{eB}(k)=\dfrac{w_{\alpha,eB}(k)}{|k|} \quad \text{and} \quad C^g_{eB}(k)=\dfrac{\textrm{d} w_{\alpha,eB}(k)}{\textrm{d} k}.$$		
	
The weighted averaged error~\eqref{wae} has an absolute minimum of ($60 \% $) for $w_{\alpha,eB}$ defined in~\eqref{rdeGN} (with $(\zetabar,\vbar)=(0,0)$) in the dispersive range $0 \leq k \leq 10$. The optimal value for $\alpha$ in this case does not play any role and $\alpha=1$ is set. This very big error shows that the eB model~\eqref{eq:eGNLW3} without factorizing high order derivatives on $\zeta$ has a limited range of applicability and thus poor dispersion properties in large wave numbers regime.
	\begin{figure}[H]
\centering
	{\includegraphics[width=1\textwidth]{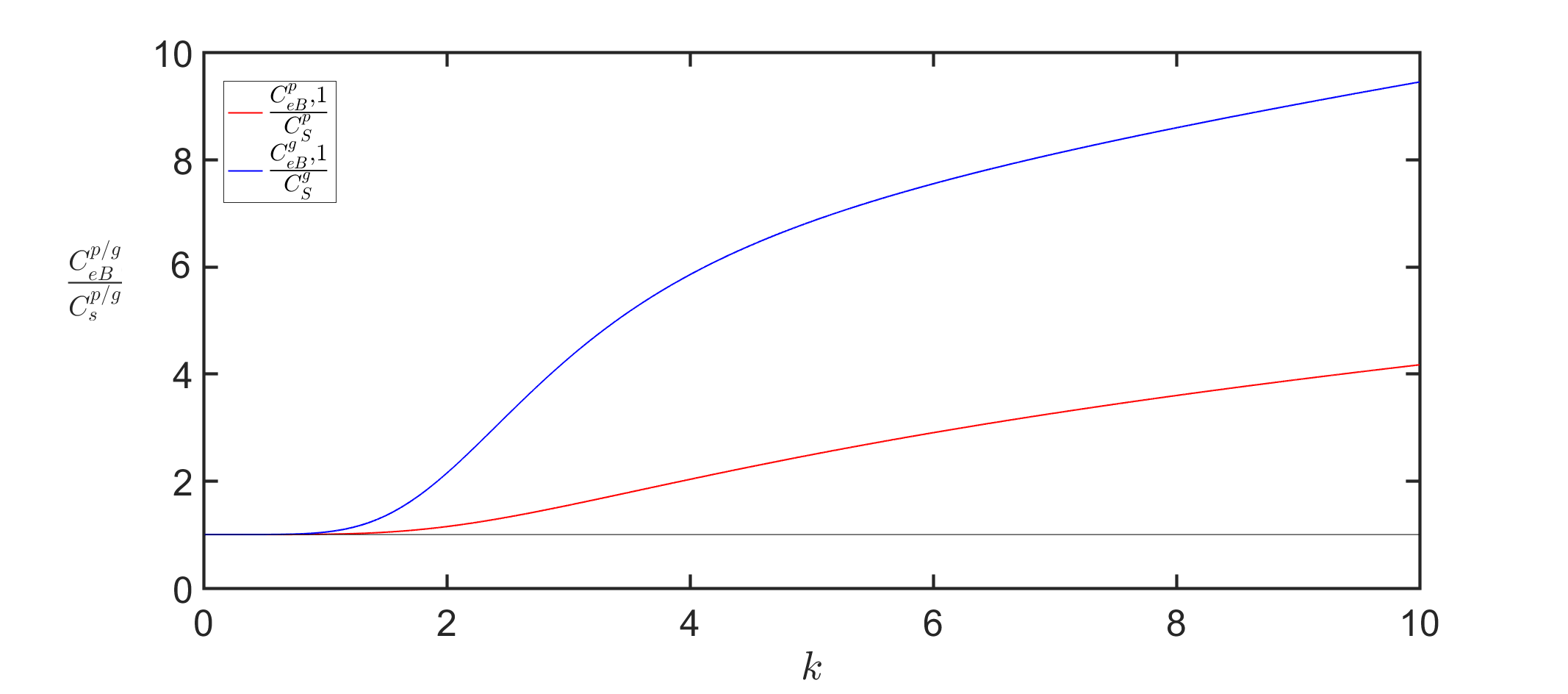}}
	\captionsetup{justification=centering}
	\caption{Errors on linear phase velocity (red) and group velocity (blue) for the eB model~\eqref{eq:eGNLW3}.}
	\label{fig1'}
\end{figure}
 In fact, in Figure~\ref{fig1'}, a clear discrepancy is observed between both ratios $\dfrac{C_{eB,1}^p}{C^p_S}$ (phase velocity, red solid line) and $\dfrac{C_{eB,1}^g}{C^g_S}$ (group velocity, blue solid line) with $\alpha=1$ when compared with the reference from Stokes theory (black solid line) when $k > 2$. On the contrary, a very good correspondence is observed in small wave numbers regime (i.e when $k \leq 1$). This show\red{s} that the eB model~\eqref{eq:eGNLW3} without factorizing high order derivatives on $\zeta$ has good dispersion properties in small wave numbers regime ($k\leq 1$). However, this is not the case in larger wave numbers regime ($k >2$). 

	Now, we discuss the dispersive properties of model~\eqref{eq:eGNLW3f5}.
The dispersion relation associated to~\eqref{eq:eGNLW3f5} is the following (see Appendix~\ref{appendix} for detailed calculation):
\begin{equation}\label{rdeGNf5}
		\dfrac{(\t{w}_{\alpha,eB}-\eps k\vbar)^2}{\hbarr k^2}=\dfrac{\Big(1+\dfrac{\eps(\alpha-1)k^2}{3} +\dfrac{\eps^2(\alpha-1)k^4}{45} +\dfrac{ (7-5\alpha)\eps^2k^4}{45(1+\frac{\eps \alpha}{3}k^2)}-\dfrac{ 2\eps^2k^2\zetabar}{3(1+\frac{\eps \alpha}{3}k^2)} \Big)}{\Big(1+\dfrac{\eps \alpha}{3}k^2+\dfrac{\eps^2 \alpha}{45}k^4\Big)}.
	\end{equation}

\begin{figure}[H]
	\centering
	{\includegraphics[width=1.0\textwidth]{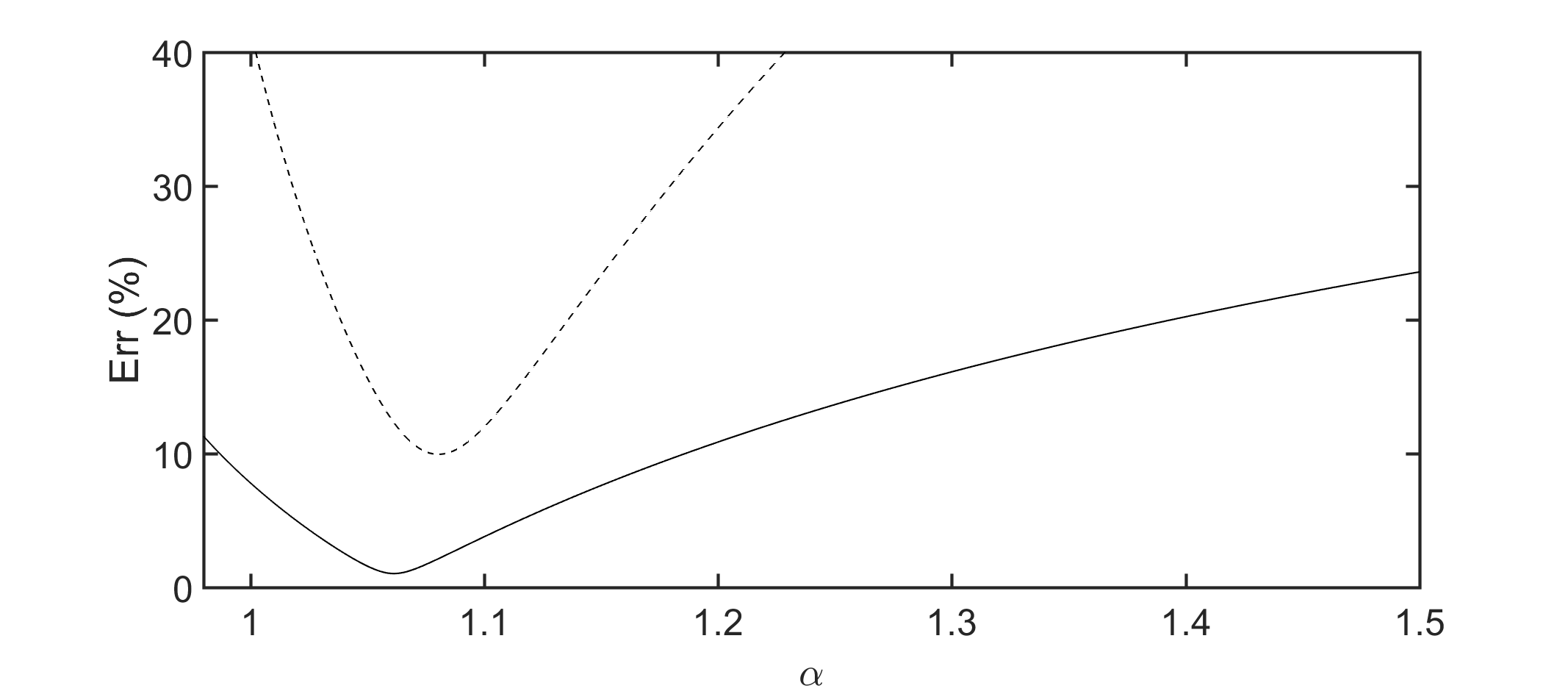}}
	\caption{Phase and group velocities weighted averaged error as a function of $\alpha$ for $0 \leq k \leq 10$. The eB model~\eqref{eq:eGNLW3f5} is in solid line, the Green-Naghdi model in the Camassa-Holm regime (GN-CH)~\cite{BGL17} is in dashes.}
	\label{fig2t}
	\end{figure}

	\noindent
	At this stage, one has to minimize the error function~\eqref{wae} for $\tilde{w}_{\alpha,eB}$ defined in~\eqref{rdeGNf5} with $(\zetabar,\vbar)=(0,0)$.  In Figure~\ref{fig2t} we plot the associated error in terms of $\alpha$. One can clearly see that the weighted average error has an absolute minimum ($\approx 1 \%$) in the dispersive range $0\leq k \leq 10$. The optimal value for $\alpha$ is 1.0610. Meanwhile, the absolute minimum of the weighted averaged error associated with the lower-order Green-Naghdi model in the Camassa-Holm regime (GN-CH) (precise up to $\mathcal{O}(\mu^2,\mu\eps^2)$)~\cite{BGL17} is much larger ($ \approx 10 \%$) with an optimal value $\alpha=1.0800$.
		\begin{figure}[H]
	\centering
		{\includegraphics[width=1\textwidth]{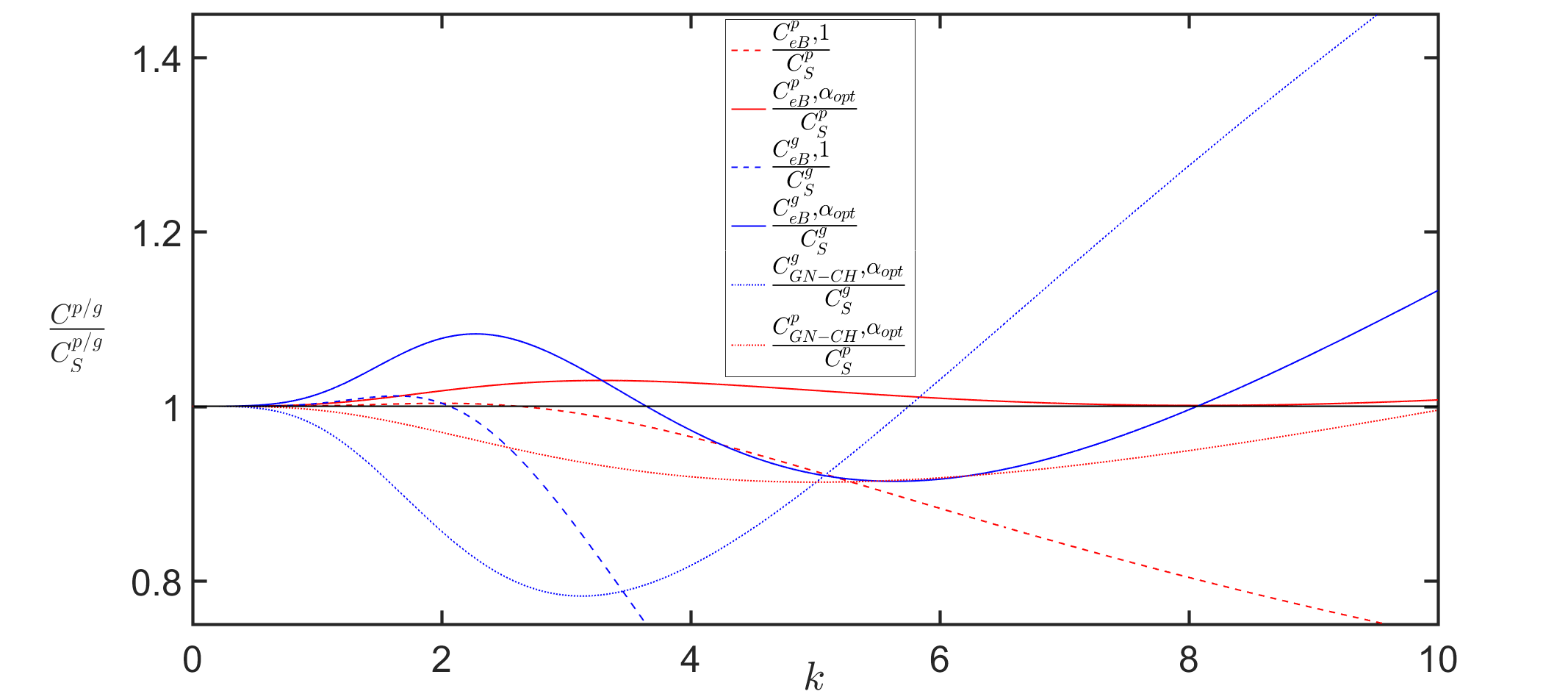}}
		\caption{Errors on linear phase velocity (red) and group velocity (blue). The reference from Stokes theory (black solid line), the eB model~\eqref{eq:eGNLW3f5} ($\alpha=1.0610$) in solid lines,  the eB model~\eqref{eq:eGNLW3f5} ($\alpha=1$) in dashes, the GN-CH model~\cite{BGL17} ($\alpha=1.0800$) in dots.}
	\label{fig2b}
	\end{figure}
	In Figure~\ref{fig2b}, errors on linear phase (red) and group (blue) velocities are plotted. The ratio $\dfrac{C_{eB,\alpha_{opt}}^p}{C^p_S}$ with an optimal choice $\alpha=1.0610$ (red solid line) is very close to the reference from Stokes theory (black solid line)
	which shows a very good correspondence between the dispersion relation obtained using the eB model~\eqref{eq:eGNLW3f5} 
	with factorized high order derivatives on $\zeta$ with $\alpha=1.0610$ and
	the theoretical one over $0\leq k \leq 10$. Larger difference exists between the group velocity (blue solid line) and the reference from Stokes theory (black solid line).
	This difference starts to proliferate when $K > 10$, showing an overestimation of this property.
%
When the modeled dispersion relation is obtained using the lower-order GN-CH model~\cite{BGL17} with $\alpha=1.0800$, a clear discrepancy exists between both linear group and phase velocity errors (blue and red dot lines) when compared with the reference from Stokes theory (black solid line). 
 We also highlight that comparisons of errors on both linear phase and group velocities of the eB model~\eqref{eq:eGNLW3f5} against other lower ordered models (for instance, sB or GN type models) are successfully reproduced and a clear superiority of our higher order model is proved when it comes to approximating properly phase and group speeds. We do not include these comparisons in the present study for the sake of shortness. In fact similar results are obtained since the improved formulation (with dispersion correction parameter) of all second order models (sB, GN-CH and GN) is the same when linearized around some rest state and thus one expects a similar linear dispersion relation.%

 In conclusion, the eB model~\eqref{eq:eGNLW3f5} contains factorized higher-order dispersive terms that are neglected in lower-order models which enlarge the application scope remarkably to cover the area from deep water (long wavelength regime) into the area of breaking waves (short wavelength regime \textit{i.e.}  high frequency), see numerical test in section~\ref{SecNumTst1}. This is the reason why, in the numerical experiments, we choose to work with the model~\eqref{eq:eGNLW3f5}, where high order derivatives are factorized, which as seen above, has an extended range of applicability and good dispersive properties in large wave numbers regime.
	\begin{remark}
In the case of variable topography, the study of dispersive properties must be supplemented with some hints about linear shoaling. In fact, the minimization of phase and group velocities errors is quite problematic in the variable topography case as a better phase and group velocities involves a deterioration of shoaling properties, which are a paramount for near-shore oceanography. 
For this reason, and in order to consider stiff configurations that include high harmonics while keeping the improved dispersive properties of the model, a three-parameter family can be derived by adding two additional parameters $\theta$ and $\gamma$ using a change of variables for the velocity~\cite{CLM}. The range of validity of the three parameter family of equations is extended to deeper water and considerable improvements in challenging configurations are obtained~\cite{LannesMarche14}. This issue is out of the scope of this article and the variable topography case will be considered in a forthcoming paper.

\end{remark}
	\subsection{Stability of the extended Boussinesq models}\label{Stabilitysec}
The high frequency instabilities of an improved GN-CH model are studied in~\cite{BGL17}. These instabilities are due to the third order derivative existing in the equation. One of the advantages of the eB is its stability in high frequency regime due to the presence of 
	higher order derivatives, namely derivatives of order five in $\zeta$. These terms seems to have a stabilizing effect. In what follows, we discuss qualitatively the stability of both models~\eqref{eq:eGNLW3} and~\eqref{eq:eGNLW3f5}. Since the parameter $\eps$ do not play a direct role in the stability results and for notational convenience we set $\eps=1$ throughout this section.

	\subsubsection{Stability of the extended Boussinesq model with high order derivatives}
	Before discussing the stability of the eB  models in high frequency regime, we would like to mention that the choice of $\alpha$ in the model~\eqref{eq:eGNLW3} in high frequency regime does not play an important role. In fact the latter model has poor dispersive properties in intermediate and large wave numbers regime, see discussion in Section~\ref{Secalphachoice}. Therefore, one has to choose $\alpha=1$.

The dispersion relation~\eqref{rdeGN} associated to~\eqref{eq:eGNLW3} with $\alpha=1$ reads:
	\begin{equation}\label{rdeGN1}
		\dfrac{(w- k\vbar)^2}{\hbarr k^2}=\dfrac{\Big(1 +\dfrac{ 2k^4}{45}-\dfrac{ 2 k^2\zetabar}{3} \Big)}{\Big(1+\dfrac{1}{3}k^2+\dfrac{1}{45}k^4\Big)}.
	\end{equation}
The perturbations of the rest state $(\zetabar,\vbar)=(0,0)$ are always stable as per the below dispersion relation:
	\begin{equation}\label{rdeGN1l}
		w^2=\dfrac{k^2\Big(1 + \dfrac{2}{45} k^4\Big)}{\Big(1+\dfrac{1}{3}k^2+\dfrac{1}{45}k^4\Big)}.
	\end{equation}
	However, a quick functional study shows that the numerator of the right-hand side of~\eqref{rdeGN1} becomes negative whenever $\zetabar>\dfrac{2k^4+45}{30k^2}$. Thus, as provided by~\eqref{rdeGN1}, $w$ remains real for large wave numbers, under the condition that $\zetabar<\dfrac{2k^4+45}{30k^2}$.  In the majority of the applications we have in mind, the overall surface deformation $\zetabar$  does not go beyond $\dfrac{2k^4+45}{30k^2}$, and this condition is satisfied. Actually, as $k$ gets large (for instance $k \approx 10$) in high frequency regime, the upper bound of $\zetabar$ gets also large ($\approx \dfrac{k^2}{15}=6.67$), hence extending the range of values of the overall surface deformation $\zetabar$ for which the condition is satisfied. This ensure\red{s} a numerical stability  in most of the situations considered for applications (see Figure~\ref{figcompHFI}). However, as we have mentioned in the introduction, from the dispersion relation~\eqref{rdeGN1l} one can remark that the phase velocity associated to the eB model with high order derivatives~\eqref{eq:eGNLW3} has the same classical velocity limit of long waves as $k\rightarrow 0$, but a finite limit equal to $\sqrt{2}$ instead of expected zero limit.
	\begin{remark}\label{remIns}
		Replacing $\partial_{x}^5  \zeta$ by $\partial_x^4\Big((1+\mathcal{T}[0])^{-1} (\partial_x \zeta)\Big)$ in the second equation of~\eqref{eq:eGNLW3} with $\alpha=1$, one gets the following model:
		\begin{equation}\label{eq:eGNLW3f51}
			\left\{ \begin{array}{l}
				\dsp \partial_{ t}\zeta +\partial_x\big(h v\big)\ =\  0,\\ \\
				\dsp \Big(1+\mathcal{T}[0] -\mathfrak{T}\Big) \Big(\partial_{ t}  v +  v \partial_x v \Big) + \partial_x \zeta+\dfrac{2}{45}  \partial_x^4\Big((1+\mathcal{T}[0])^{-1} (\partial_x \zeta)\Big) +  \dfrac{2}{3} \partial_x((\partial_x v)^2)\\+ \dfrac{2}{3}\zeta \partial_x^3 \zeta +  \partial_x \zeta \partial_x^2 \zeta=0.
			\end{array} \right. \end{equation}
		This replacement modifies the dispersion relation~\eqref{rdeGN1} into:
		\begin{equation}\label{rdeGN2}
			\dfrac{(w-k\vbar)^2}{\hbarr k^2}=\dfrac{\Big(1-\dfrac{2}{3}k^2\zetabar + \dfrac{k^4}{45}\big(\dfrac{2}{1+\frac{1}{3} k^2}\big) \Big)}{\Big(1+\dfrac{1}{3}k^2+\dfrac{1}{45}k^4\Big)}.
		\end{equation}
		A similar functional study to the previous one shows that
		the r.h.s numerator of~\eqref{rdeGN2} is negative whenever $\zetabar>\dfrac{2k^4+15k^2+45}{10k^4+30k^2}$. In high frequency regime, the upper bound of the overall surface deformation $\zetabar$ is approximately close to $0.2$, hence reducing the range of values of $\zetabar$ for which the condition is satisfied, namely $-1<\zetabar<0.2$ (keeping mind that $h=1+\zetabar$ should remain always positive). Therefore, if this condition is not satisfied the complex square root of $w$ will generate an instability in the model. Actually, a stability in high frequency regime is ensured if the condition $-1<\zetabar<0.2$ is satisfied which we believe is a limitation for the applications that we have in mind. Thus, one can deduce that factorizing only the fifth order derivative present in the second equation of~\eqref{eq:eGNLW3} does not stabilize the model, at least for a big range of values of the overall surface deformation $\zetabar$ (see Figure~\ref{figcompHFI}). 
		
		At this stage, one may wonder if a factorization of only the third order derivative on $\zeta$ may be enough to control the incriminated sign in the dispersion relation~\eqref{rdeGN1}. Indeed, this is true and $w$ remains real for large wave numbers, under the condition that $\zetabar < \dfrac{2k^6+ 6k^4+45k^2+135}{90k^2}$ (for $k \approx 10$ the upper bound of the latter inequality becomes very large). However, factoring only the third order derivative leads to a model with the same poor linear dispersive properties as~\eqref{eq:eGNLW3}. This was expected because the third order derivative term is nonlinear.  This is why we suggested in section~\ref{reformsec} to factorize second, third and fifth order derivatives present in the second equation of~\eqref{eq:eGNLW3}. The stability of the eB model with factorized high order derivatives~\eqref{eq:eGNLW3f5} is discussed in the next subsection.
	\end{remark}
	\subsubsection{Stability of the extended Boussinesq model with factorized high order derivatives}\label{stabsec}
	In what follows, we explore the stability of the eB model with factorized high order derivatives~\eqref{eq:eGNLW3f5}. We recall here the dispersion relation~\eqref{rdeGNf5} associated to the eB model~\eqref{eq:eGNLW3f5}:
	\begin{equation}\label{rdeGN3}
		\dfrac{(w-k\vbar)^2}{\hbarr k^2}=\dfrac{\Big(1+\dfrac{(\alpha-1)k^2}{3}-\dfrac{2k^2\zetabar}{3(1+\frac{\alpha}{3}k^2)} +\dfrac{(\alpha-1)k^4}{45} +\dfrac{(7-5\alpha)k^4}{45(1+\frac{\alpha}{3}k^2)} \Big)}{\Big(1+\dfrac{\alpha}{3}k^2+\dfrac{\alpha}{45}k^4\Big)}.
	\end{equation}\\The \emph{r.h.s} numerator in~\eqref{rdeGN3} is positive if and only if $$\zetabar < \dfrac{k^2\left({\alpha}\left(\left({\alpha}-1\right)k^2+15{\alpha}-12\right)+(18-15\alpha)\right)+90{\alpha}-45}{90}+\dfrac{3}{2k^2}.$$
	We recall that $\alpha >1$. Indeed, improving the dispersive properties of the model~\eqref{eq:eGNLW3f5} in large frequency regime requires an appropriate choice $\alpha=1.0610$ (see discussion in section~\ref{Secalphachoice}). With this choice of $\alpha$ and in high frequency regime, the \emph{r.h.s} of the above inequality becomes very large, namely same order as $\dfrac{\alpha(\alpha-1)k^4}{90}$. Thus, relaxing the stability condition on $\zetabar$.  This is another reason why, for the rest of the paper, we choose to work with the model~\eqref{eq:eGNLW3f5}, where high order derivatives are factorized, which as seen above is stable for the majority of applications we have in mind. Of course, one could use model~\eqref{eq:eGNLW3} which seems to be stable, but at the price of losing the improved dispersive properties 
	that model~\eqref{eq:eGNLW3f5} enjoys, see section~\ref{Secalphachoice} and test case of section~\ref{SecNumTst1}. As we have already mentioned in the introduction, one can remark from the dispersion relation~\eqref{rdeGN3} around some rest state solution $(\zetabar,\vbar)=(0,0)$ that the phase velocity associated to~\eqref{eq:eGNLW3f5} has the same classical velocity limit of long waves as $k \rightarrow 0$, namely 1, and a finite limit equal to $\sqrt{\dfrac{\alpha-1}{\alpha}}$ as $k \rightarrow \infty$ which is equal 0 (as in the case of the exact linear dispersion relation for the full Euler equation~\eqref{rdFE}) when the dispersion correction parameter is set to $\alpha=1$.
	
\begin{figure}[H]
\centering
	{\includegraphics[width=0.7\textwidth]{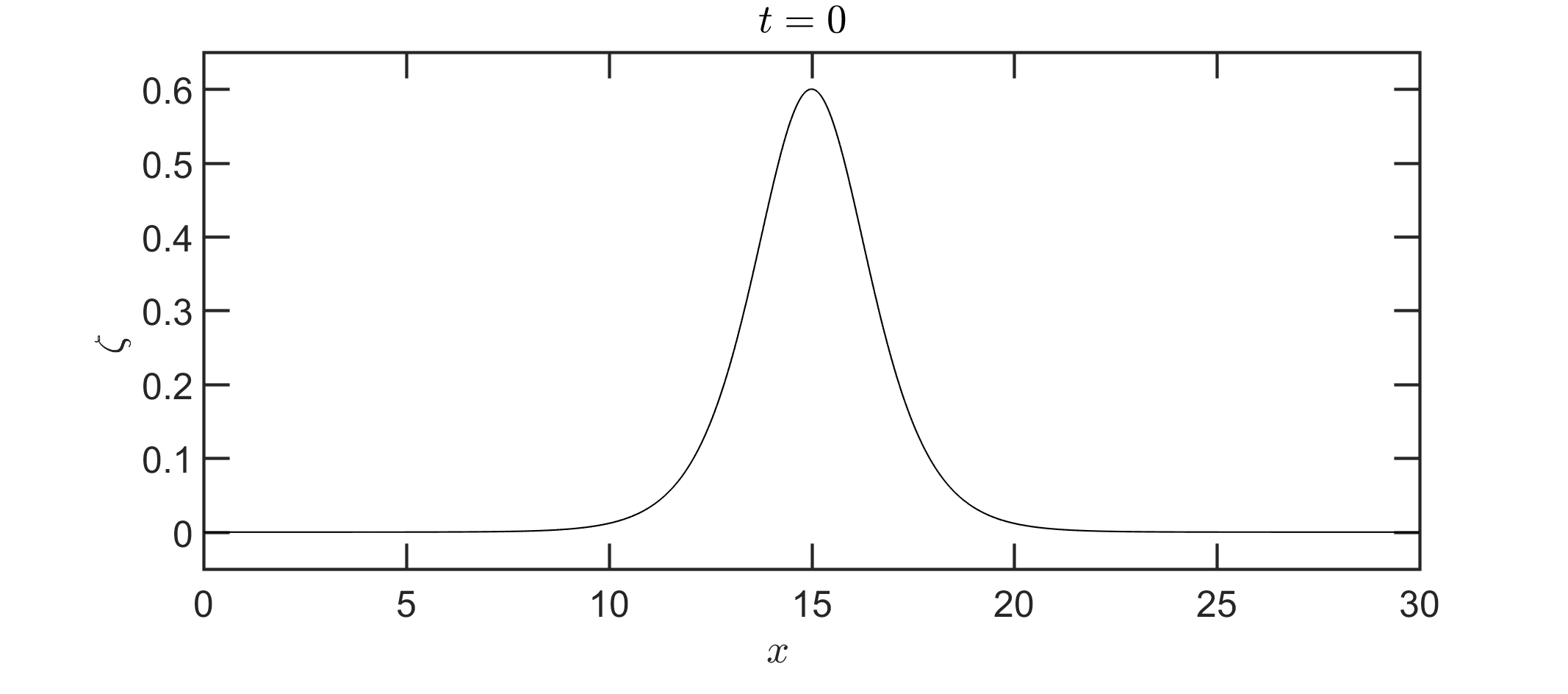}}

	{\includegraphics[width=0.7\textwidth]{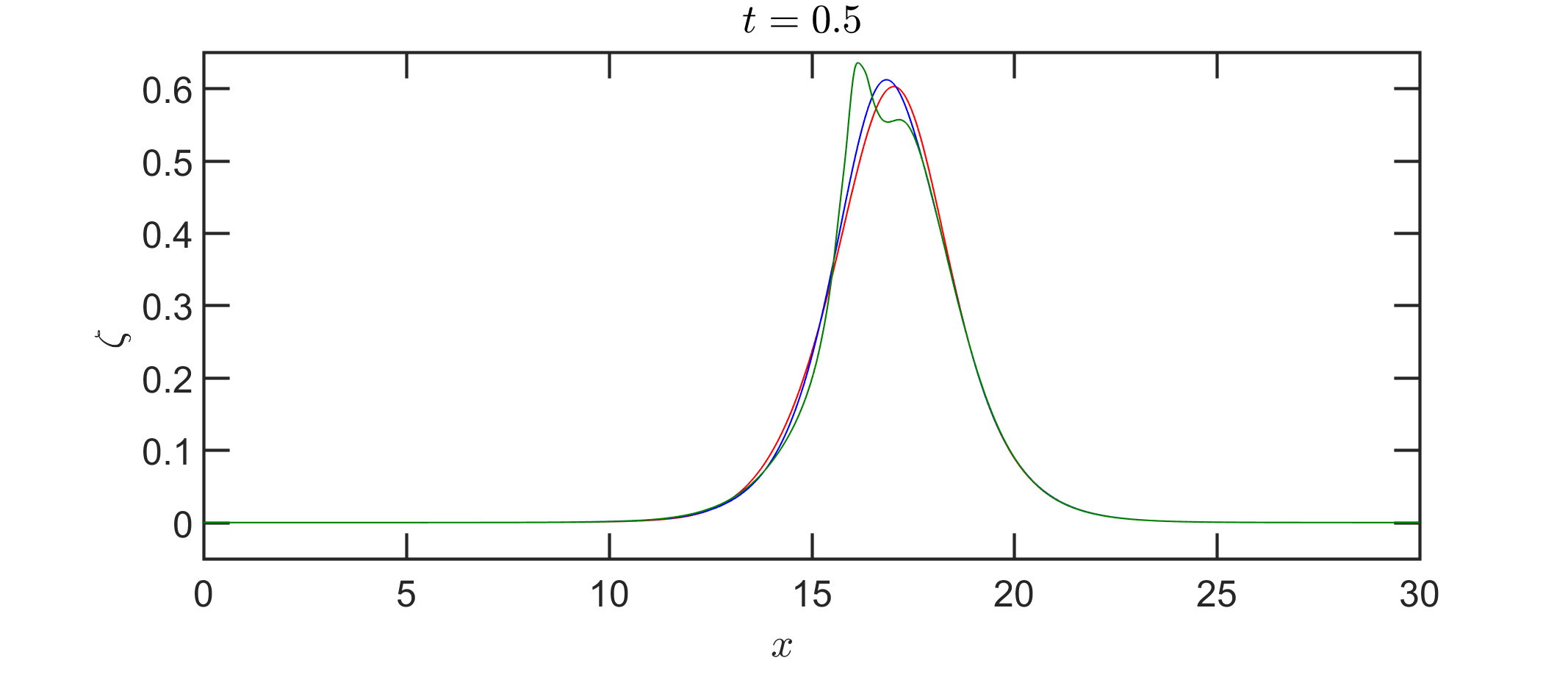}}
	
	{\includegraphics[width=0.7\textwidth]{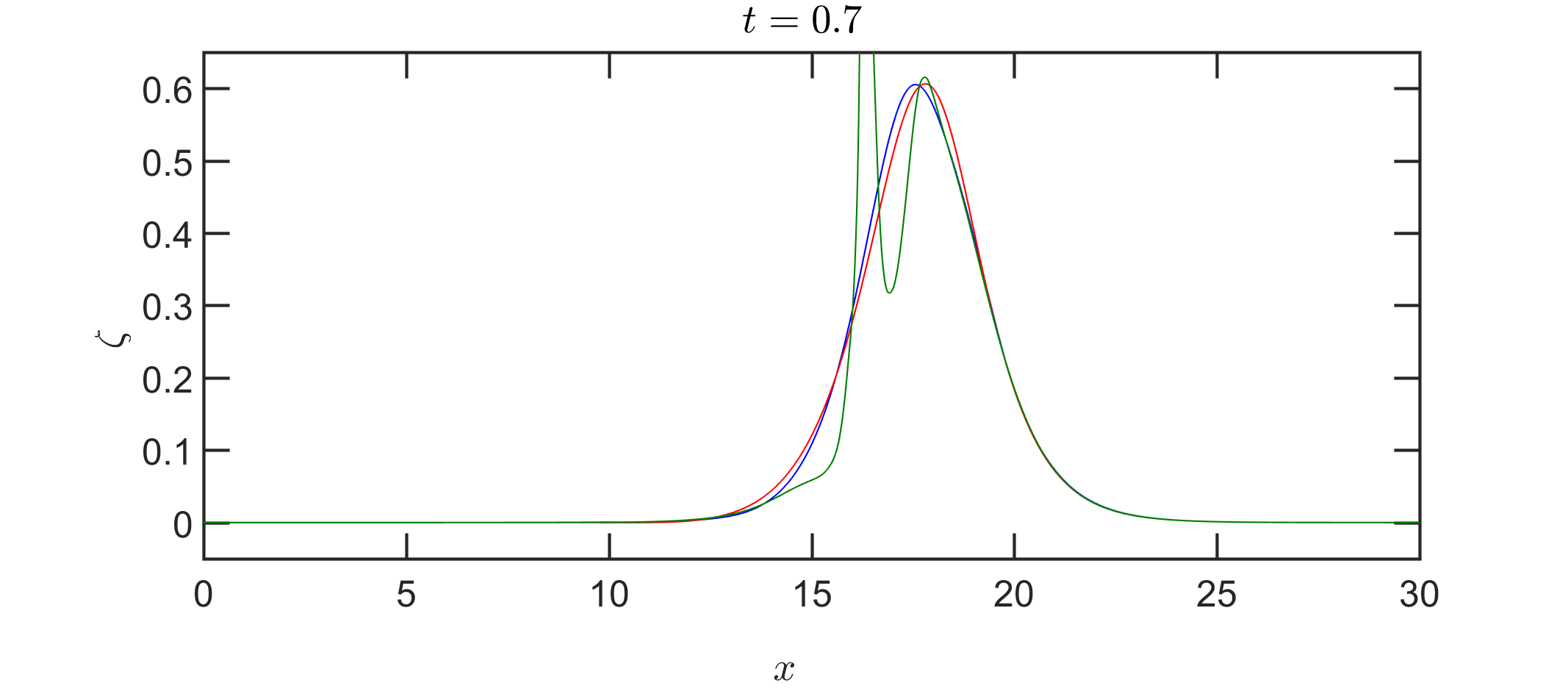}}
	\captionsetup{justification=centering}
	\caption{Comparison at different times between the solutions of the models~\eqref{eq:eGNLW3} (blue line),~\eqref{eq:eGNLW3f5} (red line) and~\eqref{eq:eGNLW3f51} (green line) in high frequency regime.}
	\label{figcompHFI}
\end{figure}

		Figure~\ref{figcompHFI} shows the left to right propagation of a solitary wave initially centered at $x_0 = 15$ of amplitude $a = 0.6$. The computational domain length is $L=30$ and discretized with 480 cells. The water surface profiles of our numerical solutions provided by the models~\eqref{eq:eGNLW3} (blue line),~\eqref{eq:eGNLW3f5} (red line) and~\eqref{eq:eGNLW3f51} (green line) are compared at $t=0.5$ and $t=0.7$ using the fifth order discretization ``WENO5-DF4-RK4" (see sections~\ref{Hyperbolic} and~\ref{Dispersive}). One can clearly see the stability of the eB models~\eqref{eq:eGNLW3} and~\eqref{eq:eGNLW3f5} while the model~\eqref{eq:eGNLW3f51} seems 
	to be unstable in high frequency regime. In fact, when implementing in model~\eqref{eq:eGNLW3f51} an initial solution that does not satisfy the limiting stability condition 
	discussed in remark~\ref{remIns}, more precisely when choosing $\zetabar=0.6 > 0.2$,  one observes a high frequency instability. We would like to mention that we tried the 
	same test but starting with an initial solution where the overall surface deformation is $\zetabar=0.1$. All models seems to be stable but we do not include this test here for the sake of simplicity.

	\section{Numerical methods}\label{NMSec}
	In what follows, we will just introduce the numerical scheme devoted to solve the eB model with factorized high order derivatives~\eqref{eq:eGNLW3f5} in order to ease the reading. A similar numerical scheme is adopted when the eB model with high order derivatives~\eqref{eq:eGNLW3} is concerned.
	
The remarkable structure of the eB models make them suitable for the implementation of a hybrid scheme splitting the hyperbolic and dispersive parts of the equations. This strategy has been initially introduced for Boussinesq-like and Green-Naghdi equations in order to handle correctly wave breaking that occurs as waves approach
	the shore, see~\cite{BCLMT,LannesMarche14}.
\red{ A computation of a half-time step of the hyperbolic part is used as a sensor to evaluate the energy loss occurring 
	during wave breaking (accurate detection of wave fronts), see~\cite{TBMCL12}. Near the breaking points, the dynamics of the waves are described correctly using the hyperbolic
	part but the dispersive components of the equation become very singular. In order to handle wave breaking, switching from the dispersive part to the hyperbolic part is indispensable.}
	In this paper, we do not investigate breaking waves. In fact, our work is limited to the flat topography case and we leave for future research 
	works 
	 the treatment of breaking waves in the variable bottom
	configuration. However, we stick here to the splitting strategy since it is computationally efficient, stable and cheap. 
	 Moreover, this splitting strategy allowed us to overcome the severe time step restriction induced due to the presence of high order derivatives by calculating the time step in the first finite-volume sub-step. The numerical investigations show that the time step restriction from the CFL condition (according to which the time step must be chosen proportional to the mesh spacing) of the finite-volume step is enough to ensure stability for the whole numerical method.%
	  The splitting scheme following the lines in~\cite{BCLMT,LannesMarche14,BGL17}
	is presented in the section below.
	\subsection{The splitting scheme}\label{SSsec}
	We recall first the eB system~\eqref{eq:eGNLW3f5} under consideration:
			\begin{equation}\label{eq:eGNLW3f5NM}
		\left\{ \begin{array}{l}
			\dsp \partial_{ t}\zeta +\partial_x\big(h v\big)\ =\  0,\\ \\
			\dsp \Big(1+\eps \alpha \mathcal{T}[0] -\eps^2\alpha\mathfrak{T}\Big) \Big(\partial_{ t}  v + \eps v \partial_x v +\dfrac{\alpha-1}{\alpha} \partial_x \zeta \Big) +\dfrac{1}{\alpha} \partial_x \zeta\\+\dfrac{(7-5\alpha)}{45} \eps^2 \partial_x^4 \Big((1+\eps \alpha \mathcal{T}[0] )^{-1} (\partial_x \zeta) \Big) +  \dfrac{2}{3} \eps^2 \partial_x((\partial_x v)^2)
			\\+\dfrac{2}{3}\eps^2 \zeta \partial_x^2 \Big((1+\eps \alpha \mathcal{T}[0] )^{-1} (\partial_x \zeta) \Big) +\eps^2 \partial_x\zeta\partial_x \Big((1+\eps \alpha \mathcal{T}[0] )^{-1} (\partial_x \zeta) \Big)=  \OO(\eps^3),
		\end{array} \right. \end{equation}
	where $h=1+\eps\zeta$, $\mathcal{T}[0]w=  -\dfrac{1}{3}\partial_x^2 w$ and $\mathfrak{T} w=-\dfrac{1}{45} \partial_x^4 w$.\\
	The solution operator $S(.)$ related to~\eqref{eq:eGNLW3f5NM} is decomposed at each time step $\Delta t$ following a hybrid Strang splitting scheme:
	$$S(\Delta t) = S_1(\Delta t/2)S_2(\Delta t)S_1(\Delta t/2).$$
	$\bullet \ S_1(t)$ is the solution operator related to the hyperbolic nonlinear shallow water equations, NSWE:
	\begin{equation}\label{hyp}
		\left\{ \begin{array}{l}
			\dsp \partial_{ t}\zeta +\partial_x\big(h v\big)\ =\  0,\\ \\
			\dsp  \partial_{ t}  v + \eps v \partial_x v +\dfrac{\alpha-1}{\alpha}  \partial_x \zeta +\dfrac{1}{\alpha} \partial_x \zeta=0.
		\end{array} \right. \end{equation}
	The NSWE system~\eqref{hyp} can be written in the following conservative form:
	\begin{equation}\label{hypcons}\left\{ \begin{array}{l}
			\dsp \partial_{ t}\zeta +\partial_x\big(h v\big)\ =\  0,\\ \\
			\dsp \partial_{ t}  v+\partial_x\Big(\dfrac{\eps}{2}v^2+\zeta\Big) =0,
		\end{array} \right. \end{equation}
	where $h=1+\eps\zeta$.\\
	\\
	$\bullet \ S_2(t)$ is the solution operator related to the remaining (dispersive) part of the equations.
	\begin{equation}\label{disp}\left\{ \begin{array}{l}
			\dsp \partial_{ t}\zeta \ =\  0,\\ \\
			\dsp \Big(1+\eps \alpha \mathcal{T}[0] -\eps^2\alpha\mathfrak{T}\Big) \Big(\partial_{ t}  v  -\dfrac{1}{\alpha} \partial_x \zeta \Big) +\dfrac{1}{\alpha} \partial_x \zeta\\+\dfrac{(7-5\alpha)}{45} \eps^2 \partial_x^4 \Big((1+\eps \alpha \mathcal{T}[0] )^{-1} (\partial_x \zeta) \Big) +  \dfrac{2}{3} \eps^2 \partial_x((\partial_x v)^2)
			\\+\dfrac{2}{3}\eps^2 \zeta \partial_x^2 \Big((1+\eps \alpha \mathcal{T}[0] )^{-1} (\partial_x \zeta) \Big) +\eps^2 \partial_x\zeta\partial_x \Big((1+\eps \alpha \mathcal{T}[0] )^{-1} (\partial_x \zeta) \Big)=0
		\end{array} \right. \end{equation}
	The hyperbolic conservative structure of system~\eqref{hypcons} allows a computation of $S_1$ following a finite-volume method. Whereas, a classical finite difference method is used to compute $S_2$. 
		\begin{remark} Treating the two-dimensional case is not the objective of this paper, however, the extension of the splitting approach to two-dimensional surface waves does not raise theoretical difficulty.
Following the steps of~\cite{LannesMarche14}, the eB model~\eqref{eq:eGNLW3f5NM} in the two-dimensional case can be reformulated in a way that is suitable to the implementation of the splitting strategy with the benefit of removing numerical obstructions. Indeed, the computation of the dispersive part in the above splitting scheme requires the inversion of a fourth order differential operator. This operator is a matricial operator that can be replaced by a new one having a diagonal structure whose inversion is not computationally demanding; numerically this is equivalent to the resolution of two sparse linear systems. Moreover, time dependency can be removed from this operator while keeping its diagonal structure so that it has not to be modified at each time step.
	\end{remark}
		\subsection{Finite volume scheme}\label{Hyperbolic}
	The hyperbolic system~\eqref{hypcons} is conveniently rewritten with conservative variables and a flux function:
	\begin{equation}\label{condform}\partial_t U +\partial_x (F(U))=0,\end{equation}
	where,
	\begin{equation}\label{consvar}
		U=\begin{pmatrix}
			\zeta\\
			v
		\end{pmatrix}
		,\quad F(U)=\begin{pmatrix}
			hv\\
			\dfrac{\eps}{2} v^2+g\zeta
		\end{pmatrix},
	\end{equation}
	with $h=1+\eps \zeta$. The Jacobian matrix is given by:
	\begin{equation}\label{matricejacob}
		A(U)=d(F(U))=\begin{pmatrix}
			\eps v& h\\
			g & \eps v
		\end{pmatrix}.
	\end{equation}
	The homogeneous system~\eqref{condform} is strictly hyperbolic if $\displaystyle{\inf_{x\in \RR} h >\ 0}$ that is to say the domain of the fluid must remain strictly connected.\\
	\\
	The Cauchy problem associated to~\eqref{condform} is the following:
	\begin{equation}\label{cauchy}\left\{ \begin{array}{l}
			\partial_t U +\partial_x (F(U)) \ =\  0, \qquad  t\geq0, x\in \RR.\\ \\
			U(0,x)=U_0(x), \qquad  x \in \RR.
		\end{array} \right. \end{equation}
	The finite volume method used to the approximation of~\eqref{cauchy} imposes conservation laws in a one-dimensional control volume $[x_{i-1/2},x_{i+1/2}] \times [t^n, t^{n+1}]$
	of dimensions $\Delta x=x_{i+1/2}-x_{i+1/2}$ and $\Delta t=t^{n+1}-t^{n}$.
	
	\begin{figure}[H]
		\begin{center}
			\includegraphics[width=0.7\textwidth]{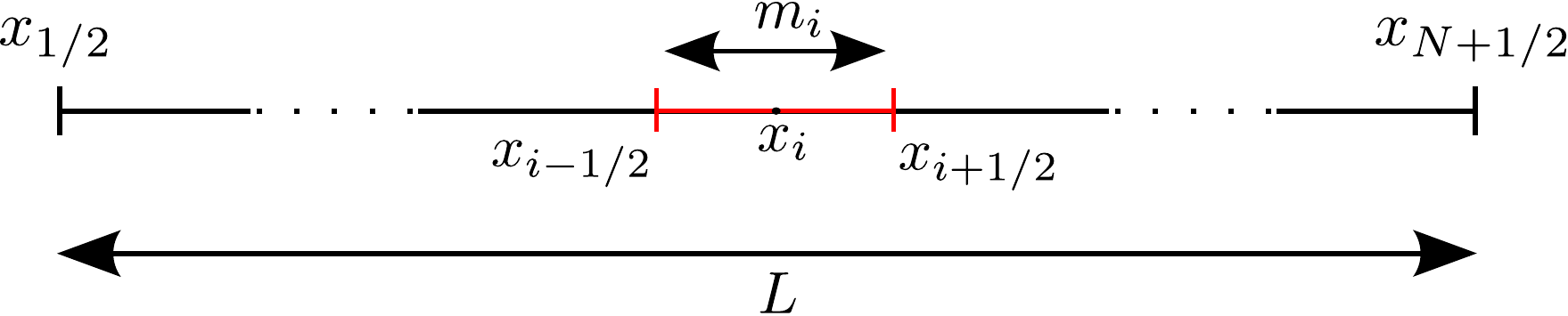}
			\caption{The space discretization.}
			\label{domain-discret}
		\end{center}
	\end{figure}
	\noindent
	The approximate cell average of $U$ on the cell $m_i=[x_{i-1/2},x_{i+1/2}] $ at time $t$ is denoted by $\overline U_i$ and given by:
	$$\overline U_i = \frac{1}{\Delta x}\int_{m_i}U(t,x)\,dx \ .$$
	The approximate cell average  of $U$ on the cell $m_i=[x_{i-1/2},x_{i+1/2}] $ at time $t^n$ is denoted by $\overline U^n_i$ and given by:
	$$\overline U^n_i = \frac{1}{\Delta x}\int_{m_i}U(t^n,x)\,dx \ .$$
	Integrating~\eqref{hyp} over the computational cell $m_i$, the semi-discrete form can be represented as:
	\begin{equation}\label{volfinidiscret}
		\frac{d \overline{U}_i(t)}{dt} + \frac{1}{\Delta x}\Big(F_{i+1/2} - F_{i-1/2} \Big) = 0
	\end{equation}
	where $F_{i\pm 1/2}$ are the numerical fluxes defined at each cell interface as:
	\begin{equation}\label{flux} F_{i+1/2}=\widetilde{F}(\overline U_i,\overline U_{i+1})\approx\dfrac{1}{\Delta x} \int_{m_i} F(U(t,x_{i+1/2}))dx.\end{equation}
	\paragraph{VFRoe method}
	In what follows, we consider the numerical approximation of the hyperbolic system of conservation laws in the form of~\eqref{condform}. To this end, we adopt the VFRoe method (see~\cite{BGH00,GHN02,GHN03}) which is an approximate Godunov scheme. It relies on the exact 
	resolution of the following linearized Riemann problem:
	\begin{equation}\label{PRL}\left\{ \begin{array}{l}
	\dsp \partial_{ t}U + \widetilde{A}(\overline U_i^n, \overline U_{i+1}^n)\partial_x U\ =\  0,\\ \\
	\dsp U(0,x)=\left \{ \begin{array}{l}
	\dsp \overline U^n_i \quad if \quad x<x_{i+1/2},\\ \\
	\dsp \overline U_{i+1}^n \quad if \quad x>x_{i+1/2},
	\end{array}\right.\end{array} \right. 
	\end{equation}
	where $\widetilde{A}(\overline U_i^n,\overline U_{i+1}^n)=A\left(\dfrac{\overline U_i^n+\overline U_{i+1}^n}{2}\right)$.\\
	By solving the linearized Riemann problem we obtain $\overline U_{i+1/2}^*= U(x=x_{i+1/2},t=t_{n})$, the  interface value between two neighboring cells. 

	\paragraph{CFL condition}\label{CFLpar}
	It is always necessary to impose what is called a CFL condition (for
	Courant, Friedrichs, Levy) on the time step to prevent the blow up of
	the numerical values. It comes usually under the form
	\begin{equation}\label{CFL}a_{i+1/2}\Delta t \leq \Delta x, \quad i=1,\ldots,N,\end{equation}
	where $a_{i+1/2}=\dsp \max_{i\in[1, N]}(j=1,2, |\lambda_j(\widetilde U_i)|)$ and $\lambda_j(\widetilde U_i)$ are the eigenvalues of $A\big(\widetilde U_i=\dfrac{\overline U_i^n+ \overline U_{i+1}^n}{2}\big)$.\\
	The restriction~\eqref{CFL} enables in practice to compute the time step
	at each time level $t_n$, in order to determine the new time level $t_{n+1} =
	t_n + \Delta t$ (within this view, $\Delta t$ is not constant, it is computed in an
	adaptive fashion).
	\paragraph{Consistency}
	The numerical flux $\widetilde{F}(U_l,U_r)$ is called consistent with~\eqref{condform} if
	\begin{equation}\label{consis}\widetilde{F}(U,U)=F(U) \quad \mbox{for all U}.\end{equation}
	


	\subsubsection{First order finite-volume scheme}\label{FV1sec}
	The semi-discrete equation \eqref{volfinidiscret} is discretized by an explicit Euler (in time) method to obtain:
	\begin{equation}\label{volfini}
	\overline U^{n+1}_i= \overline U^n_i-\dfrac{\Delta t}{h_i}(F_{i+1/2}^n-F_{i-1/2}^n),
	\end{equation}
	where the numerical flux is defined directly as the value of the exact flux at the interface value, namely: 
	\begin{eqnarray}\label{numflux}
	F_{i+1/2}^n=\widetilde{F}(\overline U_i^n,\overline U_{i+1}^n)=F(\overline U_{i+1/2}^*)\nn \\ 
	\\F_{i-1/2}^n=\widetilde{F}(\overline U_{i-1}^n,\overline U_i^n)=F(\overline U_{i-1/2}^*)\nn.
	\end{eqnarray}	
Let us remark that by construction the numerical flux given by \eqref{numflux} ensures the consistency property.
	
	In the sequel, we will suppose that the space discretization is uniform.
	\paragraph{Algorithm}\label{algo} In the following, we state the algorithm  for computing the discrete values $\overline U_i^{n+1}$ at $t^{n+1}$. Given the initial data and boundary conditions and the number $CFL \leq 1$, we start with the known discrete averaged values $(\overline{ U}_i^n)$ for $i=0,...,N+1$ at $t^n$. As long as ($t<T$) one has to do:\\
	\\
	1) Computation of $\widetilde{A}_i$ for $i=1,...,N+1$ where $\widetilde{A}_i=A\left(\dfrac{\overline U_{i-1}^n+\overline U_{i}^n}{2}\right)$.\\
	\\
	2) Computation of $r_i^1$,$r_i^2$ and $\lambda_i^1$, $\lambda_i^2$ set respectively as the eigenvectors and eigenvalues of $\widetilde{A}_i$.\\
	\\
	3) Computation of $\Delta t$, such that $\dfrac{\Delta t}{\Delta x} \leq \dfrac{CFL}{a_{i+1/2}}$.\\
	\\
	4) Computation of ${\overline U_{i-1/2}^*}$ for $i=1,...,N+1$ by solving the linearized Riemann problem.\\
	\\
	In fact we have 3 cases:\\
	\\
	$\bullet$ if $\lambda_i^1$,$\lambda_i^2$<0 then ${\overline U_{i-1/2}^*}=\overline U_i^n$.\\
	\\
	$\bullet$ if $\lambda_i^1$,$\lambda_i^2$>0 then ${\overline U_{i-1/2}^*}=\overline U_{i-1}^n$.\\
	\\
	$\bullet$ if $\lambda_i^1<0$, $\lambda_i^2>0$ then for:\begin{equation}\left\{ \begin{array}{l}\dsp x<\lambda_i^1t \quad \mbox{one has} \quad {\overline U_{i-1/2}^*}=\overline U_{i-1}^n, \\ \\
	        \dsp x>\lambda_i^1t \quad \mbox{or} \quad x<\lambda_i^2t  \quad \mbox{one has} \quad {\overline U_{i-1/2}^*}=\overline U_{i}^n-(R^{-1}[U])_2r_i^2=\overline U_{i-1}^n+(R^{-1}[U])_1r_i^1,\\ \\
	        \dsp x>\lambda_i^2t \quad \mbox{one has} \quad {\overline U_{i-1/2}^*}=\overline U_{i}^n,\end{array} \right. \end{equation}
	with $R=(r_i^1|r_i^2)$ and $[U]=\overline U^n_{i}-\overline U^n_{i-1}$.\\
	\\
	5) Computation of $F({\overline U_{i-1/2}^*})$.\\
	\\
	6) Computation of $\overline U^{n+1}_i=\overline U^n_i-\dfrac{\Delta t}{\Delta x}(F_{i+1/2}^n-F_{i-1/2}^n)$ for $i=1,...,N$.\\
	\\
	We repeat this algorithm for the new level of time $(t^{n+1}+\Delta t)$, until we reach the required final time $T$. 
	 
	 In what follows, the computation of high-order accurate numerical fluxes is reached by reconstructing left and right constant averaged values using a fifth-order WENO scheme, before applying the numerical flux.  The only change is in the computation of the interface values $\overline U_{i+1/2}^*$ which will depend on the high order reconstructed right and left values when solving the linearized Riemann problem.

	\subsubsection{High order finite-volume scheme: WENO5-RK4}\label{WENO5sec}


Considering numerical approximations of the hyperbolic system~\eqref{hypcons}, we seek a numerical scheme that reach high order accuracy in smooth regions  while avoiding the spurious oscillations around discontinuity. In fact, we aim at dealing with discontinuous initial data (dam-break problem) generating dispersive shock waves~\cite{MGH10,MID14} which needs a special treatment at the numerical scheme level. This can be achieved by using a fifth-order accuracy WENO reconstruction for hyperbolic conservation laws, following~\cite{JiangShu96,Shu98}. Second-order schemes tend to alter the dispersive properties of the model due to dispersive truncation errors. To prevent this in the study of dispersive waves, high order schemes are imperative~\cite{BCLMT,CLM,LannesMarche14,BGL17}.  For the sake of simplicity, the reader is referred to~\cite{BGL17} for more
details concerning the high order discretization of the hyperbolic system.

Regarding time discretization, fourth-order explicit Runge–Kutta ``RK4" method is used and one gets the ``WENO5-RK4" scheme.

	\subsection{Finite difference scheme for the dispersive part}\label{Dispersive}
	The splitting scheme is a mix between a finite volume discretization and a finite difference method. 
	This mix induces a switching between cell-averaged values defined by the finite volume discretization and 
	nodal values used by the finite difference discretization for each unknown and at each time step.
	Using fifth-order accuracy WENO reconstruction, one can approximate the nodal values (i.e finite difference unknowns) 
	$(U_i^n)_{i =1,N+1}$ in terms of the cell-averaged values (i.e finite volume unknowns) $(\overline U_i^n)_{i =1,N}$ by the following relation:
	\begin{equation}\label{switchFVDF}
		U_i^n=\dfrac{1}{30}\overline U_{i-2}^n-\dfrac{13}{60}\overline U_{i-1}^n+\dfrac{47}{60}\overline U_i^n+\dfrac{9}{20}\overline U_{i+1}^n-\dfrac{1}{20}\overline U_{i+2}^n +\OO(\Delta x^5), \quad 1\leq i \leq N+1,
	\end{equation}
	 One can easily recover the relation that allows to determine the cell-averaged values $(\overline U_i^n)_{i =1,N}$ in terms of the nodal values $(U_i^n)_{i =1,N+1}$ by inverting the relation~\eqref{switchFVDF}.
	The global order of the scheme is preserved. In fact, the formula is precise up to order $\OO(\Delta x^5)$ terms.
	Before proceeding by the computation of $S_2$, we recall first the remaining (dispersive part) of the equations, given in section~\ref{SSsec}.
	\begin{equation}\label{disp1}\left\{ \begin{array}{l}
			\dsp \partial_{ t}\zeta \ =\  0,\\ \\
			\dsp \Big(1+\eps \alpha \mathcal{T}[0] -\eps^2\alpha\mathfrak{T}\Big) \Big(\partial_{ t}  v  -\dfrac{1}{\alpha} \partial_x \zeta \Big) +\dfrac{1}{\alpha} \partial_x \zeta\\+\dfrac{(7-5\alpha)}{45} \eps^2 \partial_x^4 \Big((1+\eps \alpha \mathcal{T}[0] )^{-1} (\partial_x \zeta) \Big) +  \dfrac{2}{3} \eps^2 \partial_x((\partial_x v)^2)
			\\+\dfrac{2}{3}\eps^2 \zeta \partial_x^2 \Big((1+\eps \alpha \mathcal{T}[0] )^{-1} (\partial_x \zeta) \Big) +\eps^2 \partial_x\zeta\partial_x \Big((1+\eps \alpha \mathcal{T}[0] )^{-1} (\partial_x \zeta) \Big)=0.
		\end{array} \right. \end{equation}
Inverting the operator $(1+\eps \alpha \mathcal{T}[0])$ requires the discretization of $(1+\eps \alpha \mathcal{T}[0])X=B$ (i.e resolution of linear systems). The linear system is solved using the matrix division operator of Matlab, $X=(1+\eps \alpha \mathcal{T}[0]) \setminus B$. This produces the solution using Gaussian elimination, without forming the inverse.

	Using an explicit Euler in time scheme, the finite discretization of the system~\eqref{disp1} using classical finite difference methods leads to the following discrete problem:
	\begin{equation}\label{disp1disc}\left\{ \begin{array}{l}
			\dsp \dfrac{\zeta^{n+1}-\zeta^n}{\Delta t} \ =\  0,\\ \\
			\dsp \dfrac{v^{n+1}-v^n}{\Delta t}-\dfrac{1}{\alpha} D_1 (\zeta^{n})
			+\Big(1-\frac{\eps \alpha}{3} D_2+\frac{\eps^2 \alpha}{45}D_4\Big)^{-1} \Big[ \frac{1}{\alpha}D_1(\zeta^{n})\\+\dfrac{(7-5\alpha)}{45}\eps^2 D_4\Big( \big(1-\dfrac{\eps \alpha}{3}D_2\big)^{-1} (D_1(\zeta^n))\Big)+  \dfrac{2}{3}\eps^2 D_1((D_1 v^n)^2)\\+ \dfrac{2}{3}\eps^2 \zeta^n  D_2\Big( \big(1-\dfrac{\eps \alpha}{3}D_2\big)^{-1} ( D_1(\zeta^n))\Big)\\+\eps^2 D_1(\zeta^n)  D_1\Big( \big(1-\dfrac{\eps \alpha}{3}D_2\big)^{-1} ( D_1(\zeta^n))\Big) \Big]=0.
		\end{array} \right. \end{equation}
	The matrices $D_1$, $D_2$ and $D_4$ are respectively the classical centered discretizations of the derivatives $\partial_x$, $\partial^2_x$, and $\partial^4_x$ given below. The spatial derivatives are discretized using fourth-order formulas, ``DF4":
	\begin{equation*}\label{Diff1-4}
		(\partial_x U)_i=\dfrac{1}{12\Delta x}(-U_{i+2}+8U_{i+1}-8U_{i-1}+U_{i-2}),
	\end{equation*}
	\begin{equation*}\label{Diff2-4}
		(\partial_x^2 U)_i=\dfrac{1}{12\Delta x^2}(-U_{i+2}+16U_{i+1}-30U_{i}+16U_{i-1}-U_{i-2}),
	\end{equation*}
	
	\begin{equation*}\label{Diff4-4}
		(\partial_x^4 U)_i=\dfrac{1}{6\Delta x^4}(-U_{i+3}+12U_{i+2}-39U_{i+1}+56U_i-39U_{i-1}+12U_{i-2}-U_{i-3}).
	\end{equation*}
	
	A standard extension to fourth-order classical Runge-Kutta ``RK4" scheme is used, and thus one obtains the ``DF4-RK4" scheme.
	\begin{remark}\label{remf5}
		At this stage, it is worth mentioning that the numerical scheme of the eB system with high order derivatives~\eqref{eq:eGNLW3} is similar to the one developed for model~\eqref{eq:eGNLW3f5}. In fact, the hyperbolic part of the system is the same as in~\eqref{hyp}, but the high order derivatives involved in the second equation of the remaining (dispersive) part should  be treated accordingly. More precisely, third and fifth order derivatives are discretized  using the following fourth-order formulas:
		\begin{equation*}\label{Diff3-4}
			(\partial_x^3 U)_i=\dfrac{1}{8\Delta x^3}(-U_{i+3}+8U_{i+2}-13U_{i+1}+13U_{i-1}-8U_{i-2}+U_{i-3}),
		\end{equation*}
		\begin{equation*}\label{Diff5-4}
			(\partial_x^5 U)_i=\dfrac{1}{6\Delta x^5}(-U_{i+4}+9U_{i+3}-26U_{i+2}+29U_{i+1}-29U_{i-1}+26U_{i-2}-9U_{i-3}+U_{i-4}).
		\end{equation*}
	\end{remark}
Boundary conditions are imposed using the method presented in Section~\ref{BCsec}. %
	\begin{remark}
An investigation of the dispersive properties of the splitting numerical scheme adopted here is done in~\cite{BCLMT} for a Green-Naghdi type model. The splitting in time of the hyperbolic and dispersive parts is the main originality of this approach. Considering the semi-discretized (in time) version of the splitting scheme adopted in this paper, we believe that the corresponding dispersion relation can be determined following the same classical steps as in~\cite{BCLMT}. An extension to the fully-discretized scheme is possible but very technical and would not lead to any important insight on the dispersive properties of the hyperbolic/dispersive splitting.  To avoid repetition, we do not include in this paper the analysis of the discrete dispersion relation, the interested reader is referred to~\cite[Section 3.4.2]{BCLMT}.
\end{remark}
	
	\subsection{ Boundary conditions}\label{BCsec}
	To close the differential problems, boundary conditions need to be imposed. Boundary conditions for both the hyperbolic and dispersive parts of the splitting scheme are treated by imposing suitable relations on both cell-averaged and nodal quantities. In this paper, we only consider periodic boundary conditions as it was already done in \cite{BGL17} for the study of  internal waves.
	
	For the hyperbolic part, ``ghosts" cells are introduced respectively at the western and eastern boundaries of the domain. The imposed relations on the cell-averaged quantities  are the following:
	\\
	$\bullet$ $\overline U_{-k+1}=\overline U_{N-k+1}$, and $\overline U_{N+k}=\overline U_{k}$, $k \geq 1$, for periodic conditions on western and eastern boundaries.
	
	For the dispersive part, we simply impose the boundary conditions on the nodal values located outside of the domain. In this way, we maintain centered formula at the boundaries, while keeping a regular structure in the discretized model:
	
	$\bullet$ $U_{-k+1}= U_{N-k+1}$, and $U_{N+k}=U_{k}$, $k \geq 1$, for periodic conditions on western and eastern boundaries.
	\red{
\begin{remark}	
As we have already said, in the proposed model we do not try to totally avoid the computations of high order derivatives as in~\cite{Favrie2017,Bassi2020,Escalante2020}, instead we avoid the direct computation of high order derivatives on the flow variable $\zeta$. There are still fourth order derivatives (but this time not directly on $\zeta$) and a large stencil is still needed though. Approximating high order derivatives is one drawback of the proposed numerical scheme due to the large stencil needed and the development of high-order schemes is not an easy task. In the dispersive part of our model, the spatial derivatives are discretized using fourth-order formulas and one disadvantage is that the treatment of boundary conditions is more complex in that case. Moreover, in the case of variable topography, several challenging situations arise involving the design of robust, well-balanced and positive preserving numerical schemes for high-order partial differential equations (PDEs).
	\end{remark}}
	

	\section{Numerical validations}
	This part is devoted to the numerical validations of the model and the numerical scheme. 
	
	We begin by examining the numerical solution of the eB model for the case of solitary waves and \red{show that the latter enjoys better approximate solution with respect to lower-order models.} 
	
	Secondly, we study the propagation of a solitary wave solution with correctors of order $\OO(\eps^3)$ established in~\cite{KLIG}. We compare our numerical solution with an analytic one (up to an $\OO(\eps^3)$ remainder) at several times and
show that our numerical scheme is very efficient and accurate. 

As usual when dealing with a model in oceanographical science, one has to test the ability of the model and the numerical method for the test of a head on collision of 
counter propagating solitary waves. A very good agreement is observed.

	An important fact to reveal is whether or not the improved model is pertinent. 
The third and fourth numerical tests reveal that in presence of large wave number, the choice of the parameter $\alpha$ and the high order derivatives factorization are crucial.
In the fifth test, we study the dam-break problem supplemented by a comparison between the standard and extended Boussinesq models. %
This test is build to test the ability of the eB model and the numerical method to deal with irregular solution. A very good behavior is observed. \red{In the last numerical test, we consider the well-known ``Favre waves" resulting after the impact of a wave on a vertical wall. The numerical results obtained by the eB model are compared with experimental results of Favre~\cite{Favre35} and Treske~\cite{Treske94} and a good agreement is observed.}

In all the test cases, we use the ``WENO5-DF4-RK4" discretization and a CFL number of 1. Our numerical investigations highlight that the time step restriction from the CFL condition of the finite-volume step is enough to ensure stability for the whole numerical method. 

\subsection{Numerical solitary wave solution of the extended Boussinesq system}\label{NumSWSec}
Solitary waves consist of steadily translating disturbances where the nonlinear and dispersive effects counterbalance each other to create a permanent-form with a single crest solution. In this section, we numerically study solitary waves solutions of the eB model and compare it to other models including the accurate full solution of the water-waves model~\cite{DC14,Tanaka86} which will be considered as an ``exact" solution in our comparisons. We will show that the eB model enjoys a better approximate solitary wave solution (when compared to the ``exact" solution) with respect to other lower-order models.

In the last decades, there have been several works on nonlinear PDEs modeling solitary waves. The famous Korteweg-de Vries (KdV) scalar equation or the coupled Boussinesq and Green-Naghdi evolution equations describe the shallow water waves and admit explicit families of solitary wave solutions~\cite{Boussi1872,KDV,Serre53,Chen98}. The calculation of the exact analytic expressions for solitary wave solutions can be done in many ways, one of which, is the direct integration method. This method seeks traveling wave solutions. The PDE is replaced by an ordinary differential equation (ODE) by working in a traveling frame of reference. Hence, one looks for closed-form solutions in terms of special functions.

In~\cite{KLIG}, a careful examination of the direct integration method revealed that the third order non linear ODE associated to the eB model does not admit an explicit analytic solution. The explicit solution of the eB system remains an open problem. To this effect, we employ the Matlab solver \texttt{ode45} to compute numerically the solution of the ODE. The obtained numerical solutions are compared with the accurate full solution of the water-waves model computed using the Matlab script of Clamond and Dutykh~\cite{CD2013} where a fast and precise approach for computing solitary waves solution is introduced.
The fast and accurate Matlab script in~\cite{CD2013} is limited to relatively small velocities. To this end, three values of velocity are used, namely $c=1.025, \ c=1.01$ and $c=1.002$. We compare the obtained solutions with other models including the accurate full solution of the water-waves model, the original Green-Naghdi system ($\zeta_{GN}$) and the Boussinesq system ($\zeta_{B}$). The explicit solution of the original Green-Naghi model~\cite{Serre53,SuGardner69} reads:
\begin{equation}\label{solGN}
 \eps \zeta_{GN}(x)=(c^2-1)\  \text{sech}^2 \Big(\sqrt{\dfrac{3(c^2-1)}{4c^2\eps }}\ x\Big)= \eps c^2 \zeta_{B}(x)\;.
\end{equation}
The waves are re-scaled so that the Boussinesq solution do not depend on $c$. Consistently, we set $\eps=1$.

\begin{figure}[H]
	\centering
	\subcaptionbox{Re-sized waves, $c=1.025,\ 1.01,\ 1.002$}
	{\includegraphics[width=0.47\textwidth]{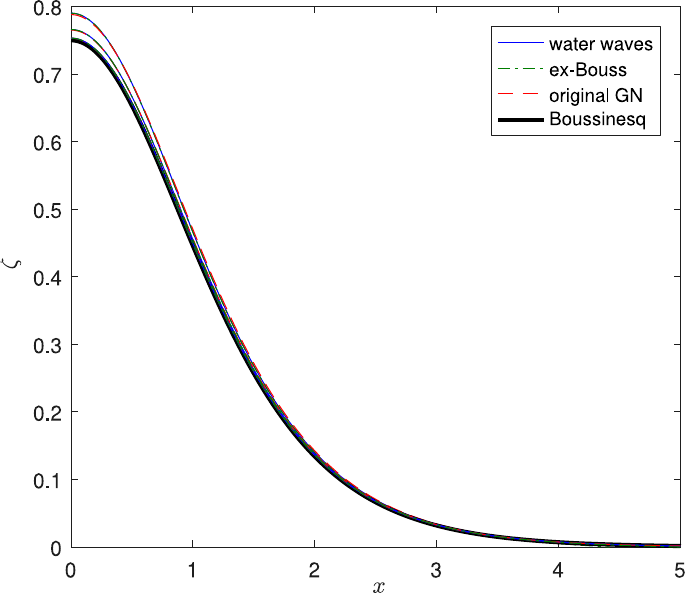}}
	\subcaptionbox{Zoom in}
	{\includegraphics[width=0.48\textwidth]{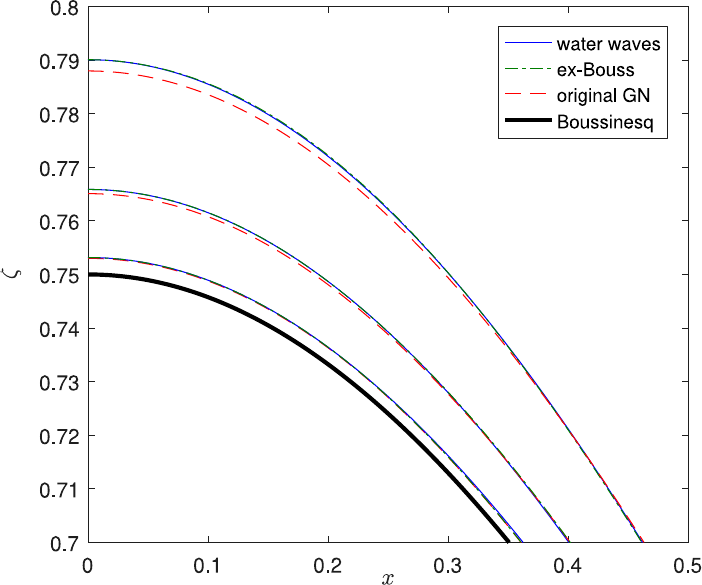}}
	\captionsetup{justification=centering}
	\caption{Comparison of the solitary waves solutions.}
	\label{SWcomp}
\end{figure}

After re-scaling, Figure~\ref{SWcomp} shows clearly that the solitary waves tend towards the Boussinesq solution $(\zeta_{B})$ as $c-1 \rightarrow 0$. A zoom-in on the crest of half of the solitary waves shows that, of all models, the eB model has the best match with the full Euler system (water-waves) solution.
\begin{table}[h]
\centering
\caption{The normalized $l^2$-norm of the error for the Boussinesq, original GN and eB models.}
 \label{l2norm}
\begin{tabular}{@{}c|cc|cc|cc@{}}
\toprule
&  Boussinesq   &  & original GN &  & ex-Bouss \\
\toprule
$c-1$ & Error & Conv. rate & Error & Conv. rate  & Error & Conv. rate\\
\midrule
0.025  &$0.0459$ & -- & $0.0058$ & --& $0.0019$ & -- \\\
0.01& 0.0188 & 0.9743& 0.0023 & 1.0361& $3.07\times 10^{-4}$ & 1.9984\\
0.002& 0.0038 & 0.9864 & $4.45\times 10^{-4}$ & 1.0186& $1.23\times 10^{-5}$ & 1.9990\\
\bottomrule
\end{tabular}
\end{table}
\begin{figure}[H]
\centering
\includegraphics[scale=0.6]{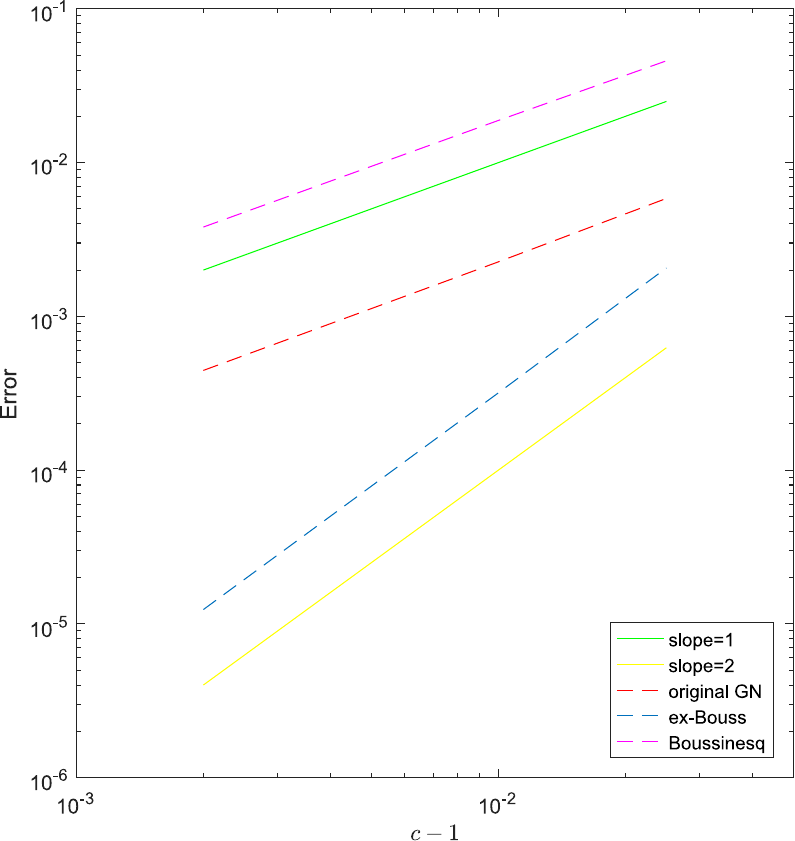}
\caption{Errors as a function of
$c-1$ (log-log plot).}
\label{convratefig}
\end{figure}
The normalized $l^2$-norm of the error as a function of $c-1$ is plotted in a log-log scale in Figure~\ref{convratefig}. The convergence rate is indeed quadratic for the extended Boussinesq model whereas it is only linear for the original Green-Naghdi model. Results are gathered in Table~\ref{l2norm}. \red{This highlights the fact that the higher-order extended Boussinesq model have a better approximate solution when compared to lower-order models.} 

	\subsection{Propagation of a solitary wave solution with correctors}\label{PSWCsec}
\begin{figure}[H]
\begin{center}
     \includegraphics[width=1.0\textwidth]{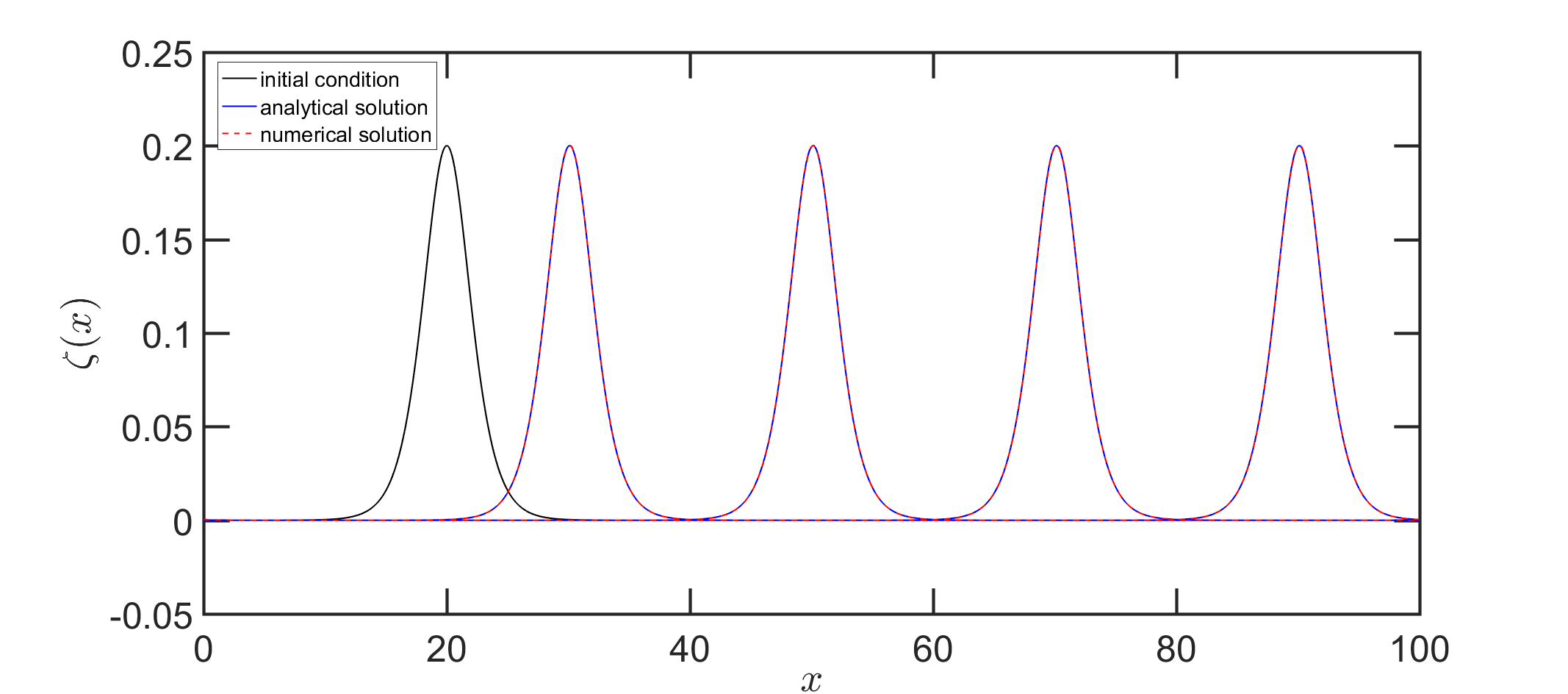}
  \caption{Propagation of a solitary wave: water surface profiles at $t=0,10,30,50$ and $70 $.}
\label{figcompSWP}
\end{center}
\end{figure}

\begin{figure}[H]
\begin{center}
     \includegraphics[width=1\textwidth]{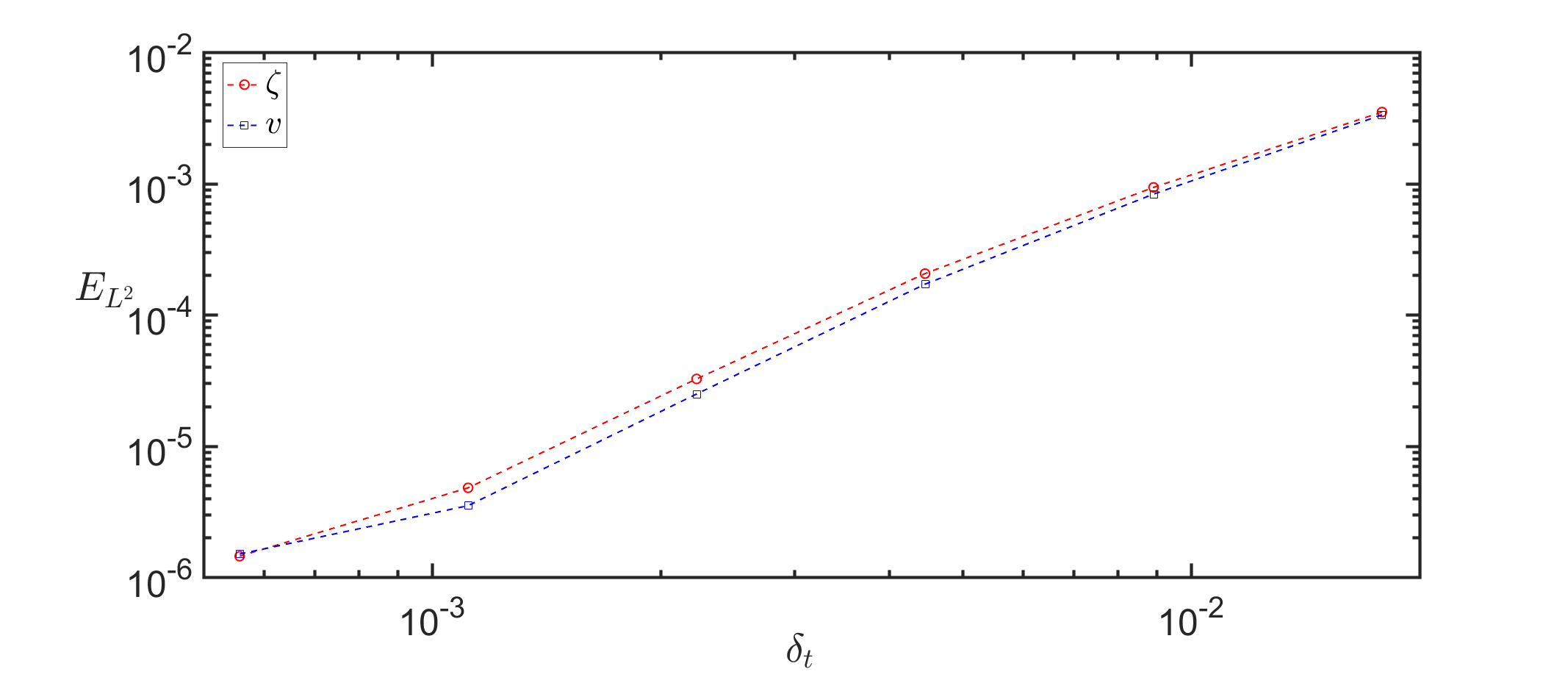}
  \caption{Propagation of a solitary wave: $L^2$-error on the water surface deformation and the averaged velocity for $a=0.2$.}
\label{loglogplot}
\end{center}
\end{figure}

A careful examination revealed that the extended Boussinesq system~\eqref{ex-boussinesq} does not admit an exact solitary wave solution, see~\cite[Section 4]{KLIG}. In order to validate our numerical scheme we use the explicit solution with correctors of order $\mathcal{O}(\eps^3)$ found in~\cite[Section 5]{KLIG} that we disclose below. Such solitary waves are analytical solutions of the extended Boussinesq system~\eqref{ex-boussinesq} up to $\mathcal{O}(\eps^3)$ remainders. Therefore, this family of solutions can be used as a validation tool for our present numerical scheme and its given by ($\zeta,v)$ with
\begin{equation}\label{zetadef}
\zeta=\zeta_1 + \dfrac{\eps^2}{ 2 }\Big[(\zeta_2^0+v_2^0)(x-t)+(\zeta_2^0-v_2^0)(x+t)+ \int_0^t f(s,x-t+s) ds-\int_0^t f(s,x+t-s) ds\Big],
\end{equation}
and
\begin{equation}\label{vdef}
v= v_1+\dfrac{\eps^2}{2}\Big[(\zeta_2^0+v_2^0)(x-t)-(\zeta_2^0-v_2^0)(x+t)+ \int_0^t f(s,x-t+s) ds+ \int_0^t f(s,x+t-s) ds \Big],
\end{equation}
where $(\zeta_1,v_1)$ is the well known explicit solution of solitary traveling wave of the sB system~\eqref{standard-bouss} given by:
\begin{equation*}\label{stdbousssol}
\left\{
\begin{array}{lcl}
\displaystyle\zeta_{1}(t,x)= a \ \text{sech}^2 \Big(k \ (x-ct)\Big)\vspace{1mm}\; ,\\
\displaystyle v_{1}(t,x)=\dfrac{ c \zeta_1(t,x)}{1+\eps\zeta_1(t,x)}\; ,
\end{array}
\right.
\end{equation*}
where $k=\sqrt{\dfrac{3a}{4}}$ and $c=\sqrt{\dfrac{1}{1-a\eps}}$ and $a$ is an arbitrary chosen constant. The initial conditions $\zeta_2^0$ and $v_2^0$ are both given in $C^\infty(\mathbb{R})$ and set $\zeta_2^0=v_2^0=\exp\Big(-\Big(\dfrac{3\pi x}{10}\Big)^2 \Big)$. The function $f(t,x)$ is defined by:
 $$f(\zeta_1,v_1)=\partial_x \zeta_1 \partial_x \partial_t v_1 +\dfrac{2}{3}\zeta_1 \partial_x^2 \partial_t v_1 + \dfrac{1}{45}\partial_x^4 \partial_t v_1  +\dfrac{1}{3} \partial_x\big(v_1 (v_1)_{xx}-(v_1)_x^2\big).$$ 
In this test, we investigate the left to right propagation of a solitary wave initially centered at $x_0 = 20 $, of amplitude $a = 0.2 $. The computational domain length is $L=100 $ and discretized with 1600 cells. The solitary wave is initially far from boundaries, thus the periodic boundary conditions do not affect the computation. The water surface profile of our numerical solution provided by the model~\eqref{eq:eGNLW3f5} with $\alpha=1$, is compared with the analytical one given by~\eqref{zetadef}-\eqref{vdef} at several times using the fifth order discretization ``WENO5-DF4-RK4". An excellent agreement between numerical and analytical solutions is observed in Figure~\ref{figcompSWP}.  The amplitude and shape of the computed solitary wave are accurately preserved during the propagation, indicating an accurate discretization of the governing equations in both space and time.
To complete the picture and assess the convergence of our numerical scheme, we compute the numerical solution for this particular test case for an increasing number of cells $N$, over a duration $T=1 $. We start with $N=400$ number of cells ($\delta_x=\frac{L}{N} =0.25 $) and successively multiply the number of cells by two. The relative errors  $E_{L^2}(\zeta)$ and $E_{L^2}(v)$ on the water surface deformation and the averaged velocity are computed at $t= 1$, using the discrete $L^2$ norm $\|.\|_{2}$:
\begin{equation}
E_{L^2}(\zeta)=\dfrac{\|\zeta_{num}-\zeta_{sol}\|_{2}}{\|\zeta_{sol}\|_{2}}; \qquad E_{L^2}(v)=\dfrac{\|v_{num}-v_{sol}\|_{2}}{\|v_{sol}\|_{2}},
\end{equation}
where $(\zeta_{num}, v_{num})$ are the numerical solutions and $(\zeta_{sol}, v_{sol})$ are the analytical ones coming from~\eqref{zetadef}-\eqref{vdef}. Results are presented in Table~\ref{L2err} and Figure~\ref{loglogplot} where $E_{L^2}(\zeta)$ and $E_{L^2}(v)$ are plotted against $\delta_x$ 
in log scales, for the considered relative amplitude $a=0.2  $.
Very accurate results are obtained, indicating that the employed numerical method is capable of computing in a stable way the propagation of a solitary wave. Moreover, computing a linear regression on all points yields a slope equal to $2.33$ for $\zeta$ and $2.34$ for the averaged velocity $v$. This result sounds rational because the global (time and space) order of our scheme  may be limited by the order of the splitting method used here, which is of order two as already discussed by Bonneton \textit{et al.} in~\cite{BCLMT}.

\begin{table}[h]
\centering
\caption{Propagation of a solitary wave: relative $L^2$-error table for the conservative variables.}
 \label{L2err}
\begin{tabular}{@{}lllll@{}}
\toprule
N & $E_{L^2}(\zeta)$ & Conv. rate &$E_{L^2}(v)$ & Conv. rate\\
\midrule
400 & $3.50\times 10^{-3}$ &-- & $3.33\times 10^{-3}$ & --\\
 800& $9.32\times 10^{-4}$ &1.9070& $8.29\times 10^{-4}$ &2.0064 \\
 1600 & $2.05\times 10^{-4}$&2.0466&  $1.70\times 10^{-4}$ &2.1433  \\
3200 & $3.23\times 10^{-5}$ &2.2449& $2.48\times 10^{-5}$&2.3479  \\
 6400& $4.79\times 10^{-6}$ &2.3873& $3.50\times 10^{-6}$&2.4845\\
  12800& $1.44\times 10^{-6}$ &2.3344& $1.49\times 10^{-6}$&2.3436\\
\bottomrule
\end{tabular}
\end{table}

\subsection{Head on collision of counter propagating solitary waves}
\begin{figure}[H]
	\centering
	{\includegraphics[scale=0.35]{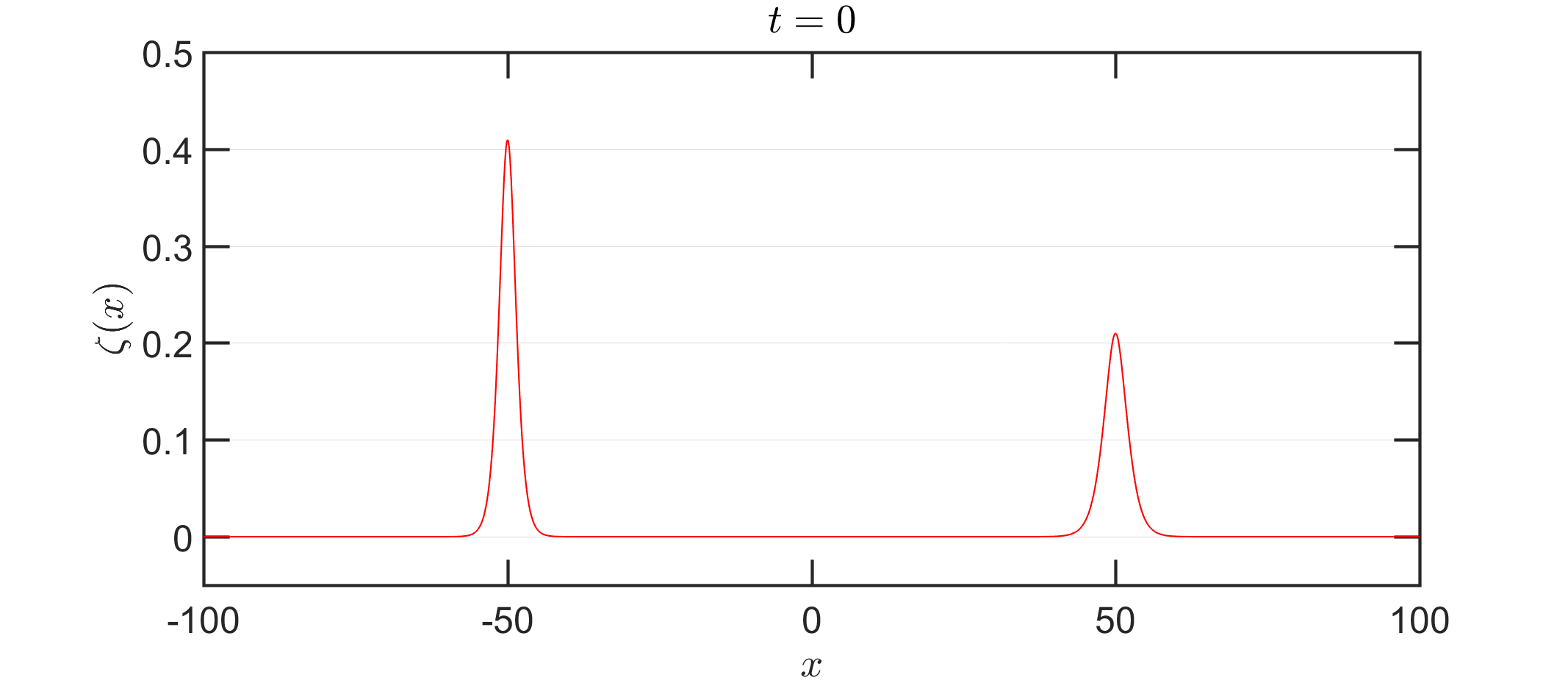}}
	{\includegraphics[scale=0.35]{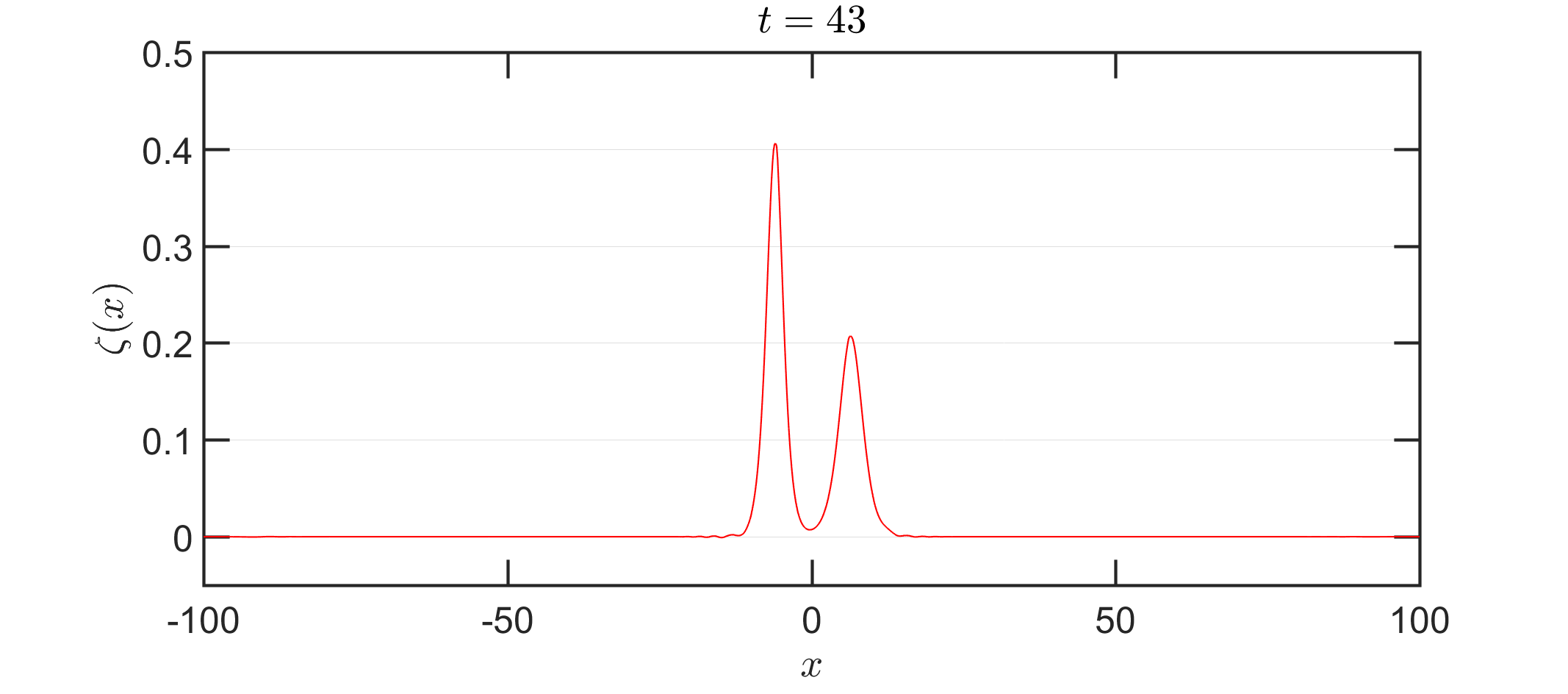}}
				{\includegraphics[scale=0.35]{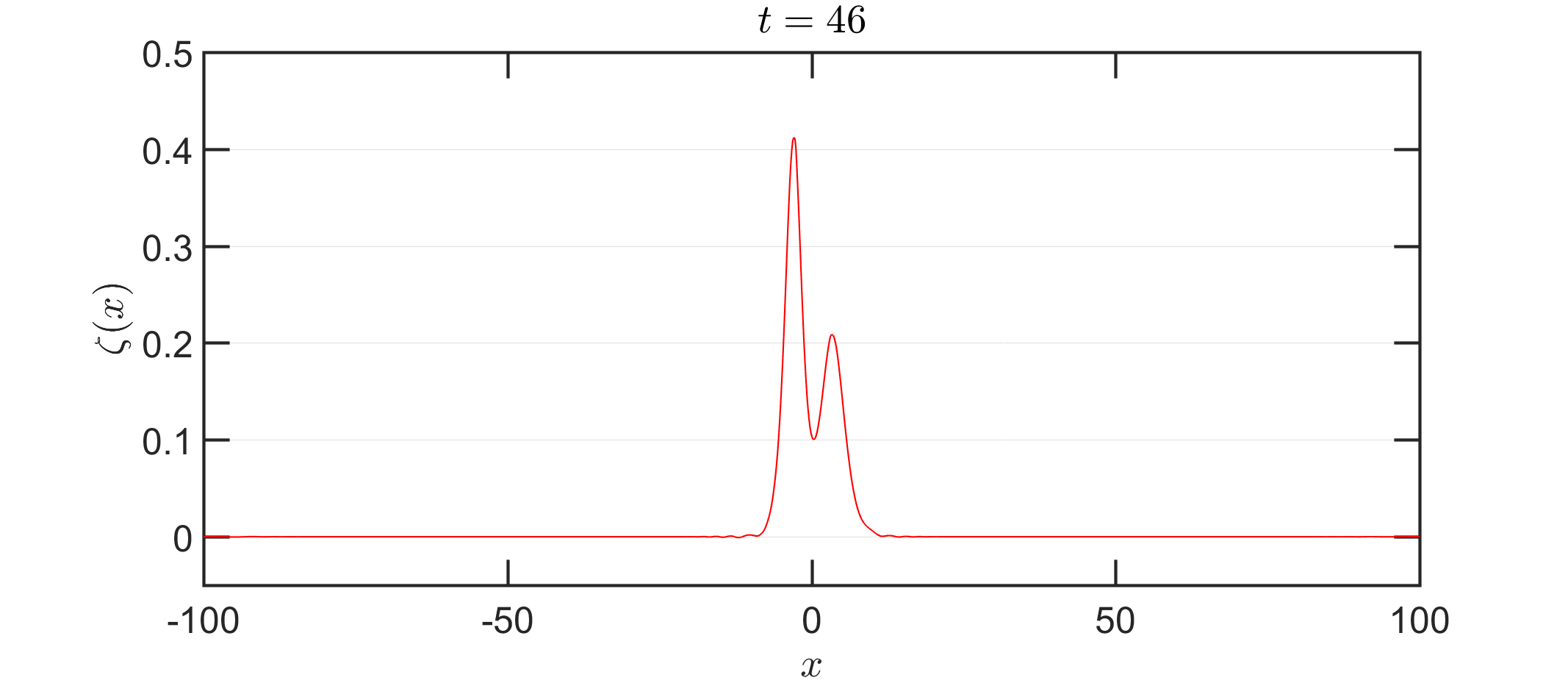}}
		{\includegraphics[scale=0.35]{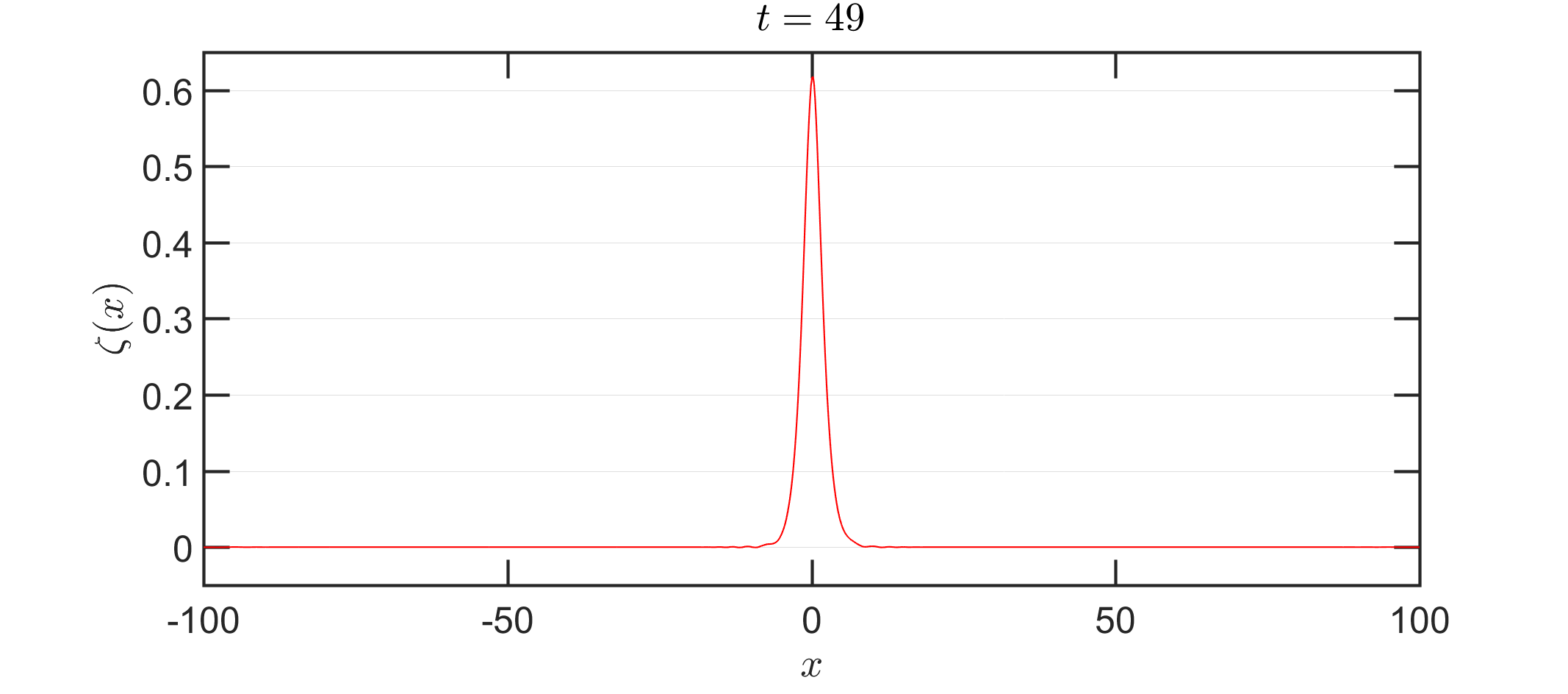}}
	\captionsetup{justification=centering}
	\caption{Head on collisions: surface wave shape at $t=0, 43, 46$ and $t=49$.}
	\label{HOCSW}
\end{figure}

\begin{figure}[H]
		\centering
	{\includegraphics[scale=0.35]{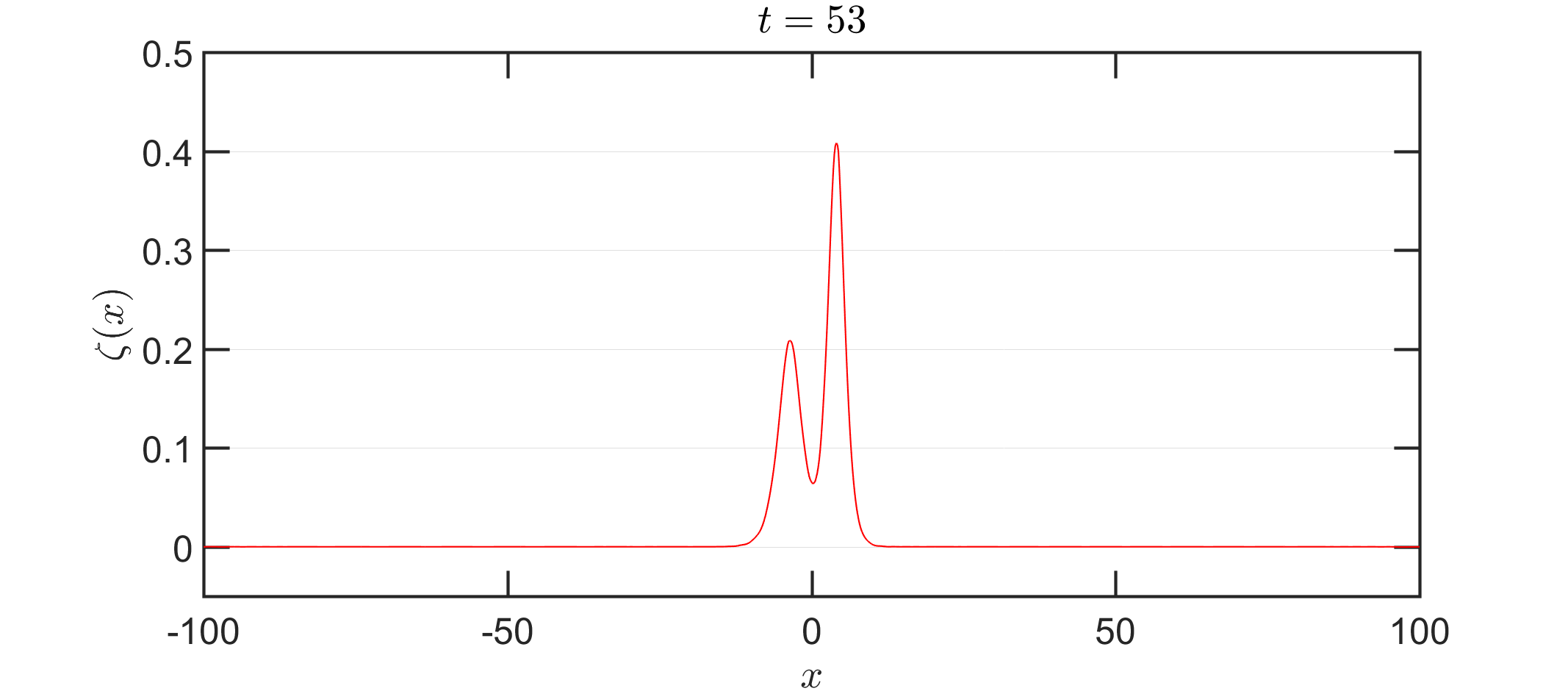}}
		{\includegraphics[scale=0.35]{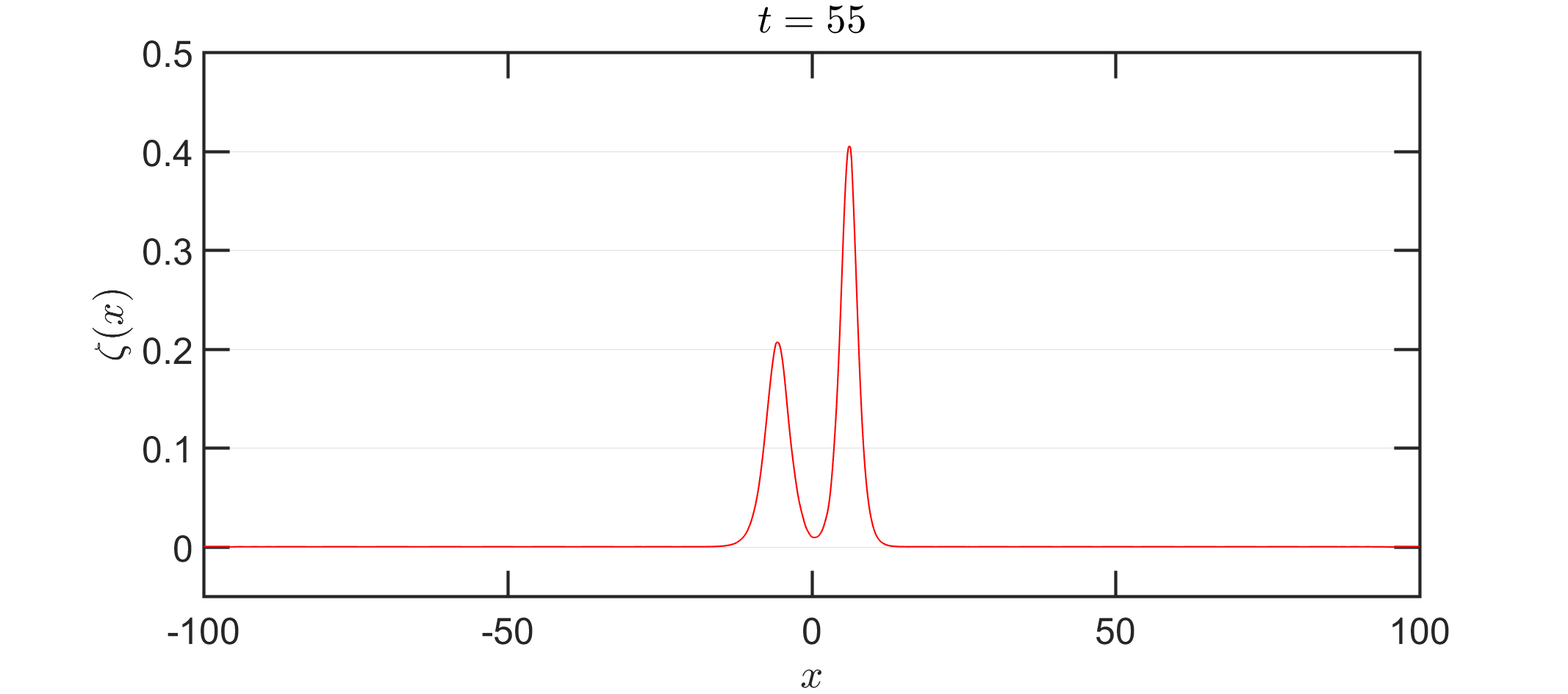}}
	{\includegraphics[scale=0.35]{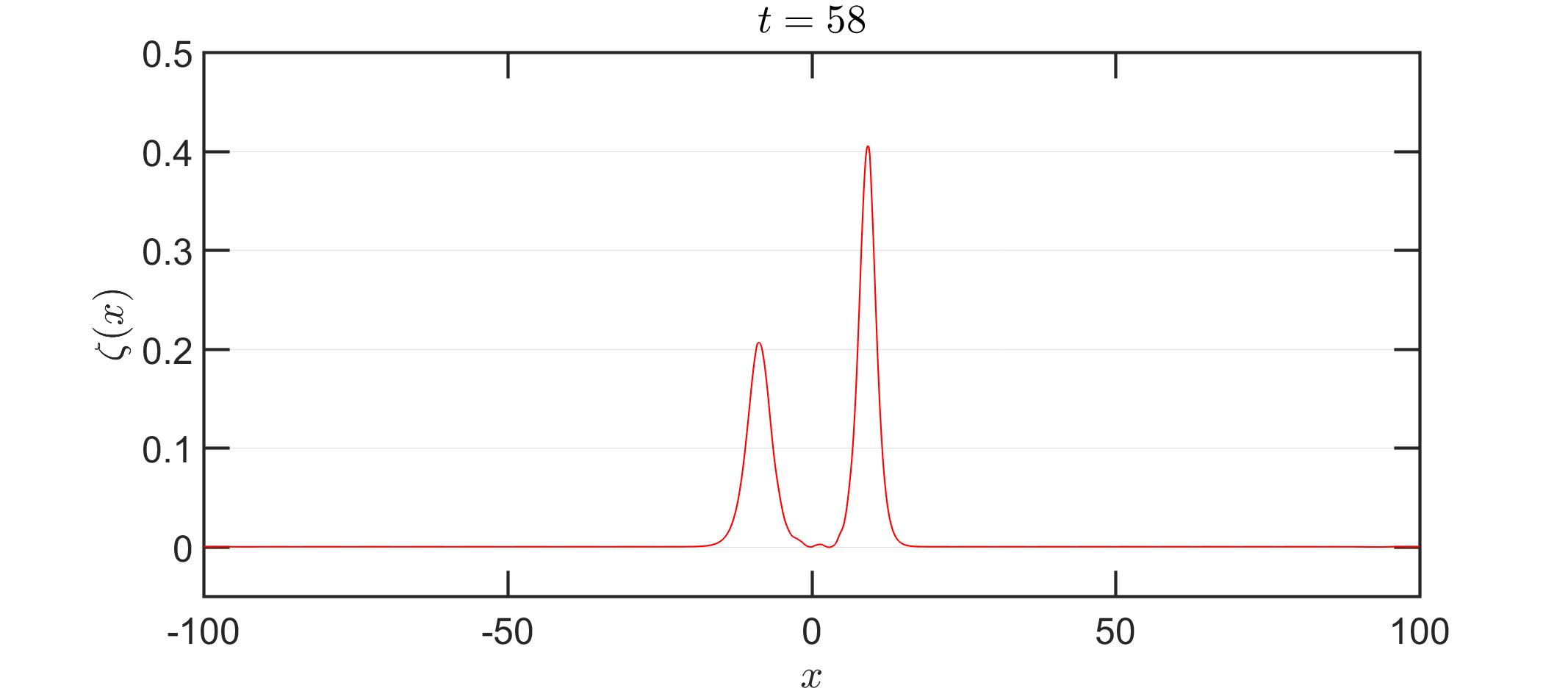}}
	{\includegraphics[scale=0.35]{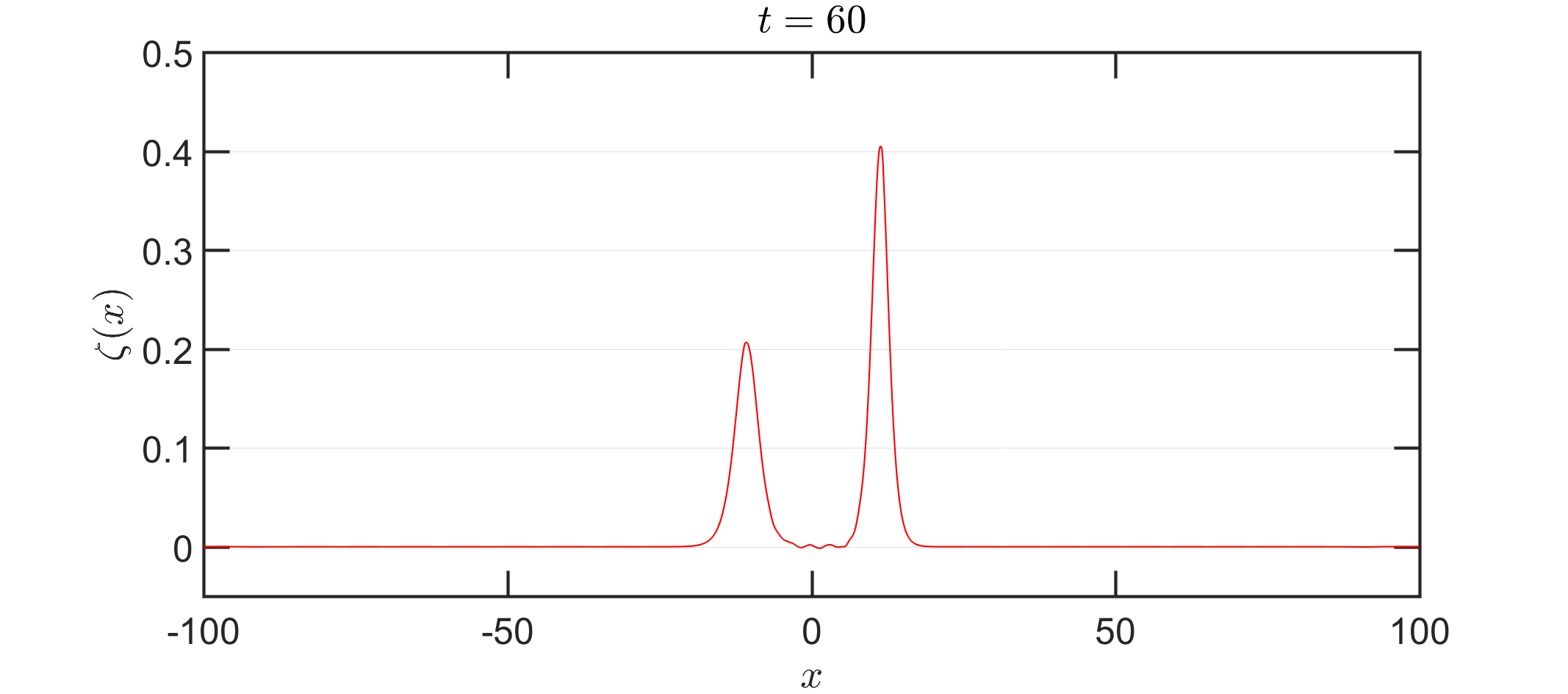}}
	\captionsetup{justification=centering}
	\caption{Head on collisions: surface wave shape at $t=53, 55, 58$ and $60$.}
	\label{HOCSW1}
\end{figure}

\begin{figure}[H]
		\centering
		{\includegraphics[scale=0.35]{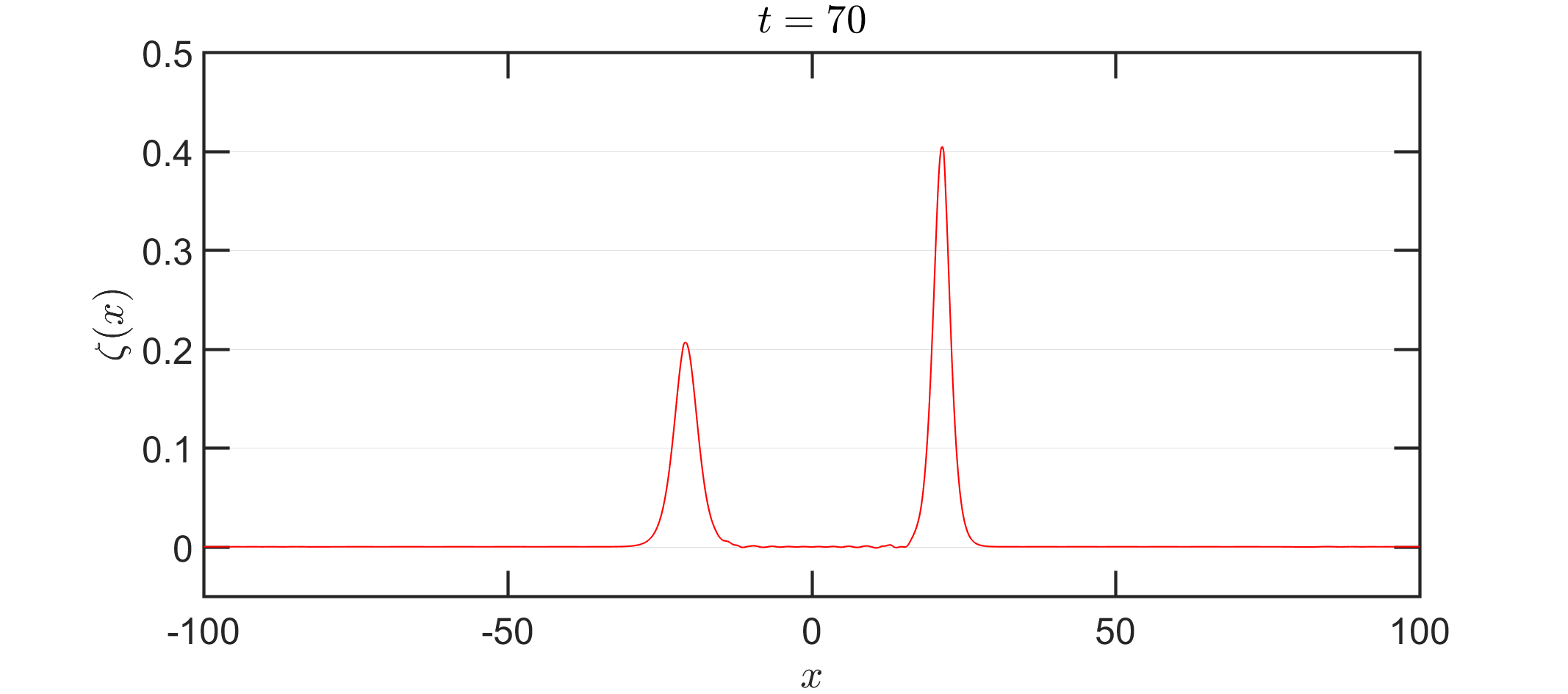}}
	{\includegraphics[scale=0.35]{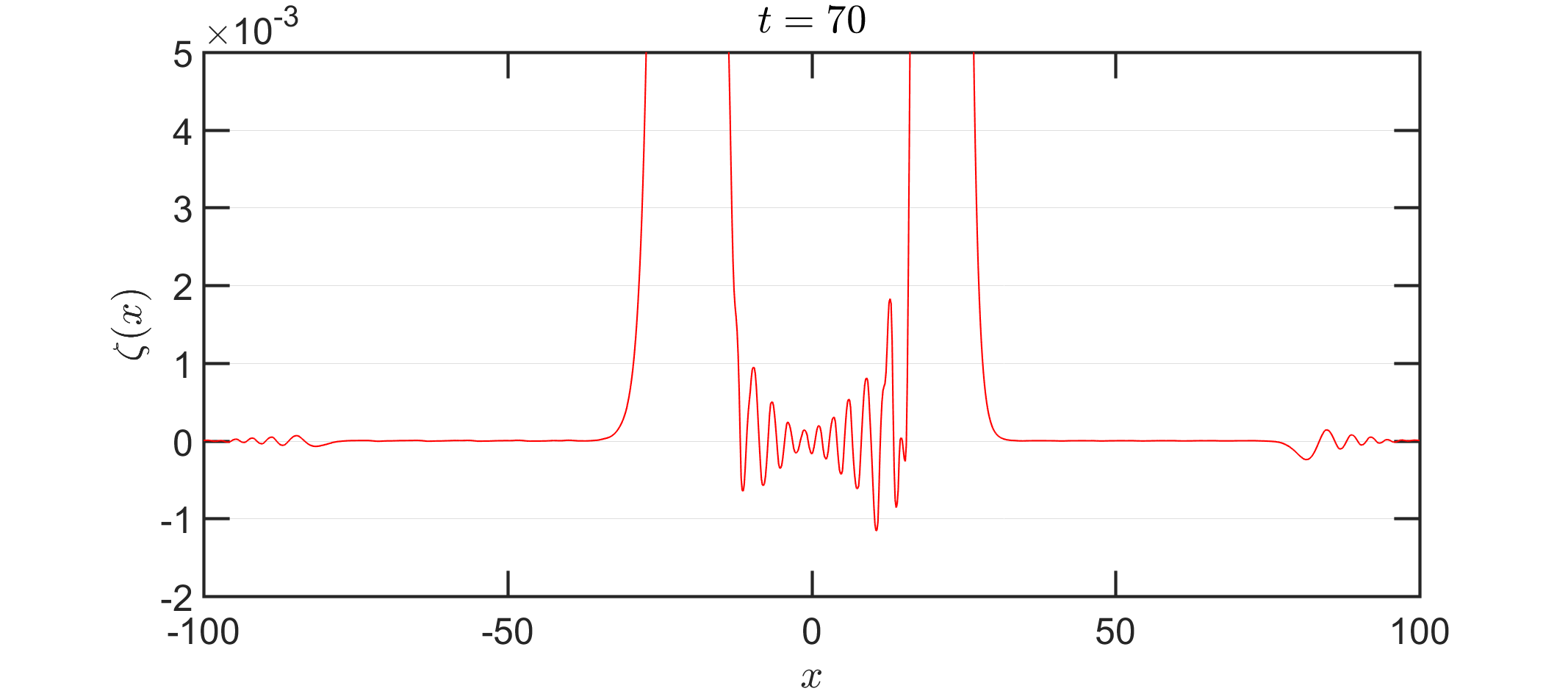}}
	\captionsetup{justification=centering}
	\caption{Head on collisions: surface wave shape at $t=70$.}
	\label{HOCSW2}
\end{figure}

	A standard nonlinear test case for numerical methods is the interaction of solitary wave. In this numerical test, we study an important phenomenon in the study of nonlinear dispesive waves, the head on collision of two counter propagating waves with different amplitudes. We used solitary wave solutions with correctors of order $\mathcal{O}(\eps^3)$ for the eB system established in~\cite{KLIG} and defined in~\eqref{zetadef}-\eqref{vdef}. We consider two solitary waves centered at $x=-50 $ and $x=50 $ at $t=0 $ on a spatial domain $L=200$ with a constant depth
$h_0 =2 $, see Figure~\ref{HOCSW}. The solitary wave centered at $x=-50$ travels to the right with a speed $c_{s,1}=1.0206$ and an initial amplitude $a_1=0.4$ while the one centered at $x=50$ travels to the left with a speed $c_{s,2}=1.0102$ and an initial amplitude $a_2=0.2$. The domain is discretized using 1200 cells and periodic boundary conditions are imposed. The numerical solutions are computed using model~\eqref{eq:eGNLW3f5} with $\alpha=1$. The collision of the two waves starts at about $t = 43$, see Figure~\ref{HOCSW}. After the interaction, each wave continue moving in its own direction and turn up to be unaffected by the collision, see Figure~\ref{HOCSW1}. A proper description of the distinctive nature of nonlinear interactions is illustrated when zooming at the oscillating dispersive tails of very small amplitude appearing at the center of the domain at $t = 70$ in Figure~\ref{HOCSW2}. One can also observe two dispersive tails with smaller amplitudes located to the left and right boundaries. The generation of such dispersive tails is due to the $\OO(\eps^3)$ remainder terms as mentioned in the beginning of Section~\ref{PSWCsec}. The high precision of our numerical scheme is verified after accurately capturing this phenomenon and inducing similar observations to earlier works~\cite{EGS06,MGH10} where the head-on collision is carried out.

	\subsection{Breaking of a regular heap of water with a large wave number}\label{SecNumTst1}
	In this numerical test, we highlight the importance of factorizing high order derivatives present in the improved eB model~\eqref{eq:eGNLW3f5} together with the appropriate choice of the parameter $\alpha$ in improving the frequency dispersion in high frequency regimes. To this end, we consider a sufficiently regular heap of water with a large wave number represented by the initial data:
	\begin{equation*}
		\zeta(0,x)=0.7 e^{-80x^2}, \ v(0,x)=0,
	\end{equation*} 
	(dashed lines) with a domain of computation $x \in (-2,2)$ discretized with 512 cells and under periodic boundary conditions.
	The non-linearity parameter is set as follow: $\eps=0.1$ (non-dimensional setting). Our numerical solutions are computed using models~\eqref{eq:eGNLW3} (without factorization) and~\eqref{eq:eGNLW3f5} (with factorization).
	We compare our numerical solutions with the numerical solutions computed using the Matlab script of Duch\^{e}ne, Israwi and Talhouk~\cite{DucheneIsrawiTalhouk16}
	and with the lower-order GN-CH model obtained in~\cite{BGL17}. In~\cite{DucheneIsrawiTalhouk16}, the original Green-Naghdi (GN) model describing a two-layer flow~\footnote{One can easily recover the one-layer configuration by setting $\gamma=0$ and $\delta=1$.}\label{1layerFN}  is improved in terms of frequency dispersion by introducing a new class of tailored GN models with a slight modification of the dispersion components using a class of Fourier multipliers. In particular, an ``improved" GN model is derived sharing the same dispersion relation as the full Euler (FE) system. These type of models are commonly called full dispersion models~\cite{BonaLannesSaut08}. The ``improved" GN model is a full dispersion model and its numerical solution will be used as a reference solution. Eventually, the goal of this test case is to show that the optimized eB model~\eqref{eq:eGNLW3f5} may provide some results which are in the same league of full dispersion models.  
	\begin{figure}[H]
	\centering
	\subcaptionbox{Comparison of the numerical solutions of the eB model~\eqref{eq:eGNLW3f5} (blue) with the ``improved" GN model (green) and the GN-CH model (red).}
	{\includegraphics[scale=0.35]{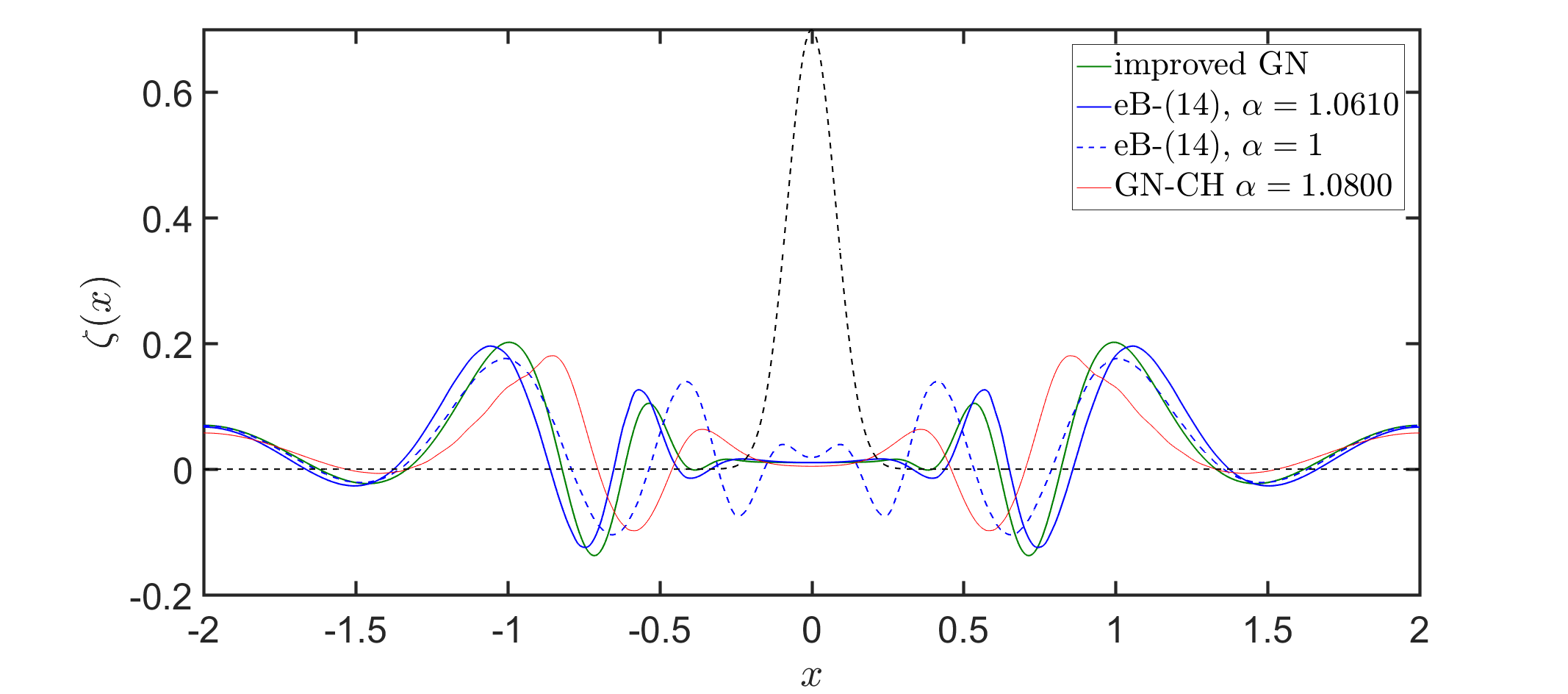}}
	\subcaptionbox{Comparison of the numerical solutions of the eB model~\eqref{eq:eGNLW3} (blue) with the ``improved" GN model (green) and the GN-CH model (red).}
	{\includegraphics[scale=0.35]{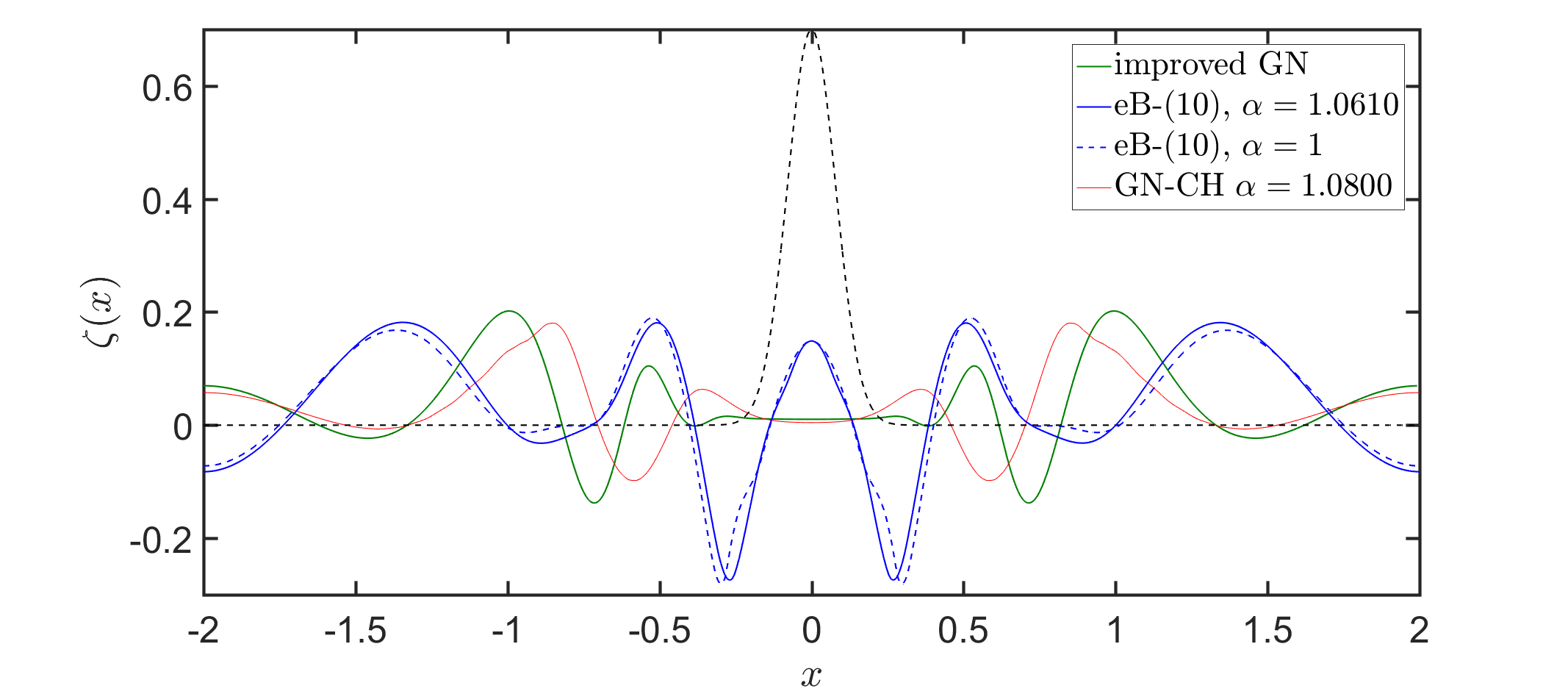}}
	\captionsetup{justification=centering}
	\caption{Comparison of the numerical solutions of the eB model~\eqref{eq:eGNLW3f5} and the eB model~\eqref{eq:eGNLW3} with the ``improved" GN model and the GN-CH model~\cite{BGL17} at $t=3$.}
\label{figcomphf}
\end{figure}
	Figure~\ref{figcomphf}(a) shows when $\alpha$ is chosen appropriately as discussed in section~\ref{Secalphachoice}, namely $\alpha_{opt}=1.0610$, our numerical solution computed over a sufficient duration $t=3$ using model~\eqref{eq:eGNLW3f5} behave similarly to the one computed with the ``improved" GN model (sharing the same dispersion relation as the FE system). In contrary, when choosing $\alpha_{opt}=1$ or when using the model~\eqref{eq:eGNLW3} to compute the solution (Figure~\ref{figcomphf} (b)), the behavior is different than the ``improved" GN model.  Note that the GN-CH model has an improved frequency dispersion due to the careful choice of the parameter $\alpha$. Nevertheless, the numerical solution computed using the GN-CH model is far from the numerical solution of the ``improved" GN model. The observed agreement in high frequency regime between the numerical solutions of the ``improved" GN model and the eB model~\eqref{eq:eGNLW3f5} rather than the GN-CH model is due to the factorized high order dispersion terms existing in~\eqref{eq:eGNLW3f5}.  In fact, the eB model~\eqref{eq:eGNLW3f5} is precise up to $\OO(\eps^3)$ order and thus contains factorized high-order dispersive terms that do not exist in the GN-CH model. As already stated in Section~\ref{Secalphachoice}, one can see that the choice of an optimal value of $\alpha$ when using the model~\eqref{eq:eGNLW3} has no beneficial effect due to the high frequency regime setting. This numerical test confirms that the model~\eqref{eq:eGNLW3} has a range of applicability limited to $k\leq 1$ and thus has poor dispersion properties in intermediate and large wave numbers regime. 
	\subsection{Breaking of a regular heap of water with a small wave number}\label{BRHWswn}
\begin{figure}[H]
		\begin{center}
			\includegraphics[width=1\textwidth]{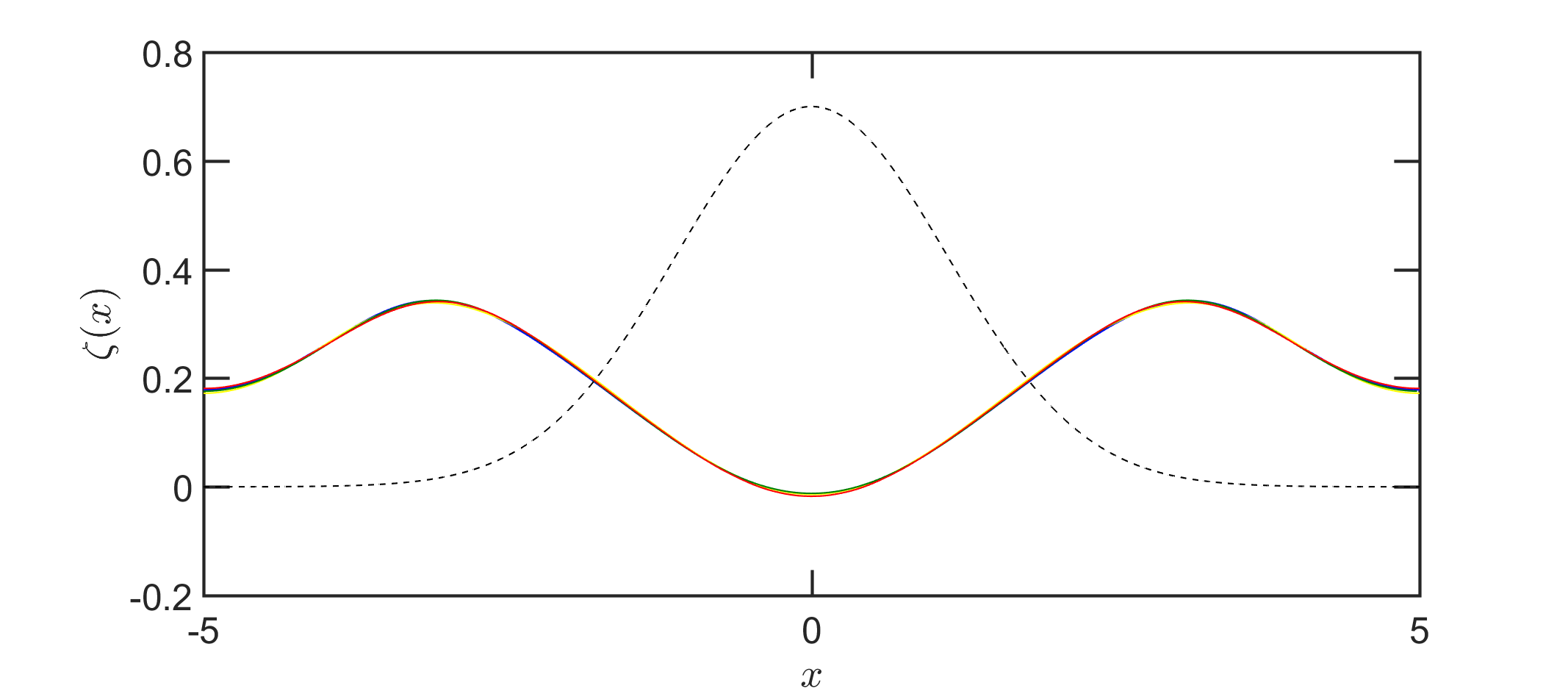}
			\caption{Comparison of the numerical solutions of model~\eqref{eq:eGNLW3f5} (blue line) and model~\eqref{eq:eGNLW3} (yellow line) and the GN-CH model~\cite{BGL17} (red line) with the ``improved" Green-Naghdi model (green line) at $t=3$.}
			\label{figcomplf}
		\end{center}
	\end{figure}

	In this numerical test, we consider the breaking of a sufficiently regular heap of water with a small wave number whose initial data is:
	\begin{equation*}
		\zeta(0,x)=0.7 e^{-0.4x^2}, \ v(0,x)=0,
	\end{equation*} 
	(dashed lines) within a domain of computation $x \in (-2,2)$ discretized with 512 cells and under periodic boundary conditions. The non-linearity parameter is set as follow: $\eps=0.5$ (non-dimensional setting). Our numerical solutions are computed using models~\eqref{eq:eGNLW3} and~\eqref{eq:eGNLW3f5} and compared with the numerical solution of the ``improved" Green-Naghdi model and the numerical solution of the GN-CH model~\cite{BGL17} over a sufficient duration $t=3$. The parameter $\alpha$ is fixed as 1  
	since varying $\alpha$ does not yield significant improvements. %
	 In fact, both eB models~\eqref{eq:eGNLW3} and~\eqref{eq:eGNLW3f5} have an equivalent dispersion relation to the one of the \emph{full Euler} system for small wave numbers 
	 and the choice of $\alpha$ does not play any role in the leading terms.
	 Indeed, Figure~\ref{figcomplf} shows a fairly good agreement between the solutions of the eB models~\eqref{eq:eGNLW3} (yellow line) and~\eqref{eq:eGNLW3f5} (blue line), and the solution of the GN-CH model~\cite{BGL17} (red line) and the one of the ``improved'' Green-Naghdi model (green line). This confirms the fact that, in small wave numbers regime, 
	every 
	 aforementioned model behave similarly and enjoy similar dispersive properties as the one of the \emph{full Euler} system. 
	\subsection{Dam-break problem}
\begin{figure}[H]\
		\centering
		{\includegraphics[width=1.0\textwidth]{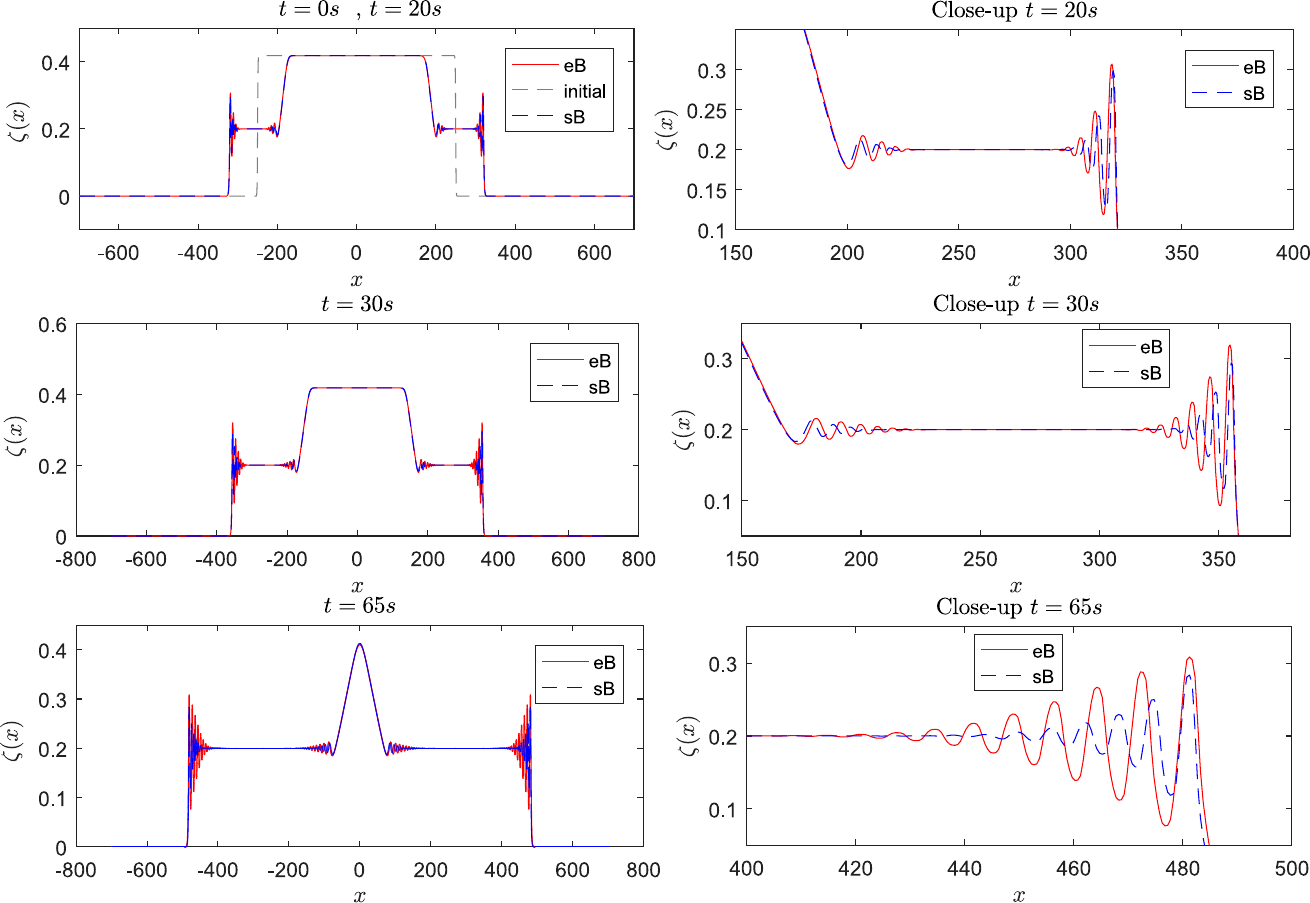}}
		\caption{Dam break: wave shape at different times, comparison between the numerical solution of the eB model (solid red line) and sB model (dashed blue lines).}
		\label{dam1}
	\end{figure}

\begin{figure}[H]
		\centering
		{\includegraphics[width=0.8\textwidth]{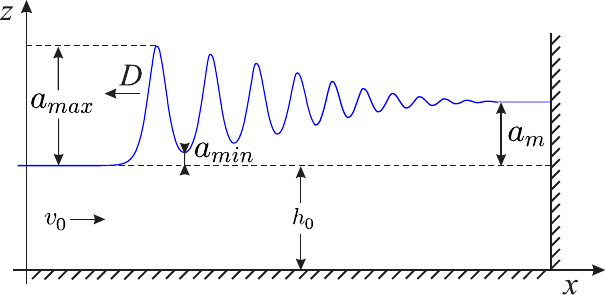}}
		\caption{A sketch of Favre waves.}
		\label{Favre_sktech}
	\end{figure}

	Dealing with non-regular solutions needs a special treatment at the numerical scheme level. Earlier works~\cite{BCLMT,CLM,LannesMarche14} have shown that the use of high-order schemes in dispersive waves study is necessary to prevent the corruption of the dispersive properties of the model. In general, dispersive shock waves are generated due to the dispersive effects~\cite{MGH10,MID14} when considering discontinuous initial data. In this numerical test, we implement a dam break problem in order to investigate the performance of our numerical scheme in handling non-regular solutions. We study the dam-break problem in the extended (eB) and standard (sB) Boussinesq models. We consider the following initial data:
	\begin{equation}\label{daminitialdata}
		\zeta(0,x)=a[1+\tanh(250-|x|)], \quad v(0,x)=0,
	\end{equation} 
	with $a=0.2091 \ m$ defined on the computational domain $x\in (-700,700)$ discretized using 2800 cells and imposed under periodic boundary conditions. 
	The solutions of the eB model are computed using the dimensional version of the eB model with factorized high order derivatives knowing that similar results were obtained when using the model without factorization. To this end, we set $\eps=1$ and add the gravity term to the equations as needed. We make this choice in order to test the scheme in the challenging conditions of the strongly dispersive regime. %
	In this test, the choice of $\alpha$ is not important, thus we choose $\alpha=1$.
The solutions of the sB model are computed accordingly using the same splitting scheme.
		The dam break wave shapes of the sB and eB models are shown at different times in Figure~\ref{dam1}. The initial data break and generate dispersive shock waves. A close-up on the profiles at $t=30s$ and $t = 65 s$ shows two dispersive shock waves counter-propagating on both sides of the ``dam", and two rarefaction waves moving in the direction of the center. The dispersive tail generated by the eB model is larger and have higher amplitude oscillations when compared to the respective tail generated by the sB model. This result is consistent with the fact that the eB model contains high order nonlinear dispersive terms not present in the sB model.
			 This numerical test shows that our high-order numerical scheme was able to capture accurately the dispersive shock waves phenomenon. Dispersive shock waves in a large class of dispersive shallow water models were carried out in several works~\cite{BGL17,EGS06,MGH10,MID14} and show a good agreement with our numerical simulation. 
\begin{remark}
The splitting strategy may limit the whole method to second order accuracy, however, for the study of dispersive waves, it is necessary to use high-order schemes to prevent the corruption of the dispersive properties of the model by some dispersive truncation errors linked to second-order schemes. The reader is referred to~\cite{BGL17}, where the study of dam-break problem in the GN-CH model is supplemented by a comparison between second (MUSCL-RK2) and fifth (WENO5-RK4) order accuracy methods. Although a splitting scheme of order 2 was adopted therein, the dispersive properties were clearly corrupted by the second-order method while the higher-order method was able to deal with discontinuous initial data and capture the rapid oscillations in dispersive shock waves. 
\end{remark}
\red{
\subsection{Favre waves}
In this numerical test, we consider the well-known ``Favre waves" resulting after the impact of a wave on a vertical wall. Experimental investigations of those waves were addressed for the first time by Favre~\cite{Favre35} in a rectangular channel. Similar experiments have been performed by Treske in~\cite{Treske94} and Soares Frazao and Zech in~\cite{Soares02} in open channels. This problem was numerically studied more recently in~\cite{Soares08,Tissier11}. Due to dispersion, the uniform free surface flow impacting a wall reflects and free surface undulations (called ``Favre waves", see Figure~\ref{Favre_sktech}) appear. The leading wave has a maximum amplitude $a_{max}$ and a minimum amplitude $a_{min}$ and is followed by waves of decreasing heights. The jump height is denoted by $a_m$ and $D$ denotes the constant velocity of the wave front. In this kind of experiments, the Froude number $F$ is defined as the ratio between the wave speed $v_0-D$ and the celerity $\sqrt{g h_0}$. Consequently, according to conservation of mass and momentum, a relation between the Froude number and the upstream and downstream water depths $h_0$ and $h_0+a_m$ can be obtained (see the B\'elanger formula Eq. (4.10) in~\cite{gavrilyuketal2016}):
\begin{equation}\label{Eq_am}
\dfrac{h_0+a_m}{h_0}=\dfrac{\sqrt{1+8F^2}-1}{2}.
\end{equation}
}
\red{
In what follows, we will compare the numerical results obtained by the extended Boussinesq model for the ``Favre waves" problem with the experimental results of Favre~\cite{Favre35} and Treske~\cite{Treske94}. On a computational domain $x\in (0,300)$, we consider a uniform initial profile defined by $\zeta(0,x)=0$ and an impact velocity $v(0,x)=v_0$ related to the relative Froude number $F$ by the formula~\cite{gavrilyuketal2016}:
$$v_0=\sqrt{g h_0} \left( F - \dfrac{1+\sqrt{1+8F^2}}{4F} \right),$$
where $g =10 \ m/s$ and $h_0=1\ m$. Reflecting boundary conditions are used on the right wall at $x=300\ m$. The solutions of the eB model are computed using the dimensional version of the eB model~\eqref{eq:eGNLW3f5} with factorized high order derivatives. To this end, we set $\eps=1$ and add the gravity term to the equations as needed. In this test, the choice of $\alpha$ is not important, thus we choose $\alpha=1$. The numerical scheme adopted in this paper for the extended Boussinesq equations can be used to solve this problem until some critical impact velocity determined in terms of the relative Froude number $F$. For higher Froude numbers ($F>\approx 1.3$), the transition from the undular bore to the bore consisting of a steep front (wave breaking) occurs. At this stage, a special treatment of the numerical scheme based on the adopted splitting strategy is required in order to handle wave breaking (see introduction of Section~\ref{NMSec}). This is not investigated in this paper and is left to a future work (see~\cite{Tissier11} for modeling of breaking waves).}

\begin{figure}[H]
		\centering
		{\includegraphics[width=0.9\textwidth]{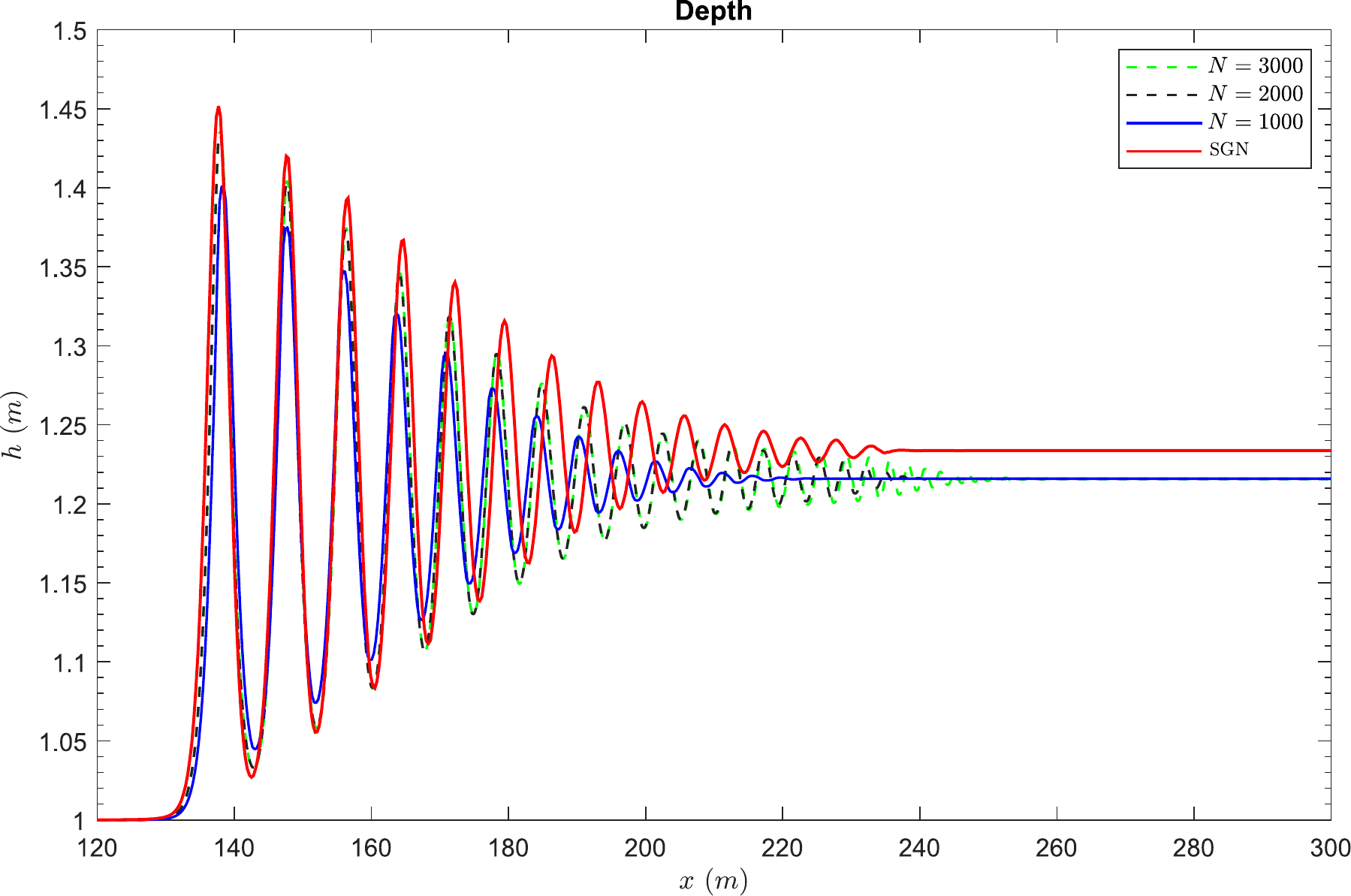}}
		\caption{Comparison of Favre waves at time $t = 54 \ s$ for the Froude number $F=1.16$. Red solid line corresponds to the numerical solution of the SGN equations obtained by the method~\cite{Favrie2017}. The results obtained with the eB model are shown for different mesh sizes: 1000 (blue solid line), 2000 (black dashed line), 3000 (green dashed line).}
		\label{compFavre}
	\end{figure}

\begin{figure}[H]
		\centering
		{\includegraphics[width=0.9\textwidth]{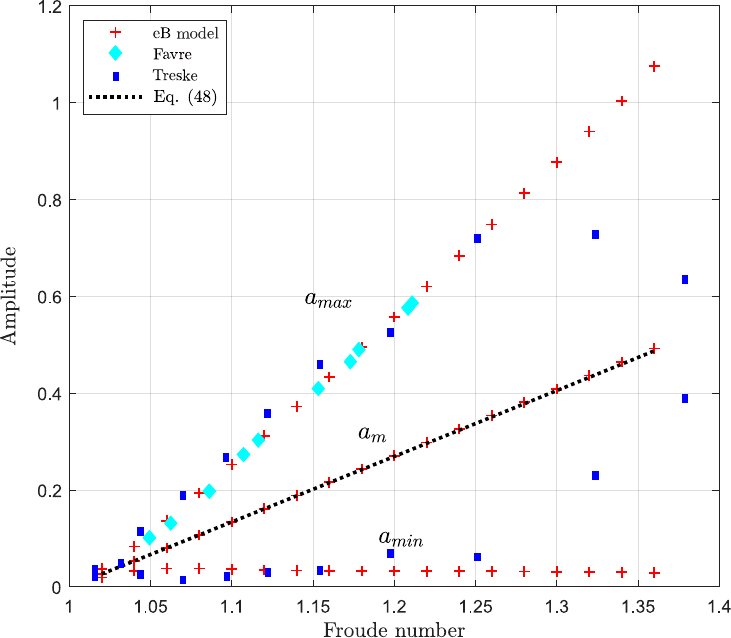}}
		\caption{Amplitude of undular bores for different Froude numbers.}
		\label{comp_a}
	\end{figure}

\red{
In Figure~\ref{compFavre}, we compare the numerical results obtained by the eB model for different mesh sizes ($N=1000, 2000$, and $3000$) with the numerical results of the Serre Green-Naghdi (SGN) equations obtained by the method in~\cite{Favrie2017}. The comparison is done at time $t=54 \ s$ for a Froude number $F=1.16$ that corresponds to the following upstream velocity $v_0=0.6490 \ m/s$. The solid red line corresponds to the numerical solution of the SGN equations obtained by the method~\cite{Favrie2017} on a $2000$ cell mesh. A good agreement is observed. Our results for $2000$ and $3000$ cells are almost superposed, thus the convergence is guaranteed. The first wave amplitude is well estimated with a finer mesh, however a small discrepancy is observed in the prediction of the jump height $a_m$ between the numerical solutions of the eB model and the SGN model. For a Froude number $F=1.16$, the corresponding jump height obtained using equation~\eqref{Eq_am} is $a_m=0.215 \ m$ which corresponds exactly to the one predicted by the eB model, see Figure~\ref{compFavre}. Thus, an accurate prediction of the jump height $a_m$ is provided by the eB model rather than the SGN equations. }
	\red{
In addition we compare the amplitudes of undular bores obtained by the eB model with experimental data of Favre~\cite{Favre35} and Treske~\cite{Treske94}. The computations are performed for different Froude numbers from the interval $F \in [1.02, 1.36]$.
The maximum $a_{max}$, the minimum $a_{min}$ amplitudes of the leading wave and the jump height $a_m$ are taken at $t =  54 \ s$ with $N=2000$. Figure~\ref{comp_a} shows that the results obtained by the eB model are in good agreement with experimental data until the wave breaking occurs corresponding to the Froude number about 1.25. After this critical value, the transition from the undular bore to the bore consisting of a steep front (breaking bore) occurs and our numerical scheme is no more valid since it does not handle wave breaking.}

	\section{Conclusion}
In this work, a numerical model is developed for a class of higher-ordered/extended Boussinesq (eB) equations describing the propagation of water-waves over a flat bottom. A reformulation with the same order of precision that avoids the calculation of high order derivatives on the surface deformation is proposed including a dispersion correction parameter $\alpha$ to be precisely chosen which improve the linear dispersion properties. Insights about linear stability and dispersion optimization are provided showing that the improvement is significant when factorizing 
every 
 high order derivative. In fact, factorizing only the fifth order derivative admits a destabilizing effect. The eB model with factorized high order derivatives provide linear dispersion characteristics which are accurate for wave numbers $k$ up to 10 which are superior to lower-order models. 

A second order time splitting scheme is then proposed relying on a combination of finite volume and finite difference methods. The hyperbolic part of the equations is discretized using a high-order finite volume WENO scheme, while the dispersive part is treated using classical high-order finite differences.

Finally several numerical computations are exhibited validating the model and the numerical methods. We began by examining the numerical solution of the eB model for the case of solitary waves and compared it to other lower-order models. \red{The eB model is found to have a better approximate solution.} 
 The propagation of a solitary wave solution with correctors is then considered allowing to study the accuracy and convergence properties of the proposed numerical scheme. In the following case, the interaction of two counter propagating solitary waves is considered as a standard nonlinear test case showing the high precision of our scheme.
The optimized higher-order eB formulation with factorization of high order derivatives provide some results which are on par with full dispersion models. In particular, the breaking of a regular heap of water  with both large and small wave numbers are studied highlighting the importance of factorizing high order derivatives in improving the frequency dispersion in high frequency regimes. \red{The dam-break problem supplemented by a comparison between the standard and extended Boussinesq models is studied}. The eB model containing high order nonlinear dispersive terms not present in the sB model generate\red{s} larger dispersive tails with higher amplitude oscillations. %
This test show\red{s} that the dispersive properties of the model are well captured thanks to the high order accuracy of the numerical scheme. \red{In the last numerical test case, a comparison with experimental data was presented for the study of ``Favre waves". The proposed scheme reproduces the measurements with a good agreement.}

   Following the steps of this work, next steps may cover the occurrence of variable topography to seriously discriminate high order models and two dimensional extension to study more real-life cases.  

\appendix
\section{Dispersion relation of the eB model~\eqref{eq:eGNLW3f5}}\label{appendix}
This appendix is devoted to the computation of the dispersion relation 
associated with the new eB model~\eqref{eq:eGNLW3f5}. We start by investigating the linear behavior of small perturbation $(\t{\zeta},\t{v})$ to a constant state solution $ (\zetabar,\vbar)$. The linear equations governing these perturbations are
	\begin{equation}\label{eq:eGNLW3f5app}
		\left\{ \begin{array}{l}
			\dsp \partial_{ t}\t{\zeta} +\vbar \eps \partial_x \t{\zeta} + \hbarr \partial_x \t{v} \ =\  0,\\ \\
			\dsp \Big(1+\eps \alpha \mathcal{T}[0] -\eps^2\alpha\mathfrak{T}\Big) \Big(\partial_{ t}  \t{v} + \eps \vbar \partial_x \t{v} +\dfrac{\alpha-1}{\alpha} \partial_x \t{\zeta} \Big) +\dfrac{1}{\alpha} \partial_x \t{\zeta}\\+\dfrac{(7-5\alpha)\eps^2}{45}  \partial_x^4 \Big((1+\eps \alpha \mathcal{T}[0] )^{-1} (\partial_x \t{\zeta}) \Big)+ \dfrac{2\eps^2}{3}\zetabar \partial_x^2 \Big((1+\eps \alpha \mathcal{T}[0] )^{-1} (\partial_x \t{\zeta}) \Big) =  \OO(\eps^3),
		\end{array} \right. \end{equation}
	where $\hbarr=1+\eps\zetabar$, $\mathcal{T}[0]w=  -\dfrac{1}{3}\partial_x^2 w$ and $\mathfrak{T} w=-\dfrac{1}{45} \partial_x^4 w$. Looking for the corresponding plane wave solutions of the form $(\zeta^0,v^0)e^{i(kx-wt)}$ with $k$
the spatial wave number and $w$ the time pulsation, one finds the dispersion relation. From the first equation of~\eqref{eq:eGNLW3f5app} one obtains:
\begin{equation}\label{eq9}-iw\t{\zeta}+\vbar(ik\eps \t{\zeta}) +\hbarr(ik\t{v})=0 \ \Rightarrow \ \t{v}=\dfrac{\t{\zeta}(w-\eps k\vbar)}{k\hbarr}.\end{equation}
The second equation of~\eqref{eq:eGNLW3f5app} becomes:
\begin{equation}\label{eq9'}\Big(1+\dfrac{\eps \alpha}{3} k^2+\dfrac{\eps^2 \alpha}{45}k^4 \Big)\Big(-iw\t{v}+\eps \vbar i k \t{v} +\dfrac{\alpha-1}{\alpha}ik\t{\zeta}\Big)+\dfrac{1}{\alpha}i k\t{\zeta}+\dfrac{(7-5\alpha)\eps^2 }{45} \dfrac{ik^5 \t{\zeta}}{1+\dfrac{\eps \alpha}{3} k^2}-\dfrac{2\eps^2 }{3}\zetabar \dfrac{ik^3 \t{\zeta}}{1+\dfrac{\eps \alpha}{3} k^2}=0 .\end{equation}
Substituting $\t{v}$ in~\eqref{eq9'} by its expression given in~\eqref{eq9} and multiplying~\eqref{eq9'} by $\dfrac{k}{i \t{\zeta}}$ yields the dispersion relation for~\eqref{eq:eGNLW3f5app}:
\begin{equation}\label{rdeGN3app}
		\dfrac{(w-\eps k\vbar)^2}{\hbarr k^2}=\dfrac{\Big(1+\dfrac{\eps(\alpha-1)k^2}{3} +\dfrac{\eps^2(\alpha-1)k^4}{45} +\dfrac{ (7-5\alpha)\eps^2k^4}{45(1+\frac{\eps \alpha}{3}k^2)}-\dfrac{ 2\eps^2k^2\zetabar}{3(1+\frac{\eps \alpha}{3}k^2)} \Big)}{\Big(1+\dfrac{\eps \alpha}{3}k^2+\dfrac{\eps^2 \alpha}{45}k^4\Big)}.
	\end{equation}
The dispersion relation obtained around some rest state solution $(\zetabar,\vbar)=(0,0)$ is the following:
\begin{equation}\label{rdeGN3apprest}
		w^2=\dfrac{k^2\Big(1+\dfrac{\eps(\alpha-1)k^2}{3} +\dfrac{\eps^2(\alpha-1)k^4}{45} +\dfrac{ (7-5\alpha)\eps^2k^4}{45(1+\frac{\eps \alpha}{3}k^2)} \Big)}{\Big(1+\dfrac{\eps \alpha}{3}k^2+\dfrac{\eps^2 \alpha}{45}k^4\Big)}.
	\end{equation}
%

\def\cprime{$'$}

%
%
%

\end{document}